\documentclass[a4paper,12pt,bibliography=totocnumbered,titlepage=false,abstracton,bookmarksnumbered=true]{scrartcl}
\usepackage{cite}
\usepackage[utf8]{inputenc}
\usepackage{amsmath}
\usepackage{amsthm}
\usepackage{amssymb}
\usepackage{amsfonts}
\usepackage{mathtools}
\usepackage{stmaryrd}
\usepackage{graphicx,color}
\usepackage[font=small, labelfont=it]{caption}
\usepackage{mathrsfs}
\usepackage[arrow,matrix,curve]{xy}
\usepackage{bbm}
\usepackage[]{hyperref}

\DeclareMathOperator{\rk}{rk}

\newcommand{\R}{{\mathbb{R}}}

\newcommand{\veps}{\varepsilon}

\newtheorem{defn}{Definition}[section]
\newtheorem{theorem}{Theorem}[section]
\newtheorem{lemma}{Lemma}[section]
\newtheorem{proposition}{Proposition}[section]

\newtheorem{corollary}{Corollary}[section]
\newtheorem{example}{Examples}[section]

\newtheorem{remark}{Remarks}[section]
\title{On manifolds with infinitely many fillable contact structures}
\author{Alexander Fauck}
\begin{document}

\maketitle

\begin{abstract}
We introduce the notion of asymptotically finitely generated contact structures, which states essentially that the Symplectic Homology in a certain degree of any filling of such contact manifolds is uniformly generated by only finitely many Reeb orbits. This property is used to generalize a famous result by Ustilovsky: We show that in a large class of manifolds (including all unit cotangent bundles and all Weinstein fillable contact manifolds with torsion first Chern class) each carries infinitely many exactly fillable contact structures. These are all different from the ones constructed recently by Lazarev.\\
Along the way, the construction of Symplectic Homology is made more general. Moreover, we give a detailed exposition of Cieliebak's Invariance Theorem for subcritical handle attaching, where we provide explicit Hamiltonians for the squeezing on the handle.
\end{abstract}

%\keywords{Exotic contact structures, exact fillings, Symplectic Homology, handle attachment}

%\ccode{Mathematics Subject Classification 2000: 53D10, 53D40}

%\tableofcontents

\section{Introduction}
\subsection{Main results}
This paper considers closed contact manifolds $(\Sigma,\xi)$ which are exactly fillable.  That means that there exists a compact manifold with boundary $V$ and a 1-form $\lambda$ on $V$, such that $\partial V = \Sigma,\; \xi=\ker \lambda|_{T\Sigma}$ and $\omega:=d\lambda$ is a symplectic form. The last condition implies that $\dim V = 2n$ is even and $\dim \Sigma= 2n{-}1$ is odd. Examples of exactly fillable contact manifolds are given by the unit cotangent bundle $S^\ast M=\big\{(x,v)\big|{x\in M}, \linebreak[1]{v\in T^\ast_x M}, ||v||_g=1\}$ of any closed Riemannian manifold $(M,g)$, whose filling is the unit disk bundle $D^\ast M=\{(x,v),||v||_g\leq 1\}$ with its tautological 1-form. Another example is provided by the standard spheres $S^{2n-1}\subset\mathbb{R}^n\times\mathbb{R}^n$ with the standard contact structure $\xi_{std}=\ker\big(x \text{d} y-y\text{d} x)|_{S^{2n-1}}$.\\
An almost contact structure on $\Sigma$ is a reduction of the structure group of $T\Sigma$ to $U(n{-}1)\times\mathbbm{1}$. Any coorientable contact structure $\xi$ on $\Sigma$ determines uniquely up to homotopy an almost contact structure by the choice of an almost complex structure on $\xi$ and a vector field $R$ transverse to $\xi$. Two contact structures on $\Sigma$ which define the same homotopy class of almost contact structures are called \emph{formally homotopic}.\\
For an exact symplectic manifold $(V,\lambda)$ as above, one can construct the Symplectic (Co)Homology $SH_\ast(V)$ resp.\ $SH^\ast(V)$, which are symplectic invariants of $V$. This (co)homology is closely related to the contact structure on $\Sigma$, as its generators can be identified with either critical points of a Morse function on $V$ or closed Reeb orbits on $\Sigma$ (see \ref{secsetup}, \ref{secsymhom} for details). The closed Reeb orbits alone generate a variant of $SH_\ast(V)$, the so called positive Symplectic Homology $SH_\ast^+(V)$. In this article, we consider only the Symplectic (Co)Homology generated by contractible Reeb orbits and we always use $\mathbb{Z}_2$-coefficients.\\
The connection to the Reeb dynamics was used to prove the result by I. Ustilovsky, \cite{Usti}, that the standard spheres $S^{4m+1}$ carry infinitely many different exactly fillable contact structures all formally homotopic to $\xi_{std}$. The proof was first carried out in \cite{Fauck1}\footnote{Using Rabinowitz-Floer homology, which can be thought of as a relative version of $SH$.} and later by Kwon and van Koert, \cite{KwKo}, and Gutt, \cite{Gutt2}, using $S^1$-equivariant Symplectic Homology. These so called exotic contact structures on $S^{4m+1}$ can all be described as Brieskorn manifolds $\Sigma_p, p\equiv \pm 1\bmod 8$ (see Sec.\ \ref{secBries}).\\
One can construct new exactly fillable contact manifolds from given ones by attaching a symplectic handle to the filling $V$ along the boundary $\Sigma$ (see Sec.\ \ref{secsur}). On $(\Sigma,\xi)$, this has the effect of contact surgery, in particular contact connected sums. The connected sum of any exactly fillable contact manifold with a standard sphere carrying an exotic contact structure is particularly interesting, as the underlying differentiable manifold stays unchanged. This construction was first explored by I. Ustilovsky in his thesis (cf.\ \cite{Usti2}). See also \cite{Esp} for more recent results. We show in this paper that for a subclass of exactly fillable contact manifolds (which is considerably larger compared to previous results) one obtains infinitely many different contact structures via connected sums.\\
Let us briefly describe this subclass (see \ref{secAsmFini} for the exact definition): A contact form $\alpha$ for $\xi$ is a 1-form such that $\xi=\ker\alpha$. The Reeb vector field $R$ of $\alpha$ is the unique vector field on $\Sigma$ such that $\alpha(R)=1$ and $d\alpha(R,\cdot)=0$. A closed $R$-orbit is transversely non-degenerate, if it satisfies a Morse-Bott assumption (see Sec. \ref{secsetup}). In particular, such orbits form a smooth submanifold $\mathcal{N}$ of the free loop space of $\Sigma$. A generator $\gamma$ of the Floer chain complex defining $SH_\ast(V)$ is a critical point of a Morse function $h$ on $\mathcal{N}$. Its degree $\mu(\gamma)$ is given by the Morse-Bott index (see \ref{secConZeh}):
\[\mu(\gamma)=\mu_{CZ}(\gamma)+\mu_{Morse}(\gamma)-{\textstyle \frac{1}{2}\dim_\gamma\mathcal{N}+\frac{1}{2}},\]
where $\mu_{CZ}$ and $\mu_{Morse}$ are the Conley-Zehnder and Morse index\footnote{$\mu_{Morse}(\gamma)-{\textstyle \frac{1}{2}}\dim_\gamma\mathcal{N}$ is also known as the signature index of $\gamma$.} of $\gamma$ and $\dim_\gamma\mathcal{N}$ is the local dimension of $\mathcal{N}$ near $\gamma$.\\
Fix a reference contact form $\alpha$. We call a contact structure \emph{asymptotically finitely generated (a.f.g.) in degree $k$ with bound $b_k(\xi)$} if there exist sequences of functions $f_l:\Sigma\rightarrow \mathbb{R}$ and numbers $(\mathfrak{a}_l)\subset\mathbb{R}$ with $\log \mathfrak{a}_l-f_l\leq \log \mathfrak{a}_{l+1}-f_{l+1}$ and $\lim_{l\rightarrow\infty}(\log\mathfrak{a}_l-f_l)=\infty$, such that for each contact form $\alpha_l:=e^{f_l}{\cdot}\alpha$ all contractible $R_l$-orbits of length less than $\mathfrak{a}_l$ are transversely non-degenerate and that for some Morse function $h$ on these orbits at most $b_k(\xi)$ have index $k$. Note that $\alpha_l=\alpha\,\forall l$, i.e.\ $f_l\equiv 0$ is permitted.\pagebreak[3]
\begin{remark}~
\begin{itemize} 
 \item The notion of asymptotically finitely generated contact structures generalizes the related notions of dynamical convex or index positive contact structures (cf.\ \cite{CieOan}, 9.5), asymptotically dynamically convex
 contact structures (cf.\ \cite{Lazarev}, 3.6) and contact structures with convenient dynamics (cf.\ \cite{KwKo}, 5.14). In particular, there exist exactly fillable contact manifolds which are a.f.g.\ in almost all degrees, but satisfy none of the other three notions (examples are the Brieskorn manifolds $\Sigma(n{+}1,...,n{+}1)\subset\mathbb{C}^{n+1}$, cf.\ \cite{FauckThesis}, 7.2).
 \item There exist exactly fillable contact manifolds $(\Sigma,\xi)$ which are not a.f.g.. In fact, McLean showed in \cite{McLean} that if $\Sigma$ admits a Stein fillable contact structure, then there exists on $\Sigma$ a contact structure $\xi$ with a filling $(V,\lambda)$, such that $\rk SH_k(V;\mathbb{Q})=\infty\;\forall\, k$. This shows by Prop.\ \ref{PropAsympfiniteGene}, that $\xi$ cannot be a.f.g.\ in any degree $k$. However, we have the following three general existence results.
\end{itemize}                                                                                                                                                                                                                                                                                                                                                                                                                                                                                                                                                                                                                                                                                                                                                                                                                                                                                                                                                                                                                                                                                                                                                                                                                                                                                                                                                                                                                                                                                                                                                                                                                                                                                                                                                                                                                                                                                                                                                                                                                                                                                                                                                                                                                                                                                                                \end{remark}
\begin{proposition}\label{PropafgExists1}
 Let $M$ be any closed manifold, $\dim M\geq 3$, and let $(S^\ast M, \xi_{std})$ be its unit cotangent bundle (provided by a Riemannian metric $g$) with its standard contact structure. Then $\xi_{std}$ is a.f.g.\  in every degree $k<0$ with bound $b_k(\xi_{std})=0$.
\end{proposition}
\begin{proof}
 This is closely related to the fact, that the (full)\footnote{The full SH is generated by all closed Reeb orbits. As mentioned above, we consider elsewhere in this article only the Symplectic Homology generated by the contractible Reeb orbits.} Symplectic Homology $SH_\ast\mspace{-2mu}(D^\ast\mspace{-4mu} M)$ of the unit disk bundle of $M$ is isomorphic to the singular homology (with $\mathbb{Z}_2$-coefficients) of the free loop space $\Lambda (M)$ of $M$ (see \cite{AbbSchw}). More explicitly, closed Reeb orbits $\gamma$ on $S^\ast M$ can be understood as closed geodesics on $M$. Then, $\gamma$ can be seen as a critical point of a Lagrangian action functional $\mathcal{L}$ on $\Lambda(M)$. Its Morse index $\mu_{Morse}(\gamma)$ is the dimension of a maximal subspace of $T_\gamma\Lambda(M)$ on which the Hessian of $\mathcal{L}$ is negative definite. As such, $\mu_{Morse}(\gamma)\geq 0$. Moreover, it was shown by Duistermaat in \cite{Duist} (in the non-orientable case by Weber in \cite{Weber}, Thm.\ 1.2) that $\mu_{Morse}(\gamma)=\mu_{CZ}(\gamma)$, so that $\mu_{CZ}(\gamma)\geq 0$ for all closed Reeb orbits $\gamma$ (in the non-orientable case we have to assume that $\gamma$ is contractible).\\
 If we choose a contact form $\alpha$ (resp.\ a Riemannian metric $g$) for $\xi_{std}$ such that all closed Reeb orbits are transversely non-degenerate and set $\alpha_l=\alpha\,\forall\, l$, then no closed Reeb orbit of $\alpha_l$ has index $k<0$. Hence, $\xi_{std}$ is a.f.g.\ in all degrees $k<0$ with bound $b_k(\xi_{std})=0$.
\end{proof}
\begin{proposition}[Lazarev, \cite{Lazarev}, Cor.\ 4.1]
 If $(\Sigma,\xi), \dim \Sigma\geq 3$, has a flexible Weinstein filling, then $(\Sigma, \xi)$ is asymptotically dynamically convex and in particular a.f.g.\ in degree $k\leq 3-n$ (Lazarev uses the reduced index $\mu_{CZ}+n{-}3$). Moreover, any Weinstein fillable contact manifold carries a flexible Weinstein fillable contact structure (cf.\ \cite{Lazarev}, proof of Thm.\ 1.9).
\end{proposition}
\begin{proposition}\label{PropafgExists}
 Let $(\Sigma,\xi)$ be a closed exactly fillable contact manifold such that there exists a contact form $\alpha$ for $\xi$ whose Reeb flow is periodic. Then $\xi$ is a.f.g.\ in every degree $k\in\mathbb{Z}\setminus[-(2n{-}1){+}1,2n{-}1]$. If the mean index $\Delta(\Sigma)$ of the principal Reeb orbit is non-zero, then $\xi$ is a.f.g.\ in every degree (see \ref{appendixA} for a proof).
\end{proposition}
Examples for exactly fillable contact manifolds which admit periodic Reeb flows are given by the standard spheres $S^{2n-1}$, all Brieskorn manifolds (see \ref{secGenBrie}), prequantization bundles (Boothby-Wang bundles) and unit cotangent bundles of manifolds admitting a periodic geodesic flow such as $S^n,\,\mathbb{C}P^n$ or all lens spaces.\\
Starting from these initial examples of a.f.g.\ contact structures, one obtains more examples by taking iteratively contact connected sums, as shown in the following proposition proved in \ref{secHandleAt}. In particular all boundaries of subcritical Weinstein domains are a.f.g.\ (see also \cite{Yau}).
\begin{proposition}\label{propcofinalHamilt}
 Let $(\Sigma,\xi)$ be a $(2n{-}1)$-dimensional closed contact manifold such that the contact structure $\xi$ is a.f.g.\ in a degree $j$ with bound $b_j(\xi)$. Let $\alpha$ be a contact form for $\xi$ and let $f_l:\Sigma\rightarrow\mathbb{R}, \alpha_l=e^{f_l}{\cdot}\alpha$ and $(\mathfrak{a}_l)\subset\mathbb{R}$ be sequences showing that $\xi$ is a.f.g.. Assume that $(\widehat{\Sigma},\hat{\xi})$ is obtained from $(\Sigma,\xi)$ by a $(k{-}1)$-surgery, $1\leq k\leq n$, $k\neq 2$, along an embedded isotropic sphere $S\subset\Sigma, \dim S=k{-}1$, such that any $\alpha_l$-Reeb chord\footnote{A Reeb chord of $S$ is a Reeb trajectory of positive length with starting and end point on $S$.} of $S$ is longer than $\mathfrak{a}_l$. Then $\hat{\xi}$ also a.f.g.\ in degree $j$ with bound 
 \[b_j(\hat{\xi})=\begin{cases}b_j(\xi)+1&\text{if }j=(n{-}k)(2N{-}1)+\left\lbrace\begin{smallmatrix}1\phantom{(n-k)}\\2(n-k)\end{smallmatrix}\right.,\;N\in\mathbb{N},\\
 b_j(\xi)&\text{otherwise.}\end{cases}\] 
 Moreover, the contact forms $(\beta_l)$, that show that $\hat{\xi}$ is a.f.g., can be chosen to agree with $(\alpha_l)$ outside an arbitrary small neighborhood of $S$.
\end{proposition}
\begin{remark}~
\begin{itemize}
 \item The condition $k\neq 2$ is needed, as 1-surgery can change the first Chern class. This can make a $\mathbb{Z}$-grading via the Conley-Zehnder index problematic (cf.\ Rem.\ \ref{ref3}). Moreover, 1-surgery can change contractible loops into non-contractible and vice versa. As we use the a.f.g.\ condition in this paper only for contractible orbits, 1-surgery is excluded for simplicity. However, for subcritical Weinstein manifolds, this restriction can be dropped as the number of all closed Reeb orbits of $\alpha_l$ of length less than $\mathfrak{a}_l$ can be uniformly bounded.
 \item The Reeb chord condition cannot be satisfied if $\alpha_l=\alpha\;\forall l$ is a fixed contact form with periodic Reeb flow. However, we show in \ref{appendixA}, that in this situation, we can perturb $\alpha$ to obtain a sequence $\alpha_l=e^{f_l}{\cdot}\alpha$ which also shows that $\xi$ is a.f.g.\ (maybe with different bound $b_j(\xi)$) and no $\alpha_l$ has a periodic Reeb flow. 
 \item If $S$ is subcritical, i.e.\ $\dim S=k{-}1<n{-}1$, then the Reeb chord condition can be achieved by a small perturbation of $S$ using a transversality argument (see \ref{appendixB}). In \ref{appendixA}, we give in the case that $(\Sigma,\xi)$ admits a periodic Reeb flow an explicit subcritical isotropic sphere $S$ which satisfies the Reeb chord condition. 
 \item If $S= S^0$ (two points), i.e.\ for $k=1$, the Reeb chord condition can always be achieved: It suffices to choose for $S$ two points that do not lie on any closed Reeb orbit and not on the same Reeb trajectory for any $\alpha_l$. In particular Prop.\ \ref{propcofinalHamilt} can always be applied to connected sums.
\end{itemize}
\end{remark}
 The notion of a.f.g.\ contact structures allows us to state the following generalization of Ustilovsky's result.
\begin{theorem}\label{theoinfinitecontactstr} 
 Suppose that $(\Sigma,\xi)$ is a closed contact manifold, $\dim \Sigma\geq 5$ and has an exact filling $(V,\lambda)$ such that for the inclusion $i:\Sigma\rightarrow V$ holds $i_\ast: \pi_1(\Sigma)\rightarrow\pi_1(V)$ is injective and that the Conley-Zehnder index is well-defined for all contractible Reeb orbits on $V$.\\
 If $\xi$ is asymptotically finitely generated in at least one degree $k\in\mathbb{Z}$, then $\Sigma$ carries infinitely many pairwise non-contactomorphic exactly fillable contact structures. If $\dim \Sigma = 4m{+}1$, then all these contact structures can be chosen to be formally homotopic to $\xi$.
\end{theorem}\pagebreak[1]
\begin{remark}~\label{ref3}
\begin{itemize}
 \item The Conley-Zehnder index $\mu_{CZ}$ is well-defined on contractible Reeb orbits, if the integral of the first Chern class $c_1(TV)$ vanishes on spheres.
 \item The conditions on $\pi_1(\Sigma)$ and $\mu_{CZ}$ on $V$ ensure that the Conley-Zehnder index is well-defined on $\Sigma$. Note that both conditions are invariant under attaching a symplectic $k$-handle ($k{-}1$-surgery on $\Sigma$) if $k\neq 2$ (see \cite{FauckThesis}, Lem.\ 66 and 67).
 \item It follows from the proof of Thm.\ \ref{theoinfinitecontactstr} that all the different contact structures $\xi'$ have $SH(V,\xi')\neq 0$ for some exact filling $V$ of $(\Sigma,\xi')$. Hence, $(\Sigma,\xi')$ is not flexible Weinstein fillable (see \cite{Zhou}, Cor.\ 1.3). It follows that all our fillable contact structures $\xi'$ on $\Sigma$ are different from the infinitely many flexibly fillable contact structures on $\Sigma$ constructed by Lazarev (cf.\ \cite{Lazarev}).
 \item For $\dim \Sigma = 4m{+}3$, it should in principle also be possible to construct infinitely many contact structures on $\Sigma$ which are all formally homotopic to $\xi$. The difficulty is to find as in Prop.\ \ref{propTwoConStr} for any $k$ a contact structure $\xi_k$ on $S^{4m+3}$ that is a.f.g.\ in degree $k$, has an exact filling $(V_k,\lambda_k)$ with $\rk SH^+_k(V_k)>0$ and is formally homotopic to $\xi_{std}$. Consult \cite{Ueb} by P. Uebele, where some results in this direction are shown. In particular, $S^7, S^{11}, S^{15}$ all carry infinitely many exotic contact structures all formally homotopic to $\xi_{std}$.
 \item Note that in general the connected sum of $(\Sigma,\xi)$ with an exotic contact sphere alters the homotopy class of almost contact structure.
 \item If $(\Sigma,\xi)=(S^\ast M, \xi_{std}), \dim M\geq 3,$ and $V=D^\ast M$, then $c_1(TD^\ast M)=0$\footnote{Note that only the real first Chern class vanishes. The integer first Chern class might be non-zero if $M$ is not orientable} (see Weber \cite{WeberDis}, B.1.7). Moreover, if $(V,\lambda)$ is any Weinstein filling of $(\Sigma,\xi)$ and $\dim V=2n\geq 6$, then $i_\ast:\pi_1(\Sigma)\rightarrow \pi_1(V)$ is always injective (cf.\ \cite{Zhou}, Prop.\ 2.3): $V$ can be built from $\Sigma$ by attaching handles with index greater then or equal to $n$. This implies $\pi_2(V,\Sigma)=0$, as $n\geq 3$. The injectivity of $i_\ast$ follows then from the homotopy long exact sequence
 \[...\rightarrow\pi_2(V,\Sigma)\rightarrow \pi_1(\Sigma)\rightarrow \pi_1(V)\rightarrow...\]
\end{itemize}
\end{remark}
In his thesis, \cite{Usti2} Thm. II , I. Ustilovsky also considered arbitrary connected sums of certain Brieskorn manifolds to construct exotic fillable contact structures on $S^{4m+1}$. Using similar arguments as for the proof of Thm.\ \ref{theoinfinitecontactstr}, we can give a new proof (not using Contact Homology) that the following subclass of the above mentioned exotic contact structures on $S^{4m+1}$ are all different:
\begin{theorem}\label{differentcontacttheo}~\\
 Let $\Sigma_p\cong S^{4m+1}$ denote the Brieskorn/Ustilovsky spheres (see Sec.\ \ref{secBries}). Let $\xi_p$ denote its (exotic) contact structure and let $j\cdot \xi_p$ denote the contact structure on the $j$-fold connected sum $j\cdot \Sigma_p$ of $(\Sigma_p,\xi_p)$ with itself. If $m\geq 2$ and $i,j,p,q\in\mathbb{N}$ are such that $p,q\equiv \pm {1}\bmod 8$ and $(j,p)\neq (i,q)$, then $j\cdot \xi_p$ and $i\cdot \xi_q$ are non-contactomorphic contact structures on $S^{4m+1}$.
\end{theorem}
Note that $i{=}j{=}1$ gives exactly Ustilovsky's original result from \cite{Usti}. As a corollary, we obtain that the mean Euler characteristic $\chi_m$ of the $S^1$-equivariant Symplectic Homology (see \cite{KoertDis}) is not a complete contact invariant. For dimension 5, this was already established in \cite{KwKo}, 5.10, while the following corollary extends this result to any dimension of the form $4m{+}1$.
\begin{corollary}\label{CorMeanEul}
 For any integers $m,k\geq 2$, there exist $c(k)\in\mathbb{Q}$ and \emph{$k$} different contact structures $\xi_k$ on $S^{4m+1}$ with the same mean Euler characteristic $\chi_m\big(S^{4m+1},\xi_k\big)=c(k)$.
\end{corollary}
The proof of Thm.\ \ref{theoinfinitecontactstr} relies on the following three results for a.f.g.\ contact manifolds $(\Sigma,\xi)$ and Symplectic Homology.
\begin{proposition}\label{PropAsympfiniteGene}
 Let $(\Sigma,\xi)$ be a closed contact manifold. If $\xi$ is a.f.g.\ in degree $k$ with bound $b_k(\xi)$, then it holds for any exact filling $(V,\lambda)$ of $(\Sigma,\xi)$ with $i_\ast:\pi_1(\Sigma)\rightarrow\pi_1(V)$ injective and $\mu_{CZ}$ well-defined, that
 \[\rk SH_k^+(V)\leq b_k(\xi)\quad\text{ and }\quad \rk SH_k(V)\leq b_k(\xi) + \rk H_{n-k}(V,\partial V).\]
\end{proposition}
\begin{proposition}\label{propTwoConStr}
 For $n\geq 3$ there exists for every degree $k\in\mathbb{Z}$ a contact manifold $(\Sigma_k,\xi_k)$ with an exact filling $(V_k,\lambda_k)$ such that $\Sigma_k$ is homeomorphic to $S^{2n-1}$, $\xi_k$ is a.f.g.\ in degree $k$ (in fact in every degree) and
 \begin{equation}\label{estV_K}
  \rk SH_k^+(V_k)=\begin{cases}2&\text{for }k\geq n+1\\1&\text{for }k\leq n\end{cases}.
 \end{equation}
\end{proposition}
\begin{theorem}[Invariance of $SH$ under subcritical handle attachment, Cieliebak, \cite{Cie}]\label{theoinvsur}~\\
 Let $W$ and $V$ be compact $2n$-dimensional symplectic manifolds with positive contact boundaries and assume that $\mu_{CZ}$ is well-defined on $W$. If $V$ is obtained from $W$ by attaching to $\partial W$ a subcritical symplectic handle $\mathcal{H},\,k<n$, then it holds that \[SH_\ast(V)\cong SH_\ast(W)\qquad\text{ and }\qquad SH^\ast(V)\cong SH^\ast(W)\qquad \forall\, \ast\in\mathbb{Z}.\]
\end{theorem}
The positive Symplectic Homology $SH^+(W)$ is not invariant under subcritical handle attachment. However, its changes can be estimated via Thm.\ \ref{theoinvsur}. In particular for the contact connected sums, we have the following result.
\begin{corollary}\label{corinvsur}
 Let $W$ and $V$ be as in Thm.\ \ref{theoinvsur} and assume that $V$ is obtained from $W$ by boundary connected sum, i.e.\ a 1-handle $\mathcal{H}$ is attached to two different path-connected components of $W$. Then
 \begin{itemize}
  \item $SH^+_\ast(V)\cong SH^+_\ast(W)$ for $\ast\neq n, n{+}1$  
  \item either $\rk SH^+_n(V)=\rk SH^+_n(W)\,{+}1\,$ and $\;\rk SH^+_{n+1}(V)=\rk SH^+_{n+1}(W)$\\
  or $\phantom{the}\rk SH^+_n(V)=\rk SH^+_n(W)\quad$ and $\;\rk SH^+_{n+1}(V)=\rk SH^+_{n+1}(W)\,{-}1$.
 \end{itemize}
\end{corollary}
\begin{proof}[\textnormal{\textbf{Proof of Theorem \ref{theoinfinitecontactstr}}}]~\\
Assume that $\xi$ is a.f.g.\ in degree $k$ with bound $b_k(\xi)$ and let $(\Sigma_k,\xi_k)$ be the contact manifold with exact filling $(V_k,\lambda_k)$ from Prop.\ \ref{propTwoConStr} for the degree $k$, i.e.\ $\Sigma_k$ is homeomorphic to $S^{2n-1}$, $\xi_k$ is a.f.g.\ in degree $k$ with bound $b_k(\xi_k)$ and $\rk SH_k^+(V_k)$ satisfies the estimate (\ref{estV_K}).\\
We denote by $\big(\Sigma\# \Sigma_k, \xi\#\xi_k\big)$ the contact connected sum of $(\Sigma,\xi)$ and $(\Sigma_k,\xi_k)$ and by ${\big(N{\cdot} \Sigma_k,N{\cdot} \xi_k\big)}$ the $N$-fold contact connected sum of $(\Sigma_k,\xi_k)$ with itself. As the differentiable structures on $S^{2n-1}$ form a finite group with respect to the connected sum operation with the standard structure as neutral element (see \cite{KerMil}), there exists an $N_0$, such that $N_0\cdot\Sigma_k$ is diffeomorphic to $S^{2n-1}$.\\
The homotopy classes of almost contact structures on the standard sphere $S^{2n-1}$ are classified by $\pi_{2n-1}\big(SO(2n)/U(n)\big)$ with 0 corresponding to the standard structure, as shown by Morita, \cite{Mor}. For $S^{2n-1}=S^{4m+1}$, i.e.\ $n$ odd, these homotopy groups are finite. Geiges, \cite{geiges1} Sec.\ 5.2, showed that taking contact connected sums corresponds to the addition in $\pi_{2n-1}\big(SO(2n)/U(n)\big)$. Thus, there exists in dimensions $4m{+}1$ an $N_0$, such that $N_0\cdot\Sigma_k$ is diffeomorphic to $S^{2n-1}=S^{4m+1}$ and ${N_0\cdot\xi_k}$ is formally homotopic to $\xi_{std}$. As contact structures are locally trivial by Darboux's Theorem, one finds that the contact connected sum of $(\Sigma,\xi)$ and $(N_0\cdot\Sigma_k,N_0\cdot\xi_k)$ has the same homotopy class of almost contact structures as $\xi$, if $\dim \Sigma=4m{+}1$.\pagebreak[3]\\
Henceforth, we assume that $N_0|N$. Taking the boundary connected sum of $(V,\Sigma)$ with $N$ copies of $(V_k,\Sigma_k)$ yields by Cor.\ \ref{corinvsur}:
\[\rk SH_k^+\big(V \,\#\,N{\cdot} V_k\big)= \rk SH_k^+\big(V\big)+N\cdot \rk SH_k^+\big(V_k\big)+c,\]
where $c=0$ if $k\neq n,n{+}1$, $0\leq c\leq N$ if $k=n$ and $-N\leq c\leq 0$ if $k=n{+}1$. Using \ref{estV_K}, we hence have the estimate:
\[\rk SH_k^+\big(V \,\#\,N{\cdot} V_k\big)\geq \rk SH_k^+\big(V\big)+N\geq N. \tag{$\ast$}\]
As $N\cdot\Sigma_k$ is diffeomorphic to $S^{2n-1}$, we have $\Sigma\,\#\big(N{\cdot}\Sigma_k\big)\cong\Sigma$ and $\big(\Sigma,\xi \,\#\,N{\cdot} \xi_k)\big)$ is exactly fillable by $V \,\#\,N{\cdot} V_k$. Prop.\ \ref{propcofinalHamilt} implies that $\big(\Sigma, \xi \,\#\,N{\cdot} \xi_k\big)$ is also a.f.g.\ in degree $k$ with bound
\[b_k\big(\xi \,\#\,N{\cdot} \xi_k\big)\leq b_k(\xi)+N\cdot b_k(\xi_k)+N.\]
This estimate with Prop.\ \ref{PropAsympfiniteGene} implies for any exact filling $(W,\lambda_W)$ of $\big(\Sigma,\xi\#N{\cdot}\xi_k)\big)$ that
\[\rk SH_k^+(W)\leq b_k\big(\xi\,\#\,N{\cdot} \xi_k\big)\leq b_k(\xi)+N\big(b_k(\xi_k){+}1\big).\tag{$\ast\ast$}\]
$(\ast)$ and $(\ast\ast)$ imply for any $N,M\in N_0{\cdot}\mathbb{N}$ with $b_k(\xi){+}M\big(b_k(\xi_k){+}1\big){<}\,N$ that $\xi \,\#\,M{\cdot} \xi_k$ and $\xi \,\#\,N{\cdot}\xi_k$ are non-contactomorphic contact structures on $\Sigma$, as the latter has a filling which cannot be a filling of the first.\\
Set $b:=\max\big\{b_k(\xi),b_k(\xi_k)\big\}$ and  $N_l:=N_0(b{+}2)^l, l\in\mathbb{N},$ and calculate
   \[b_k(\xi)+N_l\big(b_k(\xi_k)+1\big)\leq b+N_0(b+2)^l(b+1)<N_0(b+2)^{l+1}=N_{l+1}.\]
   Thus $\xi \,\#\,N_l{\cdot}\xi_k$ and $\xi \,\#\,N_m{\cdot}\xi_k$ are not contactomorphic if $l\neq m$, which provides infinitely many different exactly fillable contact structures on $\Sigma$.\medskip\\
   \underline{Bonus} : All contact structures $\xi \,\#\,N{\cdot}\xi_k$ on $\Sigma$ are different for different values of \linebreak ${N\in N_0\big(\mathbb{N}{\cup}\{0\}\big)}$. We argue by contradiction and assume that there exists $M<N$ with $\xi\,\#\,M{\cdot}\xi_k$ contactomorphic to $\xi\,\#\,N{\cdot}\xi_k$. Then we have for all $K\in\mathbb{N}$:
   \[\xi\,\#\,(N{+}K)\xi_k=\xi\,\#\, N\xi_k\,\#\,K\xi_k \cong \xi\,\#\, M\xi_k\,\#\,K\xi_k=\xi\,\#\,(M{+}K)\xi_k.\]
   In particular for $K_0=N{-}M$, we obtain
   \[\xi\,\#\,(M{+}2K_0)\xi_k=\xi\,\#\,(N{+}K_0)\xi_k\cong\xi\,\#\,(M{+}K_0)\xi_k=\xi\,\#\,N\xi_k\cong\xi\,\#\,M\xi_k.\]
   Repeating this argument, we find that $\xi\,\#\,(M{+}lK_0)\xi_k\cong\xi\,\#\,M\xi_k$ for all $l\in\mathbb{N}$. However, for $l$ sufficiently large, we already know by $(\ast)$ and $(\ast\ast)$ that these two contact structures are different, thus providing a contradiction.
\end{proof}
\subsection{Structure of the paper}
The paper is organized in three parts. In the first, we describe the construction of the Symplectic Homology $SH(V)$ of a Liouville domain $(V,\lambda)$ and discuss many of its properties, i.e.\ its positive part (\ref{secPosSym}), the canonical long exact sequence (\ref{sectrunc} (\ref{eqLongExSeq})), the transfer maps (\ref{sectransfer}) and its grading via the Conley-Zehnder index (\ref{secConZeh}).\\
Our description is more extensive than some readers might think necessary, as we define $SH(V)$ with a Morse-Bott setting on the contact boundary $\partial V$ -- an approach which is not yet standard in the literature. Moreover, our construction of $SH(V)$ and $SH^+(V)$ allows the use of a sequence of different (!) contact boundaries in the completion $\widehat{V}$ of $V$.\linebreak[1] This makes the definition of $SH(V)$ more flexible and leads naturally to the definition of asymptotically finitely generated contact structures (\ref{secAsmFini}). Another advantage of this more general approach is that $SH(V)$ and $SH^+(V)$ become almost by definition invariant under Liouville isomorphisms (see Prop.\ \ref{propSHinv}). This fact is of course well-known (see \cite{Gutt2}, Thm.\ 4.18), but requires usually a rather long proof via the transfer maps.\\
The use of different contact boundaries forced us to generalize the maximum principle (\ref{secNoEsc}). In \cite{Sei}, 3.21, Seidel already gave a maximum principle for varying contact boundaries -- yet his version requires that the slopes of the Hamiltonians involved increase exponentially when going from one boundary to another, whereas we only require that the Hamiltonians just increase at infinity. Our version of this confinement tool is based on the No-Escape lemma first presented by Abouzaid-Seidel, \cite{AbouSei} 7.2, and adapted to the setup of Symplectic Homology by Ritter, \cite{Ritter} Lem.\ 19.3.\\
The second part of the paper describes the attachment of subcritical handles and its effect on Symplectic Homology and a.f.g.\ contact structures. This section owes much to \cite{Cie} by K.\ Cieliebak. Unfortunately, this highly influential and ground-breaking paper contains
two gaps:
\begin{itemize}
 \item For the proof of the Invariance Theorem (see the second to last page in \cite{Cie}), a definition of $SH$ with varying contact boundaries is needed, as the attached handle is made iteratively thinner.
 \item The construction of the Hamiltonians involved is only sketched and it is claimed that these Hamiltonians can be chosen to have on the handle only one (!) 1-periodic orbit. Whether this last property can be achieved at all is heavily doubted by the author (see Rem.\ \ref{dis1} in \ref{ExPsi}).
\end{itemize}
As we had to describe subcritical handle attachment anyway in order to show that the a.f.g.\ property is preserved under this procedure, we seized the opportunity to give in Section \ref{secsur} a clarified proof of Cieliebak's Invariance Theorem.\\
The third part presents Brieskorn manifolds. First, we review some general results, especially on the indices of closed Reeb orbits on these manifolds. Then we give some explicit examples of Brieskorn spheres (proof of Prop.\ \ref{propTwoConStr}) and show Thm.\ \ref{differentcontacttheo}.\\
Appendices on contact manifolds with periodic Reeb flow and on avoiding Reeb chords conclude the paper. \ref{appendixB} was friendly provided by Dingyu Yang.\pagebreak[1]

\section{Symplectic Homology and cohomology}
Symplectic Homology for Liouville domains together with the transfer maps was introduced by Viterbo, \cite{Vit}. In this section we relax the standard definition to include families of Hamiltonians that are linear each with respect to a different contact boundary. For that, we generalize the No-Escape Lemma by Abouzaid/Seidel and Ritter and we discuss many of the additional features of Symplectic Homology in this generalized approach. This will allow us to estimate $SH(V)$ in the presence of the a.f.g.\ property (Prop.\ \ref{PropAsympfiniteGene}). Since Symplectic Homology is a well established theory, many details will be omitted. However, as the use of Morse-Bott techniques is less developed, a longer exposition of this approach is included. For background, we refer to the excellent papers \cite{Abou}, \cite{Sei} for Symplectic Homology in general and \cite{BourOan1} for the Morse-Bott setting. 

\subsection{Setup}\label{secsetup}
Let $(V,\omega)$ be an $2n$-dimensional compact symplectic manifold with boundary ${\partial V = \Sigma}$ such that $\omega \,{=}\,d\lambda$ is exact. The 1-form $\lambda$ defines the Liouville vector field $Y$ by ${\omega(Y,\cdot)=\lambda}$. The boundary $\Sigma$ is called a positive or negative contact boundary if $Y$ points out or into $V$ along $\Sigma$. If $\Sigma$ is a positive contact boundary, then $(V,\lambda)$ is called a Liouville domain.\\
Any hypersurface $\Sigma$ in $V$ transverse to $Y$ is a contact manifold with contact form $\alpha:=\lambda|_{T\Sigma}$. We write $\xi:=\ker \alpha$ for the contact structure and $R:=R_\alpha$ for the Reeb vector field defined by $d\alpha(R,\cdot)= 0$ and $\alpha(R)=1$. We denote by the spectrum $spec(\alpha)$ the set of periods of closed orbits of $R$. We say that $\alpha$ is \textit{transversely non-degenerate} if it satisfies the Morse-Bott assumption:
\begin{equation}\label{CondMB}\begin{aligned}
 &\textit{The set $\mathcal{N}^\eta\subset\Sigma$ formed by the $\eta$-periodic Reeb orbits is a submanifold}\\
 &\textit{for all $\eta\in\mathbb{R}$ and $T_p\,\mathcal{N}^\eta = \ker \left(D_p\phi^\eta - \mathbbm{1}\right)$ holds for all $p\in\mathcal{N}^\eta$.}
 \end{aligned}\tag{MB}
\end{equation}
A closed Reeb orbit $x$ is called transversely non-degenerate if (\ref{CondMB}) holds locally.\\
A symplectization of a contact manifold $(\Sigma,\alpha)$ is the symplectic manifold  $\big(I{\times}\Sigma, \omega:=d(e^r\alpha)\big)$, where $r\in I$ is an interval. For $I{=}\mathbb{R}, I{=}[0,\infty)$ or $I{=}(-\infty,0]$, one calls $(I{\times}\Sigma,\omega)$ the whole/positive/negative symplectization of $\Sigma$. If $\beta$ is a different contact form on $\Sigma$ defining the same contact structure and the same orientation, we find a function $f:\Sigma\rightarrow\mathbb{R}$ such that $\beta_p=e^{f(p)}{\cdot} \alpha_p \;\forall \,p{\in}\Sigma$. To such a $\beta$, we associate a hypersurface $\Sigma_\beta$ in the whole symplectization of $\Sigma$ by:
\[\Sigma_\beta :=\big\{(f(p),p)\,\big|\,p\in\Sigma\big\}.\]
Note that $(e^r\alpha)|_{\Sigma_\beta}=\beta$ and that $\Sigma_\alpha=\Sigma$. Moreover $\Sigma_\alpha$ and $\Sigma_\beta$ are naturally diffeomorphic for any two contact forms with $\ker \alpha=\ker \beta$.\\
The flow $\varphi_Y$ of $Y$ on a Liouville domain $(V,\lambda)$ with contact boundary $(\Sigma,\alpha)$ provides an exact identification $\big((-\infty,0]{\times}\Sigma, e^r\alpha \big)\rightarrow (U,\lambda)$ between the negative symplectization of $(\Sigma,\alpha)$ and a collar neighborhood $U$ of $\Sigma$ in $V$ via the map $(r,p)\mapsto \varphi_Y^r(p)$. This allows us to define the completion $(\widehat{V},\widehat{\lambda})$ of $(V,\lambda)$ by
\[\widehat{V}:=V\cup_{\varphi_Y}\big((-\delta,\infty]{\times}\Sigma\big)\qquad\widehat{\lambda}:=\begin{cases}\lambda &\text{on $V$}\\e^r\alpha & \text{on $\mathbb{R}{\times}\Sigma$.}\end{cases}\]
Write $V=V_\alpha$ and let $\beta$ be a different contact form for $\xi$. Using the obvious embedding $\mathbb{R}{\times}\Sigma\hookrightarrow\widehat{V}$, one can think of $\Sigma_\beta$ as a contact hypersurface in $\widehat{V}$ which bounds a compact region $V_\beta\subset\widehat{V}$. Note that the completions of $(V_\beta, \widehat{\lambda}|_{V_\beta})$ and $(V_\alpha,\lambda)$ are naturally identified with $(\widehat{V},\widehat{\lambda})$.\\
More generally, a diffeomorphism $\varphi : \widehat{V}{\rightarrow} \widehat{W}$ between two Liouville domains $(V,\lambda)$ and $(W,\mu)$ is a \textit{Liouville isomorphism} if $\varphi^\ast\widehat{\mu}=\widehat{\lambda}+dg$ for a compactly supported function $g$. For example $(V_\alpha,\lambda)$ and $(V_\beta,\widehat{\lambda}|_{V_\beta})$ as above with completion $(\widehat{V}, \widehat{\lambda})$ are Liouville isomorphic via the identity $Id: \widehat{V}{=}\widehat{V}_\alpha\rightarrow\widehat{V}_\beta{=}\widehat{V}$. Note that the contact forms on $\partial V_\alpha$ and $\partial V_\beta$ are different, but the contact structures agree. Therefore we think of $(V_\alpha,\widehat{\lambda}|_{V_\alpha})$ and $(V_\beta,\widehat{\lambda}|_{V_\beta})$ as the same filling for $(\Sigma,\xi)$, which leads us to:
\begin{defn}\label{filling}
 Let $(\Sigma,\xi)$ be a contact manifold. If there exists a Liouville domain $(V,\lambda)$ such that $\partial V=\Sigma$ and $\xi=\ker \lambda|_\Sigma$, then we call the equivalence class of $(V,\lambda)$ under Liouville isomorphisms an \textbf{\textit{exact (contact) filling}} of $(\Sigma,\xi)$.
\end{defn}
We will see that Symplectic Homology is invariant under Liouville isomorphisms (cf.\ Prop.\ \ref{propSHinv}) thus providing an invariant for contact structures (with a filling).\bigskip\\
A Hamiltonian on $\widehat{V}$ is a smooth $S^1$-family of functions $H_t: \widehat{V}\rightarrow \mathbb{R}$ with Hamiltonian vector field $X_H^t$ defined by $\omega(\cdot,X^t_H)=dH_t$ for each $t\in S^1$. The action of a loop $x: S^1{\rightarrow} \widehat{V}$ with respect to $H$ is defined by
\[\mathcal{A}^H(x)=\int^1_0 x^\ast\lambda - \int^1_0 H_t(x(t)) dt.\]
The critical points of the functional $\mathcal{A}^H$ are exactly the closed 1-periodic orbits of $X^t_H$. We denote the set of these 1-periodic orbits by $\mathcal{P}(H)$. Let $J_t$ denote an $S^1$-family of $\omega$-compatible almost complex structures, i.e.\ $\omega(\cdot,J_t\cdot)$ is a Riemannian metric on $V$ for every $t$. The $L^2$-gradient of $\mathcal{A}^H$ with respect to this metric is then given by $\nabla\mathcal{A}^H(x)=-J(\partial_t x-X_H^t)$. An $\mathcal{A}^H$-gradient trajectory $u:\mathbb{R}{\times} S^1\rightarrow \widehat{V}$ is hence a solution of the partial differential equation:
\begin{equation}\label{eqast}
 \partial_s u-\nabla\mathcal{A}^H(u)=\partial_s u + J(\partial_t u-X_H^t)=0\qquad\Leftrightarrow \qquad\, \big(Du-X_H^t{\otimes} dt\big)^{0,1}=0.\quad
\end{equation}
For the second equation, $Du$ is viewed as a 1-form on $\mathbb{R}{\times} S^1$ with values in $TV$. The antiholomorphic part of such forms is given by $\beta^{0,1}:=\frac{1}{2}(\beta+J\beta j)$, where $j$ is the standard almost complex structure on $\mathbb{R}{\times} S^1$, defined by $j\partial_s=\partial_t$.\\
If $H_s$ is a homotopy of Hamiltonians, then $\mathcal{A}^{H_s}$-gradient trajectories are solutions of (\ref{eqast}) with $X_H^t$ and $J$ depending on $s$.

\subsection{The No-Escape Lemma}\label{secNoEsc}
For the construction of symplectic (co)homology we look at solutions $u$ of (\ref{eqast}) satisfying $\displaystyle\lim_{s\rightarrow\pm\infty}u(s,t)=x_{\pm}(t)\in\mathcal{P}(H)$. In general, these solutions might not stay in a compact subset of $\widehat{V}$, even for $x_\pm$ fixed. So it could be that the moduli space of these solutions is neither compact nor has a suitable compactification. However, the No-Escape Lemma below shows that for certain pairs $(H,J)$ all such $u$ stay in a compact set.\\
Originally, such a result was proved by Viterbo, \cite{Vit}, as the Maximum Principle, which states that a solution $u$ cannot cross a convex hypersurface $\Sigma$ in a symplectization, if both asymptotes of $u$ are below $\Sigma$. Seidel then gave in \cite{Sei}, 3.21, a generalization which works for varying $\Sigma$, if the slope of $H$ increase exponentially between $H_+$ and $H_-$. Finally, in \cite{AbouSei}, 7.2, and \cite{Ritter}, Lem.\ 19.3, Abouzaid/Seidel and Ritter (adaptation to $SH$) introduced the No-Escape Lemma, which has the advantage that the ambient manifold does not need to be a symplectization. Here, we generalize this No-Escape Lemma to situations where $\Sigma$ is allowed to vary in $s$.\\
In order to state the lemma in full generality, let $(W,d\lambda)$ be an exact symplectic manifold with compact negative contact boundary and such that the flow $\varphi_Y^t$ of the Liouville vector field $Y$ exists for all $t\geq 0$. Then $\varphi_Y$ provides a symplectic embedding of $\big([0,\infty){\times}\partial W, d(e^r\alpha)\big)$ into $W$. For example, consider $\big([-\delta,\infty){\times}\Sigma, d(e^r\alpha)\big)$ inside $(\widehat{V},\widehat{\omega})$.\pagebreak[1]\\
Let $f_s:\partial W\rightarrow\mathbb{R}$ be a smooth family of functions, such that for some $s_0\geq 0$ holds $f_s\equiv f_{\pm s_0}$ for $|s|\geq s_0$. They define on $[0,\infty)\times\partial W$ an $s$-dependent coordinate change by $r_s:=r-f_s$ and a compact family of contact hypersurfaces $\Sigma_s$ by 
\[\Sigma_s:=\{r_s\equiv R_0\}=\big\{(R_0+f_s(p),p)\,\big|\,p\in\partial W\big\} \quad \text{ for a constant }\quad R_0\geq - \min_{s,p} f_s(p).\]
Let $J_s$ be an $s$-dependent family of almost complex structures, which are of contact type along $\Sigma_s$, meaning that $J_s^\ast\lambda = d(e^{r_s})$ holds for fixed $s$ at all points $p\in\Sigma_s$. Let $H_s: W\rightarrow \mathbb{R}$ be a homotopy of Hamiltonians such that
\begin{align}
 &\bullet& H_s(r,p)=h_s\big(e^{r-f_s(p)}\big)&=h_s(e^{r_s})&&\text{ near }\Sigma_s,\notag\\
 &\bullet&\label{eqconditiononH} \partial_s\Big(H_s-h_s(e^{R_0})+e^{R_0}\cdot h_s'(e^{R_0})\Big)&\leq 0&&\text{ everywhere on $W$.}
\end{align}
Finally, let $S\subset \mathbb{R}\times S^1$ be a compact Riemann surface with smooth boundary.
\begin{lemma}[\textbf{No-Escape Lemma}]\label{maxprinc}~\\
Let $W, S, J, r_s, R_0$ and $H_s$ be as above. Assume that for a solution $u: S\rightarrow W$ of (\ref{eqast}) holds that $u(s,t)\in\Sigma_s$ for all $(s,t)\in\partial S$ and $e^{r_s}\circ u(s,t)\geq e^{R_0}$ for all $(s,t)\in S$. Then
\[u(s,t)\in \bigcup_{s\in\mathbb{R}} \Sigma_s =\big\{p\in W\big|\exists\, s\in\mathbb{R}: p\in\Sigma_s\big\}\quad \forall\; (s,t)\in S.\]
\end{lemma}
\begin{remark}
\begin{itemize}
 \item If $H_s$ and $f_s$ are independent of $s$, then condition (\ref{eqconditiononH}) is empty, i.e.\ the No-Escape Lemma holds for all $(H,J)$ that are cylindrical along a fixed $\Sigma$.
 \item If $H_s=\mathfrak{a_s}e^{r_s}+\mathfrak{b}_s$ along $\Sigma_s$, then (\ref{eqconditiononH}) reads as $\partial_s\big(H_s-\mathfrak{b}_s)\leq 0$.
 \item If $W=[0,\infty)\times\partial W$ and $H_s=\mathfrak{a_s}e^{r_s}+\mathfrak{b}_s$ everywhere, then (\ref{eqconditiononH}) reads as $\partial_s (\mathfrak{a}_s e^{r-f_s})\leq 0$, which is equivalent to $\partial_s (\log \mathfrak{a}_s-f_s)\leq 0$.
 \item If $W=[0,\infty)\times\partial W$ and $f_s=0$, then (\ref{eqconditiononH}) can be replaced by $\partial_s h'\leq 0$, as
 \[(\partial_s h_s)(e^r\circ u)-(\partial_s h_s)(e^{R_0})+e^{R_0} (\partial_s h_s')(e^{R_0})=\int_{e^{R_0}}^{e^r\circ u}\mspace{-20mu}\partial_s h'(t)\,dt+e^{R_0} (\partial_s h_s')(e^{R_0}).\]
\end{itemize}
\end{remark}
The No-Escape Lemma and Sard's theorem imply the following corollary.
\begin{corollary}\label{Cornoesc}
 Let $V_0\subset \widehat{V}$ be a relatively compact open set with positive contact boundary $\partial V_0$, let $H: \widehat{V}\rightarrow \mathbb{R}$ be a Hamiltonian satisfying (\ref{eqconditiononH}) on $W:=\widehat{V}\setminus V_0$ and let $J$ be an almost complex structure which is cylindrical along a collar neighborhood of $\partial V_0$. Then any solution $u:\mathbb{R}\times S^1\rightarrow \widehat{V}$ to (\ref{eqast}) with asymptotes in $V_0$ stays inside $V_0$ for all time.
\end{corollary}
\begin{proof}[\textbf{Proof of the No-Escape Lemma (cf.\ \cite{Ritter}, Lem.\ 19.3):}]~\\
 At first, we calculate $\lambda$ applied to the Hamiltonian vector field on $\Sigma_s$:
 \begin{equation*}
  \begin{aligned}\lambda(X_{H_s})=d\lambda(Y,X_{H_s})=dH_s(Y)&=\partial_r H_s(r,y)\\&=h_s'(e^{r-f_s(p)}) e^{r-f_s(p)}=h'_s(e^{R_0}) e^{R_0}\end{aligned}\tag{$\ast$}
 \end{equation*}
 where the last equality holds only on $\Sigma_s$, as there $r-f_s(p)=r_s=R_0$.\pagebreak[1]\\ We define the energy $E_S(u)$ of $u$ over $S$ as $E_S(u):=\int_S ||\partial_s u||^2 ds\wedge dt$. Clearly, $E_S(u)$ is non-negative. Using a trick of M. Abouzaid, we will show that $E_S(u)\leq 0$ and hence $E_S(u)=0$, so that $\partial_s u \equiv 0$. As $u|_{\partial S}\subset \bigcup_{s\in\mathbb{R}} \Sigma_s$ and $S\subset \mathbb{R}\times S^1$, this implies that $u(s,t)\in \bigcup_{s\in\mathbb{R}} \Sigma_s$ for all $(s,t)\in S$. To prove $E_S(u)\leq 0$, we calculate:
 \begin{align*}
  &\phantom{\,=\,}\,E_S(u) =\int_S||\partial_s u||^2ds{\wedge} dt =\int_S d\lambda(\partial_s u, J\partial_s u)ds{\wedge} dt\\
  &=\int_S d\lambda(\partial_s u,\partial_t u)-d\lambda(\partial_s u, X_{H_s})ds{\wedge} dt\\
  &=\int_S u^\ast d\lambda - dH_s(\partial_s u)ds{\wedge} dt\\
  &=\int_S u^\ast d\lambda -\partial_s \big(H_s(u)\big)ds{\wedge} dt + (\partial_s H_s)(u) ds{\wedge} dt\\
  &=\int_S u^\ast d\lambda- d\big(H_s(u) dt\big) +(\partial_s H_s)(u) ds{\wedge} dt\allowdisplaybreaks\\
  &=\int_{\partial S} u^\ast\lambda- H_s(u) dt +\int_S(\partial_s H_s)(u) ds{\wedge} dt\\
  &=\int_{\partial S} u^\ast\lambda-\Big(\lambda(X_{H_s})(u)-\lambda(X_{H_s})(u)-H_s(u)\Big)dt +\int_S(\partial_s H_s)(u) ds{\wedge} dt\allowdisplaybreaks\\
  &\overset{(\ast)}{=}\int_{\partial S}\lambda\big(Du-X_{H_s}{\otimes} dt\big)+\int_{\partial S}\Big(h'_s\big(e^{R_0}(u)\big) e^{R_0}(u)-h_s(e^{R_0})\Big)dt+\int_S(\partial_s H_s)(u) ds{\wedge} dt\\
  &=\int_{\partial S}-\lambda J\big(Du - X_{H_s}{\otimes} dt\big)j + \int_S\partial s \Big(h'_s\big(e^{R_0}(u)\big) e^{R_0}(u)-h_s(e^{R_0})\Big)+(\partial_s H_s)(u) ds{\wedge} dt\\
  &\overset{(\ref{eqconditiononH})}{\leq}\int_{\partial S}-de^{r_s}(Du-X_{H_s}{\otimes} dt)j=\int_{\partial S}-de^{r_s}(Du)j.
 \end{align*}
Here, we used that orbits of $X_{H_s}$ stay inside level sets of $e^{r_s}$, so that $de^{r_s}(X_{H_s})=0$. To calculate the last integral, let $n$ be the outward normal direction along $\partial S\subset S$. Then $(n,jn)$ is an oriented frame and $\partial S$ is oriented by $jn$. So along  $\partial S$ holds
\[-de^{r_s}(Du)j(jn)=-d(e^{r_s}\circ u)(-n)\leq 0,\]
as in the inward direction $-n$, $e^{r_s}\circ u$ can only increase since $e^{r_s}\circ u$ attains its minimum $e^{R_0}$ along $\partial S$. So $E_S(u)\leq 0$ and hence $E_S(u)=0$.
\end{proof}

A 1-periodic orbit $x\in\mathcal{P}(H)$ is called \textit{non-degenerate} if the flow $\varphi^t_{X_H}$ of $X_H$ satisfies $\det \big(D\varphi^1_{X_H}(x(0))-\mathbbm{1}\big)\neq 0$. It is called \textit{transversely non-degenerate} if $\mathcal{N}:=\big\{y(0)\,\big|\,y\in\mathcal{P}(H)\big\}$ is a submanifold of $V$ near $x$ and $\ker\big(D\varphi^1_{X_H}(y(0))-\mathbbm{1}\big)=T_{y(0)}\mathcal{N}$ for all $y\in\mathcal{P}(H)$ near $x$. Note that $\mathcal{N}$ is always closed and consists of finitely many points if all orbits are non-degenerate.\\
In view of the No-Escape Lemma (Lem.\ \ref{maxprinc}) we make the following definitions:
\begin{itemize}
 \item A Hamiltonian $H$ is \emph{admissible}, writing $H\in Ad(V)$, if all 1-periodic orbits of $X_H$ are (transversely) non-degenerate and if $H$ is \emph{(weakly) linear at infinity}, that is if there exist $\mathfrak{a},\mathfrak{b},R\in\mathbb{R}$ and $f\in C^\infty(\Sigma)$ such that $\mathfrak{a}\not\in spec(e^{f(p)}{\cdot} \alpha)$ and $H$ is on $[R,\infty){\times}\Sigma\subset \widehat{V}$ of the form
 \[ H(r,p)=\mathfrak{a}\cdot e^{r-f(p)}+\mathfrak{b}.\]
 \item A homotopy $H_s$ between admissible Hamiltonians $H_\pm$ is admissible if there exist $S,R\geq0$ such that $H_s=H_\pm$ for $\pm s\geq S$ and $H_s$ has on $[R,\infty)\times\Sigma$ the form
 \[H_s=\mathfrak{a}_s\cdot e^{r-f_s(p)}+\mathfrak{b}_s \qquad \text{ with }\qquad \partial_s\big(\log \mathfrak{a}_s-f_s(p)\big)\leq 0.\]
 \item A possibly $s$-dependent almost complex structure $J$ is admissible for a Hamiltonian/homotopy $H$, if for some $\displaystyle R_0\geq \min \{R-f_s(p)\,|\,p\in\Sigma,s\in\mathbb{R}\}$ holds that  $J_s$ is of contact type near $\Sigma_s:=\big\{r{-}f_s(p)=R_0\big\} \subset\mathbb{R}{\times}\Sigma$, meaning that
 \[\lambda\circ J_s=d\big(e^{r-f_s}\big) \qquad\text{ holds for $s$ fixed and all }(r,p)\in(-\veps,\veps){\times}\Sigma_s.\]
\end{itemize}

\subsection{Symplectic Homology and cohomology}\label{secsymhom}
The following construction of Symplectic Homology in a Morse-Bott setting was first described by Bourgeois and Oancea, \cite{BourOan1}, for the case where the manifold $\mathcal{N}$ of contractible 1-periodic orbits of $H$ consists of isolated circles. Although not stated explicitly, their methods are general enough to work also if $\mathcal{N}$ is of higher dimension. See also \cite{Fra},  A and \cite{FauckThesis} for similar Morse-Bott constructions using flow lines with cascades.\\
For $H\in Ad(V)$ (with all 1-periodic orbits transversely non-degenerate), we define the Hamiltonian Floer homology $FH_\ast(H,h)$ as follows: Choose a Morse-function $h$ (and a Riemannian metric) on $\mathcal{N}$. Then let $\mathcal{P}(H,h)$ consist of the critical points of $h$, let $FC_\ast(H,h)$ be the $\mathbb{Z}_2$-vector space generated by $\mathcal{P}(H,h)$ and let $\mathcal{M}(x_-,x_+)$ consist of unparameterized flow lines with cascades between $x_\pm\in\mathcal{P}(H,h)$. Here, a flow line with cascades is a tuple $(u_1,...,u_m)$ whose components are solutions of (\ref{eqast}) and satisfy
 \begin{itemize}
  \item $\displaystyle \lim_{s\rightarrow-\infty} u_1(s,0)$ lies in the unstable manifold of $x_-$ and $\displaystyle\lim_{s\rightarrow+\infty} u_m(s,0)$ lies in the stable manifold of $x_+$, both with respect to the gradient flow of $h$ on $\mathcal{N}$,
  \item for $i=1,...,m{-}1$, the limit orbits $\displaystyle \lim_{s\rightarrow+\infty}u_i$ and $\displaystyle \lim_{s\rightarrow-\infty}u_{i+1}$ lie in the same component of $\mathcal{N}$ and are connected by a positive gradient flow line of $h$ with finite (possibly zero) length.
 \end{itemize}
For a generic Riemannian metric on $\mathcal{N}$ and a generic almost complex structure $J$, the space $\mathcal{M}(x_-,x_+)$ is a manifold. Its zero-dimensional component $\mathcal{M}^0(x_-,x_+)$ is compact and hence a finite set. Let $\#_2\mathcal{M}^0(x_-,x_+)$ denote its cardinality modulo 2. The boundary operator $\partial: FC_\ast(H,h)\rightarrow FC_\ast(H,h)$ is then defined as the linear extension of
\[\partial x_+:=\sum_{x_-\in\mathcal{P}(H)}\#_2\mathcal{M}^0(x_-,x_+)\cdot x_-.\]
A standard argument in Floer theory shows $\partial^2=0$ and the resulting homology is denoted by $FH_\ast(H,h)$. To an admissible homotopy $H_s$ between $H_\pm\in Ad(V)$ and Morse functions $h_\pm$ on the manifolds $\mathcal{N}_\pm$ of 1-periodic orbits of $H_\pm$, we consider for ${x_\pm\in\mathcal{P}(H_\pm, h_\pm)}$ the moduli space $\mathcal{M}_s(x_-,x_+)$ of unparameterized flow lines with cascades from  $x_-$ to $x_+$ where one cascade $u_i$ is an $s$-dependent $\mathcal{A}^{H_s}$-gradient trajectory, while $u_j$ is a $\mathcal{A}^{H_-}$-gradient trajectory for $j<i$ and a $\mathcal{A}^{H_+}$-gradient trajectory for $j>i$. On chain level, we define a map $\sigma_\sharp(H_-,H_+):FC_\ast(H_+,h_+)\rightarrow FC_\ast(H_-,h_-)$ via
\[\sigma_\sharp(H_-,H_+)x_+=\sum_{x_-\in\mathcal{P}(H_-,h_-)}\#_2\mathcal{M}^0_s(x_-,x_+)\cdot x_-.\]
By considering the compactification of $\mathcal{M}^1_s(x_-,x_+)$, one can show that $\partial\circ \sigma_\sharp =\sigma_\sharp\circ \partial$, so that $\sigma_\sharp(H_-,H_+)$ descends to a map  $\sigma_\ast(H_-,H_+):FH_\ast(H_+,h_+)\rightarrow FH_\ast(H_-,h_-)$, called continuation map. Considering homotopies of homotopies, one can show that $\sigma_\ast(H_-,H_+)$ is independent of the chosen homotopy. For $H_1, H_2, H_3\in Ad(V)$ the following composition rule holds
\[\sigma_\ast(H_1,H_3)=\sigma_\ast(H_1,H_2)\circ\sigma_\ast(H_2,H_3).\]
We introduce a partial order $\prec$ on $Ad(V)$ by saying $H_+\prec H_-$ if for some $R\in\mathbb{R}$ holds on $[R,\infty){\times}\Sigma$ that either $H_+{-}H_-$ is constant or $H_+\leq H_-$. Observe that admissibility of a homotopy $H_s$ between $H_-$ and $H_+$ implies that $H_+\prec H_-$. It follows that $\big(FH_\ast(H,h),\sigma_\ast\big)$ is a direct system over the directed set $\big(Ad(V),{\prec}\big)$. The Symplectic Homology $SH_\ast(V)$ is then the direct limit of this system:
\[SH_\ast(V):=\varinjlim_{H\in Ad(V)} FH_\ast(H,h).\]
One obtains Symplectic Cohomology by dualizing this construction. Explicitly, the cochain groups $FC^\ast(H,h)$ are the duals of $FC_\ast(H,h)$. As $FC_\ast(H,h)$ is $\mathbb{Z}_2$-generated by the finite set $\mathcal{P}(H,h)$, $FC^\ast(H,h)$ is also  $\mathbb{Z}_2$-generated by $\mathcal{P}(H,h)$. The coboundary operator $\delta$ is the dual of $\partial$. The analogue construction of $\sigma^\sharp(H_-,H_+)$ associated to a homotopy $H_s$ between $H_\pm\in Ad(V)$ yields continuation maps in the opposite direction $\sigma^\ast(H_-,H_+):FH^\ast(H_-,h_-)\rightarrow FH^\ast(H_+,h_+)$. By taking the same partial order on $Ad(V)$ as for homology, we obtain hence an inverse system and $SH^\ast(V)$ is its inverse limit: $SH^\ast(V):=\varprojlim FH^\ast(H,h)$.\\
Note that $SH_\ast(V)$ and $SH^\ast(V)$ do not depend on $H, J$ and $h$ as all groups $FH(H,h)$ for any of these choices are included in the direct/inverse limit.\\
A \emph{cofinal sequence} $(H_n)\subset Ad(V)$ is a sequence such that $H_n\prec H_{n+1}$ and for any $H\in Ad(V)$ holds $H\prec H_n$ for some $n\in\mathbb{N}$. More general, a set $\mathcal{F}\subset Ad(V)$ is cofinal if for any $H\in Ad(V)$ exists $F\in \mathcal{F}$ with $H\prec F$. For $\mathcal{F}\subset Ad(V)$ cofinal holds that $\displaystyle SH_\ast(V)=\varinjlim FH_\ast(F)$ and $\displaystyle SH^\ast(V)=\varprojlim FH_\ast(F)$, where $F\in\mathcal{F}$.
\begin{proposition}\label{propSHinv}
 $SH_\ast(V)$ and $SH^\ast(V)$ are invariant under Liouville isomorphisms. They are therefore invariants of the exact filling $(V,\lambda)$ of $\big(\partial V,\lambda|_{\partial V}\big)$.
\end{proposition}
\begin{proof}
 Let $(V,\lambda)$ and $(W,\mu)$ be two Liouville domains and $\varphi:\widehat{V}\rightarrow\widehat{W}$ be a Liouville isomorphism. In \cite{Sei}, page 3, it is shown that there exists $R\in\mathbb{R}$ such that $\varphi$ takes on $[R,\infty){\times}\partial V\subset \widehat{V}$ the form 
 \[\varphi(r,p)=(r{-}f(p),\psi(p)),\]
 where $f\in C^\infty(\partial V)$ and $\psi: \partial V{\rightarrow} \partial W$ satisfies $\psi^\ast \mu|_{\partial W}=e^f{\cdot} \lambda|_{\partial V}$. This implies for any $H{\in} Ad(W)$ that $\varphi^\ast H(r,p)=\mathfrak{a}e^{r-f(p)-g(\psi(p))}{+}\mathfrak{b}$ wherever $H(r,p)=\mathfrak{a}e^{r-f(p)}{+}\mathfrak{b}$ holds. Hence we have $\varphi^\ast H{\in} Ad(V)$, i.e.\ $\varphi^\ast\mspace{-3mu} Ad(W)\subset Ad(V)$. As $\varphi^{-1}$ is also a Liouville isomorphism, it follows that $\varphi^\ast\mspace{-3mu} Ad(W)= Ad(V)$. Moreover as $\varphi^\ast d\mu=d\lambda$, we have $\varphi^\ast X_H=X_{\varphi^\ast H}$, so that $\varphi^\ast$ provides a 1-1 correspondence between $P(H,h)$ and $P(\varphi^\ast H,\varphi^\ast h)$ (where $h$ is a Morse function on $\mathcal{N}$). For an admissible almost complex structure $J$ on $\widehat{W}$, we find that $\varphi^\ast J$ is also admissible on $\widehat{V}$. It follows that $FH_\ast(H,h)$ and $FH_\ast(\varphi^\ast H,\varphi^\ast h)$ are isomorphic, as the moduli spaces $\mathcal{M}^0(y,x)$ and $\mathcal{M}^0(\varphi^\ast y,\varphi^\ast x)$ are all diffeomorphic via $\varphi^\ast$ for all $x,y\in\mathcal{P}(H,h)$. The same argument shows that for any $H_\pm\in Ad(W)$ with $H_+\prec H_-$ the following diagram commutes\pagebreak[1]
 \[\begin{xy}\xymatrix{FH_\ast(H_+,h_+)\ar[d]_{\varphi^\ast}\ar[rrr]^{\sigma_\ast(H_-,H_+)} &&& FH_\ast(H_-,h_-)\ar[d]_{\varphi^\ast}\\
 FH_\ast(\varphi^\ast H_+,\varphi^\ast h_+)\ar[rrr]^{\sigma_\ast(\varphi^\ast H_-,\varphi^\ast H_+)} &&& FH_\ast(\varphi^\ast H_-,\varphi^\ast h_-)}.\end{xy}\]
 Hence we have $\displaystyle SH_\ast(W)=\hspace{-5 pt}\varinjlim_{H\in Ad(W)} \hspace{-10 pt}FH_\ast(H,h) \cong \hspace{-5 pt}\varinjlim_{\varphi^\ast H\in \varphi^\ast\mspace{-3mu} Ad(W)} \hspace{-15 pt}FH_\ast(\varphi^\ast H, \varphi^\ast h)=SH_\ast(V)$.\\
 The argument stays valid, if we replace the homology groups $FH_\ast(H,h)$ by the cohomology groups $FH^\ast(H,h)$, which shows $SH^\ast(W)=SH^\ast(V)$.
\end{proof}

\subsection{Action filtration}\label{sectrunc}
The action functional $\mathcal{A}^H$ provides filtrations of $SH_\ast(V)$ and $SH^\ast(V)$ as follows: For $H\in Ad(V)$, a Morse function $h$ on $\mathcal{N}$ and $b\in\mathbb{R}{\cup}\{\pm \infty\}$ consider the subchain groups $FC_\ast^{<b}(H){:=}FC_\ast^{<b}(H,h)\subset FC_\ast(H,h)$ generated by whose $x\in\mathcal{P}(H,h)$ with $\mathcal{A}^H(x){<}b$. For $a{<}b$, we set $FC^{[a,b)}_\ast(H):=FC^{<b}_\ast(H)\big/FC^{<a}_\ast(H)$. Similar definitions work for the intervals $(a,b), (a,b]$ etc.\ or ${\leq} b, {>}b, {\geq} b$. Note that $FC^{[a,b)}_\ast(H)=FC^{(a,b)}_\ast(H)$ if ${a\not\in\mathcal{A}^H(\mathcal{P}(H))}$. All subsequent discussions hold for any action restriction, yet we will write them  only for $(a,b)$ and ${<}b$. Lem.\ \ref{monolem} below shows that the boundary operator $\partial$ reduces the action. It induces therefore a boundary operator on $FC_\ast^{(a,b)}(H)$ which defines $FH_\ast^{(a,b)}(H)$.
\begin{lemma}\label{monolem}~\\
 Hamiltonians\big/homotopies $H$ and solutions $u$ of (\ref{eqast}) with $\displaystyle\lim_{s\rightarrow\pm\infty}u=x_\pm\in\mathcal{P}(H)$ satisfy
 \begin{equation}\label{eqAction}
  \mathcal{A}^H(x_+)-\mathcal{A}^H(x_-)\geq-\int_{-\infty}^\infty \int_0^1(\partial_s H)(u)dt\,ds.                                                                                                                                                         \end{equation}
  If $H_s\equiv H$ or $\partial_s H\leq 0$, then $\mathcal{A}^H(x_+)\geq\mathcal{A}^H(x_-)$.
\end{lemma}
\begin{proof}
\[\mathcal{A}^H(x_+)-\mathcal{A}^H(x_-)=\int^\infty_{-\infty}\mspace{-10mu}\partial_s\mathcal{A}^H(u(s))ds=\int^\infty_{-\infty} \mspace{-10mu}||\nabla\mathcal{A}^H||^2ds-\int^\infty_{-\infty}\int^1_0 \mspace{-10mu}(\partial_s  H)(u)dt\,ds.\qedhere\]
\end{proof}
Let $H_\pm\in Ad(V)$ with $H_-{\geq}H_+$ everywhere. Then we can choose an admissible homotopy $H_s$ between them with $\partial_s H\leq 0$. It follows from Lem.\ \ref{monolem} that the associated continuation map $\sigma_\sharp(H_-,H_+)$ also decreases action and induces a map
\[\sigma_\ast(H_-,H_+): FH_\ast^{(a,b)}(H_+)\rightarrow FH_\ast^{(a,b)}(H_-).\]
The truncated Symplectic Homology in the action window $(a,b)$ is then defined as the direct limit under these maps: $SH_\ast^{(a,b)}(V):=\varinjlim FH_\ast^{(a,b)}(H)$. \\
For Symplectic Cohomology, the only difference is that the coboundary operator $\delta$ increases action, so that we have to consider $FC^\ast_{>a}(H)\subset FC^\ast(H)$ as a natural subcomplex. All other truncated groups can then be defined accordingly. Then, we obtain $FH^\ast_{>a}(H)$ and $FH^\ast_{(a,b)}(H)$ as filtered Floer cohomology groups. For homotopies $H_s$ with $\partial_s H\leq 0$, the continuation maps are also well-defined on truncated groups and we obtain as inverse limits $SH^\ast_{(a,b)}(W{\subset}V)=\varprojlim FH^\ast_{(a,b)}(H)$.
\begin{remark}
 Without further restrictions, we have for all $a>-\infty$ and any $b$
\[SH^{(a,b)}_\ast(V)=0\qquad\text{ and }\qquad SH^{(-\infty,b)}_\ast(V)=SH_\ast(V).\]
Indeed, for any cofinal sequence of Hamiltonians $(H_n)$ choose an increasing sequence $(\beta_n)\subset\mathbb{R}$ such that $\displaystyle\beta_n>\max_{x\in\mathcal{P}(H_n)}\mathcal{A}^{H_n}(x)$. Define $K_n:=H_n+\beta_n-a$ and $L_n:=H_n+\beta_n-b$, which are both cofinal sequences satisfying
\[\max_{x\in\mathcal{P}(K_n)}\mathcal{A}^{K_n}(x)=\max_{x\in\mathcal{P}(H_n)}\mathcal{A}^{H_n}(x)-\beta_n+a<a\quad\text{ and }\quad \max_{x\in\mathcal{P}(L_n)}\mathcal{A}^{L_n}(x)<b.\]
It follows that $FC^{(a,b)}_\ast(K_n)=FH_\ast^{(a,b)}(K_n)=0$ for all $n$ and hence $SH^{(a,b)}_\ast(V)=0$, while $FC_\ast^{(-\infty,b)}(L_n)=FC_\ast(L_n)$ for all $n$ and hence $SH^{(-\infty,b)}_\ast(V)=SH_\ast(V)$. Similar phenomenons hold in cohomology.
\end{remark}
To obtain a meaningful action filtered version of $SH$, we have to restrict the set of admissible Hamiltonians. We will require that all $H$ are negative inside a fixed Liouville subdomain $W\subset \widehat{V}$ bounded by a contact hypersurface $\partial W$\footnote{Other possibilities are $H|_{\partial V}<0$ which leads to the $V$-shaped homology/Rabinowitz-Floer homology or $H|_{V\setminus W}<0$ which leads to the Symplectic Homology of a cobordism. See \cite{CieOan} for more details.}. In particular, one can take $W=V$. We write $SH^{(a,b)}(W{\subset}V)$ or $SH^{(a,b)}(V)$ if $W=V$, for the direct limit of these Hamiltonians\footnote{These Hamiltonians coincide with whose defining $SH_\ast(W)$ in the sense of \cite{CieOan}. However, $SH_\ast(W)\neq SH_\ast(W{\subset}V)$ in general, but $SH^{\geq 0}_\ast(W{\subset}V)=SH_\ast(W)$ as shown in Cor.\ \ref{transfer}.}. Note that different choices of $W\subset V$ give different filtrations of $SH_\ast(V)$! \\
For the definition of $FH^{(a,b)}_\ast(H)$ it suffices that only the 1-periodic orbits $x$ of $X_H$ with $\mathcal{A}^H(x)\in(a,b)$ are non-degenerate, as the others are discarded. Therefore, we call a Hamiltonian $H$ admissible for $SH^{(a,b)}_\ast(W{\subset}V)$, writing $H\in Ad^{(a,b)}(W{\subset}V)$, if it satisfies the following assumptions:
\begin{itemize}
 \item $H<0$ on $W$,
 \item $H(r,p)=\mathfrak{a}{\cdot} e^{r-f(p)}+\mathfrak{b}$ on $[R,\infty){\times}\Sigma$ for some $R$ and $f:\Sigma\rightarrow\mathbb{R}$,
 \item all $x\in\mathcal{P}(H)$ with $\mathcal{A}^H(x)\in(a,b)$ are (transversely) non-degenerate.
\end{itemize}
The partial order on $Ad^{(a,b)}(W{\subset}V)$ is given by $H\prec K$ if $H\leq K$ everywhere. Note that $Ad(W{\subset}V)=Ad^{(-\infty,\infty)}(W{\subset} V)\subset Ad(V)$ is a cofinal subset, so that $SH_\ast(W{\subset}V)=SH_\ast(V)$. For $H\in Ad(W{\subset}V)$, the projection $\pi:\;FC_\ast(H)\rightarrow FC_\ast^{\geq b}(H)$ and the short exact sequence
\[0\rightarrow FC_\ast^{(a,b)}(H)\rightarrow FC_\ast^{(a,c)}(H)\rightarrow FC_\ast^{(b,c)}(H)\rightarrow 0\]
induce in homology a map $\pi:\,FH_\ast(H)\rightarrow FH_\ast^{\geq b}(H)$ and a long exact sequence
\[\dots\rightarrow FH_\ast^{(a,b)}(H)\rightarrow FH_\ast^{(a,c)}(H)\rightarrow FH_\ast^{(b,c)}(H)\rightarrow\dots\]
Applying the direct limit yields a map $SH_\ast(V)=SH_\ast(W{\subset}V)\overset{\pi}\longrightarrow SH_\ast^{\geq b}(W{\subset}V)$ and (as $\displaystyle\varinjlim$ is an exact functor) a long exact sequence
\[\dots\rightarrow SH^{(a,b)}_\ast(W{\subset}V)\rightarrow SH^{(a,c)}_\ast(W{\subset}V)\rightarrow SH^{(b,c)}_\ast(W{\subset}V)\rightarrow\dots\]
Note that $\varprojlim$ is in general not an exact functor. However, it preserves exactness if all modules $FH^\ast(H,h)$ are finite dimensional vector spaces, which is here the case (see \cite{bourbaki2}, §6, No. 3, Prop.\ 4 or \cite{eilenberg}, Chap.\ VIII, Sec.\ 5, Thm.\ 5.7). Hence, the above exact sequences make sense also in cohomology.

\subsection{Positive Symplectic Homology}\label{secPosSym}
The most used filtered group is the positive Symplectic Homology $SH^+(V):=SH^{>\veps}(V)$, where $\veps>0$ is an arbitrary small positive constant. Although the action filtration does in general depend on the chosen subdomain $V_a\subset \widehat{V}$ and hence on the contact form $\alpha$, $SH^+(V)$ can be defined independently of $\alpha$ and is also invariant under Liouville isomorphisms. To see this, we consider the following standard type of Hamiltonians $H$ on $\widehat{V}$. Let $\alpha_0$ be the contact form on $\partial V=\Sigma$, let $\alpha=e^{f}{\cdot}\alpha_0$ be any contact form, where $f:\Sigma \rightarrow \Sigma$ is a smooth function, and let $V_\alpha, \Sigma_\alpha=\partial V_\alpha$ be the associated subdomain as in \ref{secsetup}. We call $H$ adapted to $\alpha$ if 
\begin{itemize}
 \item $H$ is inside $V_\alpha$ equal to a $C^2$-small Morse function $g$ with $-\veps<g<0$, 
 \item $H$ is on $[-\veps,\infty){\times}\Sigma_\alpha$ of the form $H(r,p)=h(e^r)$ for a function $h$ with $h''{\geq} 0$ and ${h(e^r)=\mathfrak{a} e^r-\mathfrak{a}}$ on $[0,\infty){\times}\Sigma_\alpha$ with $\mathfrak{a}\not\in spec(\alpha)$.
\end{itemize}
 As for $H(r,p)=h(e^r)$ holds that $X_H(r,p)=h'(e^r){\cdot} R_\alpha(p)$, we find that any $H$ of this form has two types of 1-periodic orbits:
\begin{itemize}
 \item[\textbf{\textit{I}} :] constant orbits inside $V_\alpha$ corresponding to critical points,
 \item[\textbf{\textit{II}} :] non-constant orbits near $\Sigma_\alpha$ corresponding to $R_\alpha$-Reeb orbits of length $h'(e^r)$.
\end{itemize}
The action $\mathcal{A}^H(x)$ of an orbit $x$ on level $r{=}r_0$ is $h'(e^{r_0}){\cdot} e^{r_0}-h(e^{r_0})$. The action of orbits of type I is hence smaller than $\veps$ and for type II close to some $h'(e^r)\in spec(\alpha)$. With $\mathfrak{a_\alpha}:=\min spec(\alpha)>0$, we find that the action difference between type I and II is at least $\frac{1}{2}\mathfrak{a}_\alpha+\veps$ for $\veps$ small.\\
Note that $H$ takes the form $H(r,p)=h(e^{r-f(p)})$ with respect to coordinates ${(r,p)\in\mathbb{R}{\times}\Sigma}$ adapted to the reference hypersurface $\Sigma=\partial V$, so $H\in Ad(V)$. For a fixed Morse function $k$ on $\mathcal{N}$, we set $FC^0(H):=FC^{<\veps}(H,k)$ and $FC^+(H):=FC^{>\veps}(H,k)=FC(H,k)\big/FC^{<\veps}(H,k)$. Note that $FC^+(H)$ is generated by the orbits of type II, i.e.\ by the closed Reeb orbits.\\
For two contact forms $\alpha_\pm$ let $H_{\alpha_\pm}$ be Hamiltonians adapted to $\alpha_\pm$ respectively. If $H_{\alpha_+}\prec H_{\alpha_-}$, it is not difficult to find an admissible homotopy $H_{\alpha_s}$, such that $\alpha_s$ is a homotopy of contact forms for $\xi$, $H_{\alpha_s}$ is adapted to $\alpha_s$ for each $s$, $H_{\alpha_s}\equiv H_{\alpha_-}$ for $s{\leq} 0$ and $H_{\alpha_s}\equiv H_{\alpha_+}$ for $s{\geq} {+}1$. Note that $\partial_s H_{\alpha_s}\leq 0$ does not necessarily hold, so that continuation maps as described above may not be defined on $FH^+(H_{\alpha_\pm})$. However, $\partial_s H_{\alpha_s}$ is bounded from above, as the admissibility implies outside of a compact set 
\[\partial_s H_{\alpha_s}=\partial_s\big(\mathfrak{a}_s{\cdot}e^{r-f_s(p)}-\mathfrak{a}_s\big)\leq 0.\]
\begin{lemma}\label{lemactiondicho}
 If $\partial_s H_{\alpha_s}\leq \frac{1}{2}\mathfrak{a}_{\alpha_-}$, then there are well-defined continuation maps $\sigma_\ast\mspace{-2mu}(H_-, \mspace{-2mu}H_+)$ mapping $FH^0(H_+)$ into $FH^0(H_-)$ and $FH^+(H_+)$ into $FH^+(H_-)$.
\end{lemma}
\begin{proof}
 A solution $u$ of (\ref{eqast}) with $\displaystyle \lim_{s\rightarrow\pm\infty} u= x_\pm\in\mathcal{P}(H_{\alpha_\pm},k_\pm)$ satisfies by Lem.\ \ref{monolem}:
\[\mathcal{A}^{H_{\alpha_+}}(x_+)-\mathcal{A}^{H_{\alpha_-}}(x_-)\geq-\mspace{-5mu}\int_{-\infty}^\infty \int_0^1(\partial_s H_{\alpha_s})(u)dt\,ds\geq -\mspace{-5mu}\int_0^1 \int_0^1\frac{\mathfrak{a}_{\alpha_-}}{2}dt\,ds =- \frac{\mathfrak{a}_{\alpha_-}}{2}.\]
This implies that if $x_+$ is an orbit of type I, i.e.\ with action less than $\veps$, then $x_-$ is also of type I, as its action is at most $\veps+\frac{1}{2}\mathfrak{a}_{\alpha_-}<\mathfrak{a}_{\alpha_-}$, smaller than any action for type II. The continuation map $\sigma_\#(H_-, H_+)$ hence maps $FC^0(H_+)$ to $FC^0(H_-)$, which implies that $\sigma_\#(H_-, H_+)$ also maps $FC^+(H_+)=FC^{>\veps}(H_+)$ to $FC^+(H_-)=FC^{>\veps}(H_-)$.
\end{proof}
If $\partial_s H_{\alpha_s}\leq \frac{1}{2}\mathfrak{a}_{\alpha_-}$ does not hold, we can use the following trick due to Cieliebak and Frauenfelder, \cite{FraCie}. Fix a smooth monotone cutoff function $\beta\in C^\infty(\mathbb{R},[0,1])$ satisfying $\beta(s)=1$ for $s{\geq}1$ and $\beta(s)=0$ for $s{\leq}0$. Set $\mathfrak{a}_{min}:=\min \mathfrak{a}_{\alpha_s}=\min \bigcup_s spec(\alpha_s)$, which is positive as $R_{\alpha_s}$ has no zero for all $s$. Fix $N\in\mathbb{N}$ and consider for $0\leq j\leq N$ the Hamiltonians $H^j:=H_{\alpha_{j/N}}$ with homotopies between them
\[\alpha^j_s:= \alpha_{\frac{j+\beta(s)}{N}}\qquad\text{ and }\qquad H^j_s:=H_{a^j_s}, \quad 0\leq j\leq N{-}1.\]
Note that all $H^j_s, H^j$ are admissible, $H^0{=}H_{\alpha_-}, H^N{=}H_{\alpha_+}$ and $\partial_s H^j=\frac{\beta'(s)}{N}(\partial_s H)_{\alpha_s^j}$. As $\partial_s H_{\alpha_s}$ is bounded from above, we can assume for $N$ large that $\partial_s H_s^j \leq \frac{1}{2}\mathfrak{a}_{min}$ for each $0\leq j\leq N{-}1$. Hence, we obtain by Lem.\ \ref{lemactiondicho} that the continuation maps $\sigma_\ast(H^{j-1}, H^j):FH^0(H^j)\rightarrow FH^0(H^{j-1})$ and $FH^+(H^j)\rightarrow FH^+(H^{j-1})$ are well-defined. Then we set
\[\sigma_\ast(H_-, H_+):=\sigma_\ast(H^0,H^1)\circ\sigma_\ast(H^1,H^2)\circ...\circ\sigma_\ast(H^{N-1},H^N),\]
which maps $FH^0(H_+)$ into $FH^0(H_-)$ and $FH^+(H_+)$ into $FH^+(H_-)$.\pagebreak[1]\\
It follows from the homotopy of homotopies argument and the composition rule that this definition of $\sigma_\ast(H_-, H_+)$ does not depend on the chosen homotopy $H_{\alpha_s}$ nor on its decomposition $\{H^j_s\}$ into slower homotopies. The map $\sigma_\ast(H_-, H_+):FH^+(H_+)\rightarrow FH^+(H_-)$ is hence well-defined and we can define the positive Symplectic Homology as the direct limit over these maps: $\displaystyle SH^+(V)=\varinjlim FH^+(H)$.\\
Its invariance under Liouville isomorphisms follows in the same way as in Prop.\ \ref{propSHinv}. Moreover, we have the long exact sequence
\begin{equation}\label{eqLongExSeq}
 ...\rightarrow SH^0_\ast(V)\rightarrow SH_\ast(V)\rightarrow SH^+_\ast(V)\rightarrow SH^0_{\ast-1}(V)\rightarrow...,
\end{equation}
where $SH^0_\ast(V)\cong H_{n-\ast}(V,\partial V)$ holds by the usual argument (see \cite{Vit} or \cite{FraCieOan}, Lem.\ 2.1), as Hamiltonians adapted to a fixed contact form are a cofinal family for this direct limit.

\subsection{The Conley-Zehnder index}\label{secConZeh}
We $\mathbb{Z}$-grade the Symplectic Homology via the Conley-Zehnder index $\mu_{CZ}$. For simplicity, we restrict ourselves to the subcomplex of $FC_\ast(H,h)$ generated by the contractible 1-periodic orbits of $X_H$, which is no restriction if $V$ is simply connected. If we assume that $i_\ast:\pi_1(\partial V)\rightarrow\pi_1(V)$ induced by the inclusion is injective and that $\mu_{CZ}$ is well-defined on $V$, then this grading is independent of $V$.\\
To define $\mu_{CZ}(v)$ for a contractible 1-periodic Hamiltonian orbit $v$ choose a map $u$ from the unit disc $D\subset\mathbb{C}$ to $\Sigma$ such that $u(e^{2\pi it})=v(t)$, which is possible, as $v$ is contractible and $i_\ast$ injective. Then choose a symplectic trivialization ${\Phi:D{\times}\mathbb{R}^{2n}\rightarrow u^\ast TV}$ of $(u^\ast TV, u^\ast \omega)$, or equivalently a trivialization of $u^\ast\xi$, as the Liouville and Reeb vector field trivialize its complement. The linearization of the Hamiltonian flow $\varphi^t_{X_H}$ along $v$ via $\Phi$ defines a path $\Psi(t)$ in the group $Sp(2n)$ starting at $\mathbbm{1}$. The Maslov index of this path (as defined in \cite{RoSa}, \cite{Sal}) is $\mu_{CZ}(v)$.\\
Let $J_0$ denote the standard almost complex structure on $\mathbb{R}^{2n}$. Then every smooth path $\Psi:[a,b]\rightarrow Sp(2n)$ can be written in the form $\frac{d}{dt}\Psi(t)=J_0 S(t)\Psi(t)$, where $t\mapsto S(t)$ is a smooth path of symmetric matrices. The Maslov index has in particular the following properties (see \cite{RoSa} and \cite{Gutt1}):\\
 \textbf{(CZ0)} If $sign(S)=0$ everywhere, then $\mu_{CZ}(\Psi)=0$.\label{CZ0}\\
 \textbf{(CZ1)} If $\Psi:[0,T]\rightarrow Sp(2),\,\Psi(t)=e^{it}$, then $\displaystyle\mu_{CZ}(\Psi)=\left\lfloor\frac{T}{2\pi}\right\rfloor + \left\lceil\frac{T}{2\pi}\right\rceil.$\\
 \textbf{(product)} For $\Psi{\oplus}\Psi': [a,b]\rightarrow Sp(2n){\oplus} Sp(2n')\subset Sp\big(2(n{+}n')\big)$ holds \\$\mu_{CZ}(\Psi{\oplus}\Psi')=\mu_{CZ}(\Psi)+\mu_{CZ}(\Psi')$.\\
 \textbf{(naturality)} $\mu_{CZ}(\Phi\Psi\Phi^{-1})=\mu_{CZ}(\Psi)$, if $\Phi: [a,b]\rightarrow Sp(2n)$ is a contractible loop.\\
 \textbf{(zero)} If $\dim \ker\big(\Psi(t){-}\mathbbm{1}\big)=k$ is constant on $[a,b]$, then $\mu_{CZ}(\Psi)=0$.\\
 \textbf{(homotopy)} $\mu_{CZ}(\Psi_0)=\mu_{CZ}(\Psi_s)$ along homotopies $\Psi_s$ with fixed endpoints.\\
 \textbf{(catenation)} $\mu_{CZ}(\Psi|_{[a,b]})=\mu_{CZ}(\Psi|_{[a,c]})+\mu_{CZ}(\Psi|_{[c,b]})$ for any $a<c<b$.\\
If $k{\cdot} v$ denotes the $k$-fold iteration of a closed orbit $v$, then holds 
\begin{flalign}\label{eqIteration}
 \text{\textbf{(iterations formula)}}&\qquad \mu_{CZ}(k{\cdot} v) = k \cdot\Delta (v)+ R(k{\cdot} v),&
\end{flalign}
where $\Delta(v)$ is the \textit{mean index} of $v$ and $R(k{\cdot} v)$ is an error term that is bounded by $|R(k{\cdot} v)|\leq \frac{1}{2}\dim(\xi)=n{-}1$. This formula can be found in \cite{SaZeh}, Lem.\ 3.4 for the definition of the Conley-Zehnder index by Salamon and Zehnder. Note that this definition requires that $\Psi(0)=0$ and $\det(\Psi(T)-\mathbbm{1})\neq 0$. That the iterations formula also holds for the more general definition of $\mu_{CZ}$ due to Robbin and Salamon, \cite{RoSa}, can be seen as follows (cf.\ \cite{FauckThesis}, Lem.\ 60): 
Let $\Psi(t)\in Sp(2n)$ be the path of symplectic matrices associated to $v$ (and $k{\cdot} v$) with $\dot{\Psi}=J_0S\Psi$. Then for $\veps>0$ small holds
\begin{align*}
2\mu_{CZ}(k{\cdot} v)=2\mu_{CZ}(\Psi|_{[0,k\cdot T]})&\overset{(1)}{=}\mu_{CZ}(\Psi|_{[0,k\cdot T-\veps]})+\mu_{CZ}(\Psi|_{[0,k\cdot T+\veps]})\\
&\overset{(2)}{=}\Delta(\Psi|_{[0,k\cdot T-\veps]})+R(-\veps)+\Delta(\Psi|_{[0,k\cdot T+\veps]})+R(+\veps).
\end{align*}
Here, we used in (1) the definition from \cite{RoSa} and in (2) the definition from \cite{SaZeh}. With $\veps\rightarrow 0$, the iterations formula follows as $\Delta(\Psi|_{[0,t]})$ is continuous in $t$, $\Delta(\Psi|_{[0,kT]})=k\cdot\Delta(\Psi|_{[0,T]})$ and $|R(\pm\veps)|\leq n{-}1$ for all $\veps$ (see \cite{SaZeh}).\pagebreak[3]\\
The Conley-Zehnder index is also defined for Reeb orbits (via the Reeb flow instead of the Hamiltonian flow). Moreover, the two indices coincide for closed $X_H$-orbits $v$ that are reparametrizations of Reeb orbits.\\
A transversely non-degenerate orbit $v\in\mathcal{P}(H,h)$ is not graded by the Conley-Zehnder index alone, but by the following Morse-Bott index (see \cite{BourOan1} or \cite{FraCie})
\begin{equation}\label{eqmu}
 \mu(v):=\mu_{CZ}(v)+\mu_{Morse}(v)-{\textstyle\frac{1}{2}}\dim_v\mathcal{N}+{\textstyle\frac{1}{2}} sign \big(h''(e^r)\big),
\end{equation}
where $\mu_{Morse}(v)$ is the Morse index of $v$, $\dim_v\mathcal{N}$ the dimension of $\mathcal{N}$ near $v$ and $sign \big(h''(e^r)\big)$ the sign of $h''$ on the level $e^r$, where $v$ lives. The following lemma asserts that $\mu(v)$ is well-defined.
\begin{lemma}
 $\mu_{CZ}(v)$ is constant on connected components of $\mathcal{N}$.
\end{lemma}
\begin{proof}
 If $v_0, v_1$ are two 1-periodic orbits of $X_H$ in the same connected component of $\mathcal{N}$ and $u_0: D{\rightarrow}\mathbb{C}$ is such that $u_0(e^{2\pi it})=v_0(t)$, then we may construct $u_1:D{\rightarrow}\mathbb{C}$ such that $u_1(e^{2\pi it})=v_1(t),\; u_0(z)=u_1(\frac{1}{2}z)\;\forall\, z {\in} D$ and $v_s(t):=u_1(\frac{1}{2}(s{+}1)e^{2\pi it}), s\in[0,1]$ defines a path of 1-periodic $X_H$-orbits between $v_0$ and $v_1$. Choosing a symplectic trivialization $\Phi_0$ of $(u_0^\ast TV, u_0^\ast\omega)$, we may extend  this to a trivialization $\Phi_1$ of $(u_1^\ast TV, u_1^\ast\omega)$. Then, we find that $\Psi_1=\Phi_1^{-1}(v_1)\circ D\varphi^t_{X_H}(v_1)\circ \Phi_1(v_1)$ is homotopic with fixed endpoints to the catenation $\Lambda_1\ast\Psi_0\ast\Lambda_0$, where 
\begin{align*}
 \Psi_0(t)&=\Phi_0^{-1}(v_0(t))\circ D\varphi^t_{X_H}(v_0(0))\circ\Phi_0(v_0(0)),\\
 \Lambda_0(s)&=\Phi^{-1}_{1-s}(v_{1-s}(0))\circ D\varphi^0_{X_H}(v_{1-s}(0))\circ\Phi_{1-s}(v_{1-s}(0)),\\
 \Lambda_1(s)&=\Phi^{-1}_s(v_s(1))\circ D\varphi^1_{X_H}(v_s(0))\circ\Phi_s(v_s(0)).
\end{align*}
Here, $\Lambda_0\equiv\mathbbm{1}$ and $\dim\ker(\Lambda_1(s){-}\mathbbm{1})=\dim_v\mathcal{N}$. Using (zero), (homotopy) and (catenation), we find $\mu_{CZ}(\Psi_1)=\mu_{CZ}(\Lambda_1)+\mu_{CZ}(\Psi_0)+\mu_{CZ}(\Lambda_0)=\mu_{CZ}(\Psi_0)$.
\end{proof}

\subsection{Asymptotically finitely generated contact structures}\label{secAsmFini}
Let $H$ be a Hamiltonian adapted to a contact form $\alpha$ as in \ref{secPosSym}, i.e.\ $H$ is a Morse function $g$ inside $V_\alpha$ and cylindrically convex on $[-\veps,\infty){\times}\Sigma_\alpha$. Assuming that $\alpha$ satisfies (\ref{CondMB}), we find that the non-constant orbits of $X_H$ form a manifold $\mathcal{N}$ which agrees with the manifold formed by the closed Reeb orbits. If we choose a Morse function $h$ on $\mathcal{N}$, then $FC(H,h)$ is generated by critical points of $h$ and $g$, while $FC^+(H,h)$ is generate by critical points of $h$. If $crit_k(g)$ denotes the critical points of $g$ with Morse-Bott index $k$, then $FC_k(H,h)$ is generated by at most $\# crit_k(h)+\#crit_k(g)$ many critical points and $FC_k^+(H,h)$ by at most $\# crit_k(h)$ both for any $H$ of the above form and hence
 \begin{align*}
  \rk FH_k(H,h)&\leq \# crit_k(h)+\# crit_k(g) &&\Rightarrow& \rk SH_k(V)&\leq \# crit_k(h)+\# crit_k(g),\\
  \rk FH_k^+(H,h)&\leq \# crit_k(h)&&\Rightarrow& \rk SH_k(V)&\leq \# crit_k(h).
 \end{align*}
 Unfortunately, $\# crit_k(h)$ may be infinite. In particular the result of contact surgery may have infinitely many new orbits. However, we can assume that we have only finitely many orbits up to a certain length which motivates the following definition.
 \begin{defn}\label{defnAsympFinite}~\\
  Let $(\Sigma,\xi)$ be a compact contact manifold and let $\alpha$ be a contact form for $\xi$. We say that $\xi$ is an \textbf{asymptotically finitely generated contact structure in degree $k$ with bound $b_k(\xi)$}, if there exist sequences of smooth functions $f_l:\Sigma\rightarrow\mathbb{R}$ and real numbers $(\mathfrak{a}_l)\subset\mathbb{R}$ such that
  \begin{itemize}
   \item  $\log \mathfrak{a}_l-f_l(p)\leq \log \mathfrak{a}_{l+1}-f_{l+1}(p)$ and $\displaystyle\lim_{l\rightarrow\infty}\big(\log\mathfrak{a}_l-f_l(p)\big)=\infty$ for all $p\in\Sigma$,
   \item all contractible Reeb orbits of length at most $\mathfrak{a}_l$ of the contact form $\alpha_l:=e^{f_l}\cdot\alpha$ are transversely non-degenerate and among these orbits at most $b_k(\xi)$ have Morse-Bott index $k$ for some choice of a Morse function $h$ on $\mathcal{N}$.
  \end{itemize}
 \end{defn}\pagebreak[3]
\begin{remark}~\label{remarksafg}
 \begin{itemize}
  \item We do not require that the $\alpha_l$ are distinct. In particular $\alpha=\alpha_l\; \forall\, l$ is possible, if for $\alpha$ itself all closed Reeb orbits are transversely non-degenerate and only finitely many have Morse-Bott index $k$. In this situation, $b_k(\xi)$ can be chosen to equal the number of closed Reeb orbits of $\alpha$ having Morse-Bott index $k$.
  \item By slightly increasing the $\mathfrak{a}_l$, we can always assume that $\mathfrak{a}_l\not\in spec\big(\alpha_l\big)$.
  \item The first condition is equivalent to $\big(\log\mathfrak{a}_l-f_l\big)\rightarrow+\infty$ uniformly in $l$. Moreover, we can assume without loss of generality that $\mathfrak{a}_{l+1}\geq \mathfrak{a}_l$ and $f_{l+1}\leq f_l\leq 0$. Indeed, if this does not hold, set $M_1:=\max_p f_1(p),\;M_{l+1}:=\max_p \big(f_{l+1}(p){-}f_l(p)\big)$ and consider the new sequences
  \[\hat{f}_l:=f_l\,{-}M_l\,{-}...{-}M_1\qquad\text{ and }\qquad \hat{\mathfrak{a}}_l:=\mathfrak{a}_l\cdot\exp\big({-}M_l\,{-}...{-}M_1\big).\]
  Then $\hat{f}_1\leq 0$ and $\hat{f}_{l+1}=f_{l+1}\,{-}M_{l+1}\,{-}M_l\,{-}...{-}M_1\leq f_l\,{-}M_l\,{-}...{-}M_1=\hat{f}_l$. Moreover, $\log \hat{\mathfrak{a}}_l{-}\hat{f}_l=\log \mathfrak{a}_l{-}f_l$, hence $\big(\log \hat{\mathfrak{a}}_l{-}\hat{f}_l\big)\rightarrow\infty$. As $M_{l+1}=f_{l+1}(p_0){-}f_l(p_0)$ for at least one $p_0\in\Sigma$, we have $\hat{f}_{l+1}(p_0)=\hat{f}_l(p_0)$ and hence by the first a.f.g.\ condition:
  \[\log \hat{\mathfrak{a}}_l-\hat{f}_l(p_0)\leq \log \hat{\mathfrak{a}}_{l+1}-\hat{f}_{l+1}(p_0)\quad\Leftrightarrow\quad \hat{\mathfrak{a}}_l\leq \hat{\mathfrak{a}}_{l+1}.\]
  Finally, $\hat{\alpha}_l=e^{-(M_l+...+M_1)}\cdot \alpha_l$, so the new Reeb vector field is ${\hat{R}_l=e^{M_l+...+M_1}{\cdot} R_l}$. Thus all closed Reeb orbits of $R_l$ of length at most $\mathfrak{a}_l$ are in 1-1 correspondence to the closed Reeb orbits of $\hat{R}_l$ of length at most $\hat{\mathfrak{a}}_l$.
  \item The existence of a.f.g.\ contact structures is guaranteed by Prop.\ \ref{PropafgExists1} - \ref{propcofinalHamilt}.
  \item Recently, similar definitions have been introduced: convenient dynamics in \cite{KwKo} and asymptotical dynamical convexity in \cite{Lazarev}. All three definitions work with sequences of contact forms with ``nice'' closed Reeb orbits below a certain length. However, our definition is a generalization of the other two as index positivity or negativity implies a.f.g.\ for $k\leq 0$ resp. $k\geq 0$.
 \end{itemize}
\end{remark}
\begin{proof}[\textnormal{\textbf{Proof of Proposition \ref{PropAsympfiniteGene}}}]~\\
 Let $(\Sigma,\xi)$ be a closed contact manifold with exact filling $(V,\lambda)$ and well-defined Conley-Zehnder index. If $\xi$ is a.f.g.\ in degree $k$ with bound $b_k(\xi)$, we have to show that $\rk SH_k^+(V)$ can be bounded in terms of $b_k(\xi)$.\\
 Set $\alpha:=\lambda|_{T\Sigma}$ and let $\alpha_l=e^{f_l}{\cdot} \alpha$ and $(\mathfrak{a}_l)\subset\mathbb{R}, \mathfrak{a}_l\not \in spec(\alpha_l),$ be sequences as in Defn.\ \ref{defnAsympFinite}, showing that $\xi$ is a.f.g.. Let $\Sigma_l:=\Sigma_{\alpha_l}$ and $V_l:=V_{\alpha_l}$ be associated to $\alpha_l$ as in \ref{secsetup}. Let $H_l$ be adapted to $\alpha_l$ as in \ref{secPosSym}, i.e.\ $H$ is a Morse function $g$ inside $V_l$ and is on $[-\veps,\infty){\times}\Sigma_l$ of the form $H_l(r,p)=h_l(e^r)$ with $h''_l\geq 0$ and $h_l(e^r)=\mathfrak{a}_le^r-\mathfrak{a}_l$ on $[0,\infty){\times}\Sigma_l$. If we express $H_l$ in the fixed coordinates $[0,\infty){\times}\Sigma$, it takes for $r$ sufficiently large the form
 \[H_l(r,p)=\mathfrak{a}_le^{r-f_l(p)}-\mathfrak{a}_l.\]
 Apparently, each $H_l$ is an admissible Hamiltonian. As $\log \mathfrak{a}_{l+1}-f_{l+1}\geq \log \mathfrak{a}_l-f_l$, we find that $H_{l+1}\prec H_l$ for all $l$ and $\lim_{l\rightarrow\infty} (\log \mathfrak{a}_l-f_l)=\infty$ implies that $(H_l)$ is cofinal in $Ad(V)$.\pagebreak[1]\\
 Hence we have that
 \[\lim_{l\rightarrow\infty} FH_\ast(H_l)=SH_\ast(V)\quad\text{ and }\quad \lim_{l\rightarrow\infty} FH^+_\ast(H_l)=SH^+_\ast(V),\]
 where we use the definition of $SH^+_\ast(V)$ as described in \ref{secPosSym}. As shown in \ref{secPosSym}, $FH^+(H_l)$ is generated by orbits of type II, i.e.\ by the closed Reeb orbits of $\alpha_l$. As there are for each $l$ at most $b_k(\xi)$ closed Reeb orbits of $\alpha_l$ with Morse-Bott index $k$, we find that
 \[ \rk FH^+_k(H_l)\leq \rk FC^+_k(H_l)\leq b_k(\xi) \qquad \overset{\varinjlim}{\Longrightarrow}\qquad \rk SH^+_k(V) \leq b_k(\xi).\]
 With the long exact sequence (\ref{eqLongExSeq}): $...\rightarrow SH_k^0(V)\rightarrow SH_k(V)\rightarrow SH_k^+(V)\rightarrow...$, where $SH^0_k(V)\cong H_{n-k}(V,\partial V)$, we obtain the estimate
 \[\rk SH_k(V)\leq \rk SH_k^0(V)+\rk SH_k^+(V)\leq \rk H_{n-k}(V,\partial V)+ b_k(\xi).\qedhere\]
\end{proof}

\subsection{The transfer morphisms}\label{sectransfer}
Following Viterbo, \cite{Vit}, we show in this section for Liouville subdomains $W\subset V$ the existence of the transfer maps \[\pi_\ast(W,V):SH_\ast(V)\rightarrow SH_\ast(W) \quad\text{ and }\quad \pi^\ast(W,V):SH^\ast(W)\rightarrow SH^\ast(V).\]
As shown in \ref{sectrunc}, we have maps $SH_\ast(V)\rightarrow  SH_\ast^{\geq0}(W{\subset}V)$ and $SH^\ast_{\geq0}(W{\subset}V)\rightarrow SH^\ast(V)$ coming from the projection in action filtered complexes. We will show in Prop.\ \ref{proptrans} and Cor.\ \ref{transfer} the identities $SH_\ast^{\geq0}(W{\subset}V)= SH_\ast(W)$ and $SH^\ast_{\geq0}(W{\subset}V)= SH^\ast(W)$, so that these projections are the transfer maps. Prop.\ \ref{proptrans} is based on ideas by Viterbo, \cite{Vit}, however the proof follows McLean, \cite{McLeanDis}. We include it here for completeness and as the use of the No-Escape Lemma provides simplifications.
\begin{proposition}[McLean,\cite{McLeanDis}]\label{proptrans}~\\
 Let $W\subset V$ be a Liouville subdomain. Then, there exists an increasing cofinal sequence $(H_n)\subset Ad^{\geq0}(W{\subset}V)$ and a sequence of decreasing admissible homotopies $(H_{n,n+1})$ between them such that:
 \begin{enumerate}
  \item $H_n|_W, \; H_{n,n+1}|_W$ are admissible Hamiltonians{\large/}homotopies on $(W,\omega)$,
  \item all 1-periodic orbits of $X_{H_n}$ in $W$ have positive action and all 1-periodic orbits of $X_{H_n}$ in $\widehat{V}\setminus W$ have negative action,
  \item all $\mathcal{A}^H$-gradient trajectories of $H_n$ or $H_{n,n+1}$ connecting 1-periodic orbits in $W$ are entirely contained in $W$ for all admissible $J$ that are of contact type near $\partial W$.
 \end{enumerate}
\end{proposition}
\begin{proof}
 It will be convenient to use $z=e^r$ rather than $r$ for the radial coordinate in the completions $(\widehat{W},\widehat{\omega})$ and $(\widehat{V},\widehat{\omega})$. We embed $\widehat{W}$ into $\widehat{V}$ using the flow of the Liouville vector field $Y$. The cylindrical end $[1,\infty)\times\partial W$ is then a subset of $\widehat{V}$. The radial coordinate is denoted $z_W$ on $\partial W\times(0,\infty)$ and  $z_V$ on $\partial V\times(0,\infty)$. We choose a constant $P$ such that $\{z_W{\leq} 1\}\subset \{z_V{\leq} P\}$, which implies $\{z_W{\leq} C\}\subset\{z_V{\leq} C{\cdot} P\}$ for any $C>0$. Let $\alpha_W:=\lambda|_{T\partial W}$, $\alpha_V:=\lambda|_{T\partial V}$ and assume that $(\partial W,\alpha_W)$ and $(\partial V,\alpha_V)$ satisfy (\ref{CondMB}).\\
 To construct $H_n$ choose an increasing sequence $(\mathfrak{a}_n)\subset\mathbb{R}^+$ with $\mathfrak{a}_n\rightarrow\infty$ and 
 \begin{align*}
  (\mathfrak{a}_n)&\not\in \Big(spec(\partial W,\alpha_W)\cup 4P{\cdot} spec(\partial V,\alpha_V)\Big)\qquad\text{ for all $n$}.\\
  \text{Let }\qquad\mu_n:&=dist\big(\mathfrak{a}_n,spec(\partial W,\alpha_W)\big)=\min_{a\in spec(\partial W,\alpha_W)}|\mathfrak{a}_n-a|>0,
 \end{align*}
 let $(\veps_n)$ be a decreasing sequence with $\veps_n\rightarrow 0$ and $\veps_1$ sufficiently small. Finally, let $Z_n$ be an increasing sequence such that: $\qquad Z_n>\frac{\mathfrak{a}_n}{\mu_n}$ and $Z_n>2$.\hfill $(\ast)$\\
To ease notation, we write only $Z,\mathfrak{a},\mu,\veps$, whenever there is no danger of confusion.\medskip\\
 \begin{figure}[htb]
\centering
 \resizebox{13cm}{!}{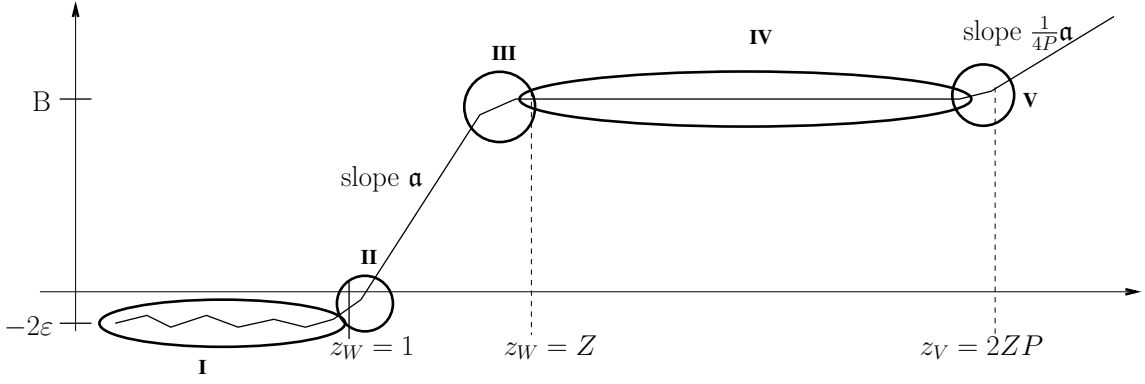}
 \caption{\label{fig7}The Hamiltonian $H_n$ and the areas of the five orbit types}
\end{figure}
 Next, we describe the Hamiltonian $H_n$ (see figure \ref{fig7} for a schematic illustration):\\
 Inside $W\setminus\big([1{-}\veps,1){\times}\partial W\big)$ let $H_n$ be a $C^2$-small Morse function with $-2\veps<H_n< \veps$. On $[1{-}\veps,Z]{\times}\partial W$ let it be of the form $H_n(z_W,p)=g(z_W)$ with $g(1)=-\veps,\;0\leq g'(z_W)\leq \mathfrak{a}$ and $g'(z_W)\equiv \mathfrak{a}$ for $1\leq z_W \leq Z{-}\veps$. On $[Z,2Z]{\times} \partial W$ let $H_n\equiv B_n\approx \mathfrak{a}_n{\cdot}(Z_n{-}1)$ be constant.\\
 On $[1,\infty){\times}\partial V$ keep $H_n$ constant until we reach the hypersurface $\{z_V=2ZP{-}\veps\}$ (recall $\{z_W{\leq} 2Z\}\subset\{z_V{\leq} 2ZP\}$). Then let $H_n$ be of the form $H_n(z_V,p)=f(z_V)$ for $z_V\geq 2ZP{-}\veps$ with $0\leq f'(z_V)\leq \frac{1}{4P}\mathfrak{a}$ and $f'(z_V)\equiv\frac{1}{4P}\mathfrak{a}$ for $z_V\geq 2ZP$.\pagebreak[3]\medskip\\
 Recall that the action of an $X_H$-orbit on a fixed $z$-level is $h'(z)\cdot z-h(z)$, if $H(z,p)=h(z)$. Hence we distinguish five types of 1-periodic orbits of $X_H$:
 \begin{itemize}
  \item[\textbf{\textit{I}} :] critical points inside $W$ of action $\geq\veps$ (as $H_n\leq-\veps$ and $C^2$-small inside $W$)
  \item[\textbf{\textit{II}} :] non-constant orbits near $z_W=1$ of action $\approx g'(z)>0$
  \item[\textbf{\textit{III}} :] non-constant orbits on $z_W=c$ for $c$ near $Z$ of action\\ $\Big.\approx g'(c)\cdot c-B<(\mathfrak{a}{-}\mu)\cdot Z-B\approx -\mu\cdot Z+\mathfrak{a}\overset{(\ast)}{<}0$ 
  \item[\textbf{\textit{IV}} :] critical points in $Z<z_W , z_V<2ZP-\veps$ of action $-B<0$
  \item[\textbf{\textit{V}} :] non-constant orbits on $z_V=c$ for $c$ near $2ZP$ of action\\ $\approx f'(c){\cdot} c-B<\frac{1}{4P}\mathfrak{a}{\cdot} 2ZP-B\approx\frac{1}{2}\mathfrak{a}Z-\mathfrak{a}(Z{-}1)=\frac{1}{2}\mathfrak{a}(2{-}Z)\overset{(\ast)}{<}0.$
 \end{itemize}
Hence, $(H_n)$ satisfies claim (2) of the proposition. For claim (1), note that $H_n|_W<H_{n+1}|_W$ (as $-2\veps_{n+1}>-2\veps_n$ and $\mathfrak{a}_{n+1}>\mathfrak{a}_n$) and that the linear extensions of $H_n|_W$ to $\widehat{W}$ form a cofinal sequence of admissible Hamiltonians on $W$.
\begin{figure}[ht]
\centering
 \resizebox{10cm}{!}{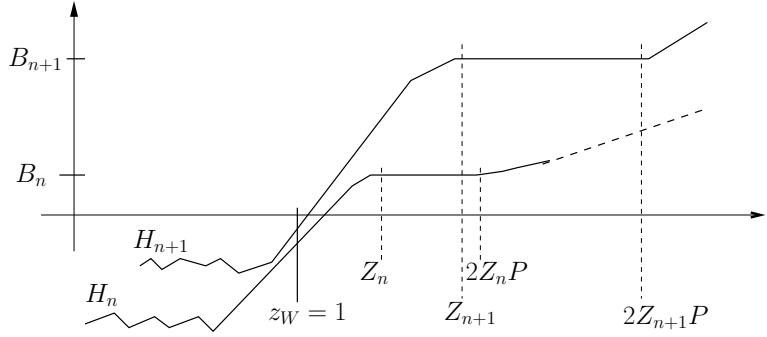}
 \caption{\label{fighomo}Two Hamiltonians $H_n$ and $H_{n+1}$}
\end{figure}\\
In fact $H_n<H_{n+1}$ holds globally: As $\frac{1}{4P}\mathfrak{a}_n<\frac{1}{4P}\mathfrak{a}_{n+1}$, we have $H_n<H_{n+1}$ for $z_V$ large enough. Yet, the most dangerous area is around $z_V=2Z_{n+1}P$ (see figure \ref{fighomo}), where we estimate:
\begin{align*}
 H_n(2Z_{n+1}P)=h_n(2Z_{n+1}P)&=\textstyle\frac{1}{4P}\mathfrak{a}_n{\cdot} 2Z_{n+1}P+B_n-\frac{1}{2}\mathfrak{a}_nZ_n\\
 &\approx \textstyle\frac{1}{2}\mathfrak{a}_nZ_{n+1}\,{+}\mathfrak{a}_n(Z_n{-}1)\,{-}\frac{1}{2}\mathfrak{a}_nZ_n\\
 &<\textstyle\frac{1}{2}\mathfrak{a}_{n+1}Z_{n+1}+\mathfrak{a}_{n+1}(\frac{Z_{n+1}}{2}\,{-}1)\\ &=\mathfrak{a}_{n+1}(Z_{n+1}{-}1)\;=\;B_{n+1}\;\approx\; H_{n+1}(2Z_{n+1}P).
\end{align*}
As $H_n<H_{n+1}$, we find a decreasing homotopy $H_{n,n+1}$ from $H_{n+1}$ to $H_n$ satisfying:
\begin{align*}
 H_{n,n+1}(z_W,p)&=\mathfrak{a}_s^W{\cdot} z_W \, {+} \mathfrak{b}^W_s & &\text{near }z_W{=}1,& \partial_s \big(H_{n,n+1}{-}\mathfrak{b}^W_s\big)&\leq 0 &&\text{for $z_W\,{\geq}\, 1$}\\
 H_{n,n+1}(z_V,p)&=\mathfrak{a}_s^V{\cdot} z_V \; {+} \mathfrak{b}^V_s & &\text{and}& \partial_s \big(H_{n,n+1}{-}\mathfrak{b}^V_s\big)&\leq 0 &&\text{for $z_V\,{\geq}\, 2Z_{n+1}P$}\mspace{-1mu}.
\end{align*}
If $J$ is admissible and of contact type near $z_W=1$, then the No-Escape Lemma (Cor.\ \ref{Cornoesc}) implies that all $\mathcal{A}^H$-gradient trajectories for $H_n$ or $H_{n,n+1}$ connecting 1-periodic orbits inside $W$ stay inside $W$ for all time. This shows claim (3).
\end{proof}
\begin{corollary}\label{transfer} $\quad\displaystyle SH^{\geq 0}_\ast(W{\subset}V)\simeq SH_\ast(W)\;\text{ and }\; SH^\ast_{\geq0}(W{\subset}V)\simeq SH^\ast(W).$
\end{corollary}
\begin{proof}
 We only prove the corollary for homology, cohomology being completely analog. Take the sequence of Hamiltonians $(H_n)$ constructed in Prop.\ \ref{proptrans}. Clearly it is cofinal and $(H_n)\subset Ad^{\geq0}(W{\subset}V)$, as 1-periodic orbits with positive action are either isolated critical points inside $W$ (as $H$ is Morse and $C^2$-small there) or isolated Reeb-orbits near $z_W=1$ -- in both cases non-degenerate. Hence we have
 \[ SH^{\geq0}_\ast(W{\subset}V)=\lim_{\longrightarrow} FH^{\geq0}_\ast(H_n).\]
 Let $\tilde{H}_n\in Ad(W)$ be the linear extension of $H_n|_W$ to $\widehat{W}$ with slope $\mathfrak{a}_n$. Then we have obviously $FC^{\geq0}_\ast(H_n)=FC_\ast(\tilde{H}_n)$. As any $\mathcal{A}^H$-gradient trajectory connecting 1-periodic orbits in $W$ stays in $W$, the boundary operators on $FC_\ast^{\geq 0}(H_n)$ and $FC_\ast(\tilde{H}_n)$ coincide and we have $FH^{\geq0}_\ast(H_n)=FH_\ast(\tilde{H}_n)$. As the $\mathcal{A}^H$-gradient trajectories for the homotopies $H_{n,n+1}$ stay inside $W$, the continuation maps $\sigma(H_{n+1},H_n): FH^{\geq0}_\ast(H_n)\rightarrow FH^{\geq0}_\ast(H_{n+1})$ coincide with the continuation maps $\sigma(\tilde{H}_{n+1},\tilde{H}_n): FH_\ast(\tilde{H}_n)\rightarrow FH_\ast(\tilde{H}_{n+1})$. Hence:
 \[SH^{\geq0}_\ast(W{\subset}V)= \varinjlim FH^{\geq0}_\ast(H_n) = \varinjlim FH_\ast(\tilde{H}_n)=SH_\ast(W).\]
\end{proof}
In the literature, there is a second description of the transfer maps which goes as follows. Let $(H_n)$ be the Hamiltonians described in Prop.\ \ref{proptrans} and let $(K_n)$ be the following sequence of Hamiltonians: Inside $V$ we require that $K_n$ is a $C^2$-small Morse function such that $K_n|_W\leq H_n|_W$ and on $[1{-}\veps,\infty)\times\partial V$ let it be of the form $K_n(z_V,p)=f(z_V)$ with $0\leq f'(z_V)\leq\frac{\mathfrak{a}_n}{4P}$, where $\mathfrak{a}_n$ and $P$ are as in Prop.\ \ref{proptrans} (see figure \ref{figHam2nd}).
\begin{figure}[ht]
\centering
 \resizebox{13cm}{!}{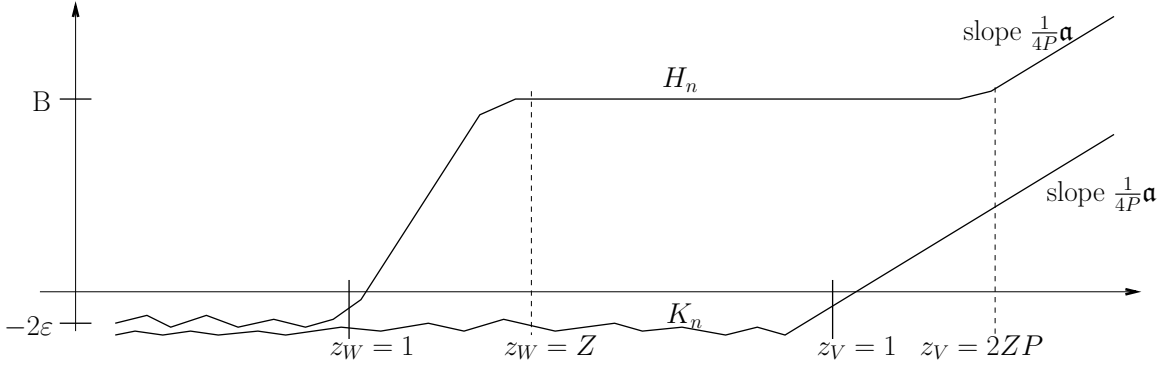}
 \caption{\label{figHam2nd}The two Hamiltonians $H_n$ and $K_n$}
\end{figure}
We find that $K_n\leq H_n$ everywhere, as in particular
\[K(2ZP)\approx 2ZP{\cdot}\frac{\mathfrak{a}}{4P}-\frac{\mathfrak{a}}{4P}=\mathfrak{a}\Big(Z-\Big(\frac{Z}{2}{+}\frac{1}{4P}\Big)\Big)\overset{(Z>2)}{<}\mathfrak{a}(Z{-}1)\approx B=H(2ZP).\]
Hence, we can find an everywhere increasing homotopy between $K_n$ and $H_n$, defining a continuation map $\sigma_\ast(K_n,H_n):FH_\ast(K_n)\rightarrow FH_\ast(H_n)$ which respects action filtration. The second version of the transfer map $\tilde{\pi}_\ast(W,V)$ is the limit of these continuation maps:
\[\begin{xy}\xymatrix{\tilde{\pi}_\ast(W,V)\;:\; SH_\ast(V)\overset{(\ast)}{=}\varinjlim FH^{\geq 0}_\ast(K_n)\ar[rr]^{\qquad\sigma_\ast(K_n,H_n)} && \varinjlim FH^{\geq 0}_\ast(H_n)=SH_\ast(W)}\end{xy}.\]
Here, the identity $(\ast)$ is due to the fact that all 1-periodic orbits of $K_n$ have positive action. The advantage of this transfer map is precisely that it respects action filtration. Hence it defines a map $\tilde{\pi}_\ast(W,V): SH^+_\ast(V)\rightarrow SH_\ast^+(W)$. However, the advantage of the first definition of $\pi_\ast(W,V)$ is that it fits into the following exact triangle, where the third group has a geometric meaning:
\[\begin{xy}\xymatrix{SH_\ast(V)=SH_\ast(W{\subset} V)\ar[rr]^{\pi_\ast(W,V)} && SH_\ast^{\geq 0}(W{\subset} V)=SH_\ast(W) \ar[dl]^{[-1]}\\ & \ar[ul] SH_\ast^{<0}(W{\subset} V)}\end{xy}.\]
\begin{proposition}\label{prop2transmaps}
 $\pi_\ast(W,V)$ and $\tilde{\pi}_\ast(W,V)$ are the same map: $SH_\ast(V)\rightarrow SH_\ast(W)$. The same holds true for the corresponding maps on cohomology.
\end{proposition}
\begin{proof}
 We can slightly modify the construction of the $H_n$, such that its 1-periodic orbits are all transversely non-degenerate and $H_n\in Ad(V)$. We just have to require on $Z\leq z_W,z_V\leq 2ZP\,{-}\veps$ that $H_n{-} B$ is a $C^2$-small Morse function instead of $H_n\equiv B$ being constant there. Then we have that $FH_\ast(H_n)$ is well-defined.\medskip\\
 \emph{\underline{Claim} :} $FH_\ast(H_n)\cong FH_\ast(K_n)$.\\
 \emph{Proof:} Note that $K_n$ and $H_n$ differ on $[R,\infty)\times\partial V$ only by a constant for $R$ large. Hence we can find a homotopy $H_s$ between them such that $H_s$ and $H_{-s}$ are both admissible. This provides continuation maps $\sigma_\ast(K_n,H_n): FH_\ast(K_n)\rightarrow FH_\ast(H_n)$ and \linebreak $\sigma_\ast(H_n,K_n): FH_\ast(H_n)\rightarrow FH_\ast(K_n)$. They satisfy $\sigma_\ast(H_n,K_n)\circ\sigma_\ast(K_n,H_n)=id_{FH(H_n)}$ and $\sigma_\ast(K_n,H_n)\circ\sigma_\ast(H_n,K_n)=id_{FH(K_n)}$, which implies that they are isomorphisms.\medskip\\
 As the homotopy from $K_n$ to $H_n$ can be chosen everywhere increasing, we find that $\sigma_\ast(K_n,H_n)$ respects action filtration and we have the following commutative diagram\pagebreak[1]
 \[\begin{xy}\xymatrix{FH_\ast(H_n)\ar[rr]^{\pi^{H_n}} && FH_\ast^{\geq0}(H_n)\\
  FH_\ast(K_n) \ar[u]^\cong_{\sigma_\ast(K_n,H_n)} \ar[rr]^{\pi^{K_n}}_\cong && FH^{\geq 0}(K_n). \ar[u]_{\sigma_\ast(K_n,H_n)}}\end{xy}\]
 Applying the direct limit yields again a commutative diagram, where isomorphisms are taken to isomorphisms:
 \[\begin{xy}\xymatrix{&\varinjlim FH_\ast(H_n)\ar[rr]^{\pi_\ast(W,V)} && \varinjlim FH_\ast^{\geq0}(H_n)\ar@{=}[r] & SH_\ast(W)\\
  SH_\ast(V)\ar@{=}[r] &\varinjlim FH_\ast(K_n) \ar[u]^\cong \ar[rr]_{\cong} &&\varinjlim FH^{\geq 0}(K_n). \ar[u]_{\tilde{\pi}_\ast(W,V)}}\end{xy}\]
  Hence, we find that $\pi_\ast(W,V)$ and $\tilde{\pi}_\ast(W,V)$ coincide on homology. The same arguments work also for cohomology, since $\varprojlim$ preserves isomorphisms, as it is left exact.
\end{proof}

\section{Contact surgery and handle attaching}\label{secsur}
In this section, we first describe the general construction for contact surgery by attaching a symplectic handle $\mathcal{H}$ to the symplectization of a contact manifold. Then, we describe explicitly symplectic handles as subsets of $\mathbb{R}^{2n}$ given by the intersection of two sublevel sets $\{\psi {<}{-}1\}\cap\{\phi{>}{-}1\}$ for functions $\phi$ and $\psi$ on $\mathbb{R}^{2n}$. Subsequently, we describe how to extend an admissible Hamiltonian over the handle to a new admissible Hamiltonian with only few new 1-periodic Hamiltonian orbit. The proofs of the Invariance Theorem (Thm.\ \ref{theoinvsur}) and of the invariance of a.f.g.\ under subcritical surgery (Prop.\ \ref{propcofinalHamilt}) conclude this section.

\subsection{Surgery along isotropic spheres}\label{secconsur}
\begin{figure}[htb]
\centering 
 \resizebox{6.5cm}{!}{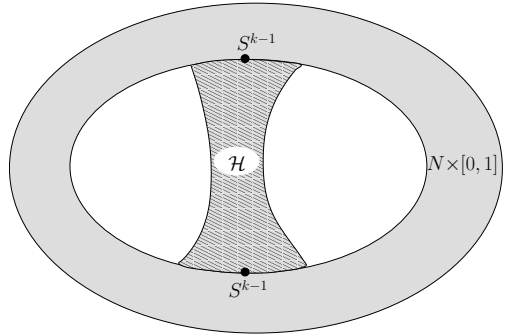}
 \caption{\label{fig2} $N\times[0,1]$ with handle attached}
\end{figure}

We briefly recall the contact surgery construction due to Weinstein, \cite{Wein}. Consider an isotropic sphere $S^{k-1}$ in a $(2n{-}1)$-dimensional contact manifold $(N,\xi)$. The 2-form $\omega=d\alpha$ for a contact form $\alpha$ for $\xi$ defines a natural conformal symplectic structure on $\xi$. Let $\perp_\omega$ be the $\omega$-orthogonal on $\xi$. As $S$ is isotropic, $TS\subset TS^{\perp_\omega}$. So, the normal bundle of $S$ in $N$ decomposes as
\[ TN/ TS= TN/ \xi \,\oplus\, \xi/(TS)^{\perp_\omega} \oplus (TS)^{\perp_\omega}/ TS.\]
The Reeb field $R$ trivializes $TN/\xi$, while $\xi/(TS)^{\perp_\omega}$ is canonically isomorphic to $T^\ast S$ via $v\mapsto \iota_v\omega$. The conformal symplectic normal bundle $CSN(S):=(TS)^{\perp_\omega}/TS$ carries a natural conformal symplectic structure induced by $\omega$. As $S$ is a sphere, the embedding $S^{k-1}\subset\mathbb{R}^k$ provides a natural trivialization of $\mathbb{R}R\oplus T^\ast S$. This trivialization together with a conformal symplectic trivialization of $CNS(S)$ specifies a framing for $S$ in $N$.\\ 
Following \cite{Wein}, an isotropic setup is a quintuple $(P,\omega,Y,\Sigma,S)$, where $(P,\omega)$ is an exact symplectic manifold, $Y$ a Liouville vector field for $\omega$, $\Sigma$ a hypersurface transverse to $Y$ and $S$ an isotropic submanifold of $\Sigma$. Isotropic setups satisfy the following variant of the neighborhood theorem for isotropic manifolds:
\begin{proposition}[\textbf{Weinstein},\cite{Wein}] \label{X}
 Let $(P_\pm,\omega_\pm,Y_\pm,\Sigma_\pm,S_\pm)$ be two isotropic setups. Given a diffeomorphism $\phi:S_-\rightarrow S_+$ covered by an isomorphism $CSN(S_-)\mspace{-5mu}\rightarrow \mspace{-3mu}CSN(S_+)$, there exist neighborhoods $U_\pm$ of $S_\pm$ in $P_\pm$ and an extension of $\phi$ to an isomorphism of isotropic setups
 \[\phi : (U_-,\omega_-,Y_-,\Sigma_-\cap U_-,S_-) \rightarrow (U_+,\omega_+,Y_+,\Sigma_+\cap U_+,S_+).\]
\end{proposition}
Contact surgery along an isotropic sphere is now defined as follows. Let $\mathcal{H}\approx D^k{\times} D^{2n-k}$ be a symplectic handle (see \ref{handle}) and let $S^{k-1}$ be an isotropic sphere in $(N,\xi)$. Given a standard framing for $S$, Prop.\ \ref{X} allows us to glue the (lower) boundary $S^k{\times} D^{2n-k}$ of $\mathcal{H}$ to the symplectization $[0,1]\times N$ along the boundary part $U_+\cap[0,1]{\times}N$ of a tubular neighborhood $U_+$ of $\{1\}{\times}S$ (see Figure \ref{fig2}). We obtain an exact symplectic manifold $P:=[0,1]{\times}N\cup_{S}\mathcal{H}$ with two boundary components  $\partial^-P :=\{0\}{\times}N$ and $\partial^+P$. Both components are contact and $\partial^+P$ is obtained from $N$ by surgery along $S$.
\vspace{-0.2cm}

\vspace{-0.4cm}\subsection{Symplectic handles}\label{handle}
\begin{figure}[htb]
\centering
 \resizebox{10cm}{!}{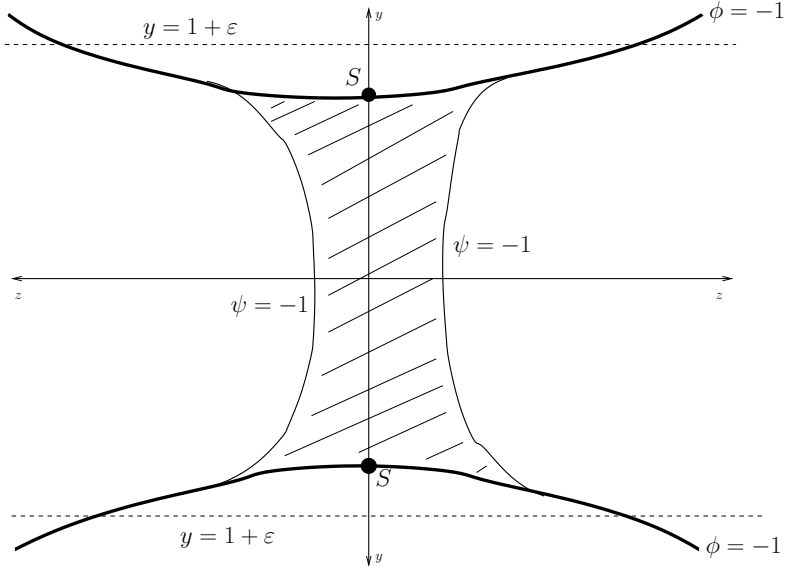}
 \\\caption{\label{fig3} The handle $\mathcal{H}$}
\end{figure}

We consider $\mathbb{R}^{2n}$ with symplectic coordinates $(q,p)=(q_1,p_1,...\,,q_n,p_n)$ and the following Weinstein structure (cf.\ \cite{Wein}):
\begin{align*}
\lambda &:= \textstyle\sum_{j=1}^k\left(2q_jdp_j + p_jdq_j\right)+\sum_{j=k+1}^n\frac{1}{2}\left(q_jdp_j -p_jdq_j\right),\; \omega\,{:=}\,d\lambda \,{=}\textstyle\sum_{j=1}^n dq_j\wedge dp_j,\\
 Y &:= \textstyle\sum_{j=1}^k \left(2q_j\partial_{q_j} - p_j\partial_{p_j}\right) + \sum_{j=k+1}^n \frac{1}{2}\left(q_j\partial_{q_j} + p_j\partial_{p_j}\right),\\
  \phi&:=\textstyle\sum_{j=1}^k\left( q_j^2-\frac{1}{2}p_j^2\right) + \sum_{j=k+1}^n \frac{1}{2}\left(q_j^2+p_j^2\right).
\end{align*}
Note that $Y$ is the Liouville vector field for $\lambda$, as $\iota_Y\omega = \lambda$.\pagebreak[1]\\
For convenience, we introduce furthermore functions $x,y,z:\mathbb{R}^{2n}\rightarrow\mathbb{R}$ with Hamiltonian vector fields $X_x, X_y, X_z$:
\begin{align*}
 x&=x(q,p)=  \sum_{j=1}^k q_j^2, & y&=y(q,p)= \sum_{j=1}^k \frac{1}{2} p_j^2, & z&=z(q,p)=  \sum_{j=k+1}^n \frac{1}{2}(q_j^2+p_j^2),\\
 X_x &= \sum_{j=1}^k 2q_j\partial_{p_j},&
 X_y &= \sum_{j=1}^k -p_j\partial_{q_j},&
 X_z &= \sum_{j=k+1}^n \left(q_j\partial_{p_j} - p_j\partial_{q_j}\right).
\end{align*}
This allows us to write $\phi=x-y+z$ and $X_\phi=X_x-X_y+X_z$.\pagebreak[3]\\
The level surface $\Sigma^- :=\{\phi{=}{-}1\}$ is transverse to $Y$, so that $(\Sigma^-,\lambda|_{\Sigma^-})$ is contact. The set $S:=\{x{=}z{=}0,\;y{=}{+}1\}$ is an isotropic sphere in $\Sigma^-$ and $(\mathbb{R}^{2n},\omega,Y,\Sigma^-\mspace{-4mu},S)$ is the isotropic setup where we glue a symplectic handle $\mathcal{H}$ to a contact manifold. In order to specify $\mathcal{H}$, we choose a different Weinstein function $\psi$ on $\mathbb{R}^{2n}$, satisfying the following assumptions:
\begin{itemize}
 \item[$(\psi 1)$] $X_\psi=C_x{\cdot} X_x-C_y{\cdot} X_y+C_z{\cdot} X_z,$ where $C_x,C_y,C_z\in C^\infty(\mathbb{R}^{2n}),\;\;C_x,C_y,C_z> 0$,
 \item[$(\psi 2)$] $\psi=\phi$ on $\{\phi\leq -1\}$ except for a small neighborhood of $S$,
 \item[$(\psi 3)$] The closure $\overline{\{\psi{<} {-}1\}\cap\{\phi{>} {-}1\}}$ is diffeomorphic to $\overline{D^k{\times} D^{2n-k}}$.
\end{itemize}
The handle is then defined as $\mathcal{H}:=\overline{\{\psi{<} {-}1\}\cap\{\phi{>} {-}1\}}\;\,$ (see Fig. \ref{fig3}).
\pagebreak[1]
\begin{remark}~
 \begin{itemize}
  \item If $\psi(0)\neq -1$, it follows from $(\psi 1)$ that the level set $\Sigma^+:=\{\psi{=}{-}1\}$ is also a contact hypersurface, as $Y{\cdot} \psi > 0$ away from 0. Due to $(\psi 2)$, we have $\Sigma^-=\Sigma^+$ away from a neighborhood of $S$. Condition $(\psi 3)$ on the other hand assures that $\Sigma^+$ is obtained from $\Sigma^-$ by surgery along $S$.
  \item Condition $(\psi1)$ is automatically satisfied if $\psi=\psi(x,y,z)$ is a function on $x,y,z$ with $\partial_x \psi,\partial_z \psi >0$ and $\partial_y \psi<0$ as its Hamiltonian vector field is then given by
  \[ X_\psi= \left(\frac{\partial \psi}{\partial x}{\cdot} X_x + \frac{\partial \psi}{\partial y}{\cdot} X_y + \frac{\partial \psi}{\partial z}{\cdot} X_z\right).\]
  \item $\mathcal{H}$ stays unchanged if we take $\phi'=\mathfrak{a}{\cdot} \phi + \mathfrak{b}$ and $\psi'=\mathfrak{a}{\cdot}\psi+\mathfrak{b}$ for $\mathfrak{a},\mathfrak{b}\in\mathbb{R}$, $\mathfrak{a}> 0$, provided we set $\bigg.\quad \mathcal{H}=\overline{\big.\{\psi'{<} {-}\mathfrak{a}{+}\mathfrak{b}\}\cap\{\phi'{>} {-}\mathfrak{a}{+}\mathfrak{b}\}}$.
  \item  Consider the Lyapunov function $L(q,p):=\sum_{j=1}^k q_jp_j$. Note that $(\psi1)$ implies $X_{\psi}{\cdot} L>0$ away from the critical points of $L$, which shows that all periodic orbits of $X_{\psi}$ are contained in the set $\{x{=}y{=}0\}$. The same holds true for $\psi'=\mathfrak{a}{\cdot}\psi+\mathfrak{b}$.
 \end{itemize}
\end{remark}
It is not difficult to find a Weinstein function $\psi:\mathbb{R}^{2n}\rightarrow\mathbb{R}$ satisfying $(\psi 1)$--$(\psi 3)$. Fix two constants $\veps,\delta>0$ and choose a smooth function $g:\mathbb{R}\rightarrow (-\infty,1]$ such that
\begin{align}
 g(t) &\phantom{:}= \begin{cases}\frac{1}{1+2\veps}\cdot t& \text{ for }\quad t\leq 1\\ 1&\text{ for }\quad t\geq 1+ 3\veps\end{cases}\quad\text{ and }\quad 0\leq g'(t) \leq {\textstyle \frac{1}{1+2\veps}}.\notag\\
 \text{Then set }\hspace{0.5cm}\psi_\delta&:=x-y+z-(1{+}\veps)+(1{+}\veps)\cdot g\big(y+{\textstyle \frac{1}{\delta}}(x{+}z)\big),\label{eqXX}\\
 \mathcal{H}_\delta&:=\overline{\big.\{\psi_\delta{<}{-}1\}\cap\{\phi{>}{-}1\}}.\notag
\end{align}
\begin{remark}~\label{remarksonpsi}
 \begin{itemize}
  \item Decreasing $\veps$ or $\delta$ makes the handle thinner, i.e. $\Sigma^-\cap \mathcal{H}_\delta$ becomes smaller. However, we will fix $\veps$ and decrease only $\delta$.
  \item For reference, let us fix $\delta_0$, $\psi_{\delta_0}$, the associated handle $\mathcal{H}_{\delta_0}$  and the hypersurface $\Sigma^+:=\big\{\psi_{\delta_0}{=}\,{-}1\big\}$. A different choice of $\delta$  and $\psi_\delta$ defines a different handle $\mathcal{H}_\delta$ (see Fig.\ \ref{fig2handles}). However, attaching these handles to a symplectic manifold $W$ yields matching completions $\widehat{W{\cup} \mathcal{H}_{\delta_0}}$ and $\widehat{W{\cup} \mathcal{H}_\delta}$.\\
  Indeed, if $\delta\leq \delta_0$ then $\psi_\delta\geq\psi_{\delta_0}$ and as $\psi_{\delta_0},\psi_\delta$ both increase along flow lines of $Y$, each flow line of $Y$ not inside $\{x{=}z{=}0\}$ first hits $\Sigma^\delta{:=}\,\{\psi_\delta{=}\,{-}1\}$ and then $\Sigma^+{:=}\,\{\psi_{\delta_0}{=}\,{-}1\}$. Hence $\Sigma^\delta$ can be identified with a hypersurface in the symplectization of $\Sigma^+$ given as a graph of a function $f:\Sigma^+\rightarrow(-\infty,0]$, where $f(p)$ is the unique time such that $\varphi^{-f(p)}(p)\in\Sigma^{\delta}$ for the flow $\varphi^t$ of $Y$.
  \item \begin{figure}[!htb]
\centering
 \resizebox{7cm}{5cm}{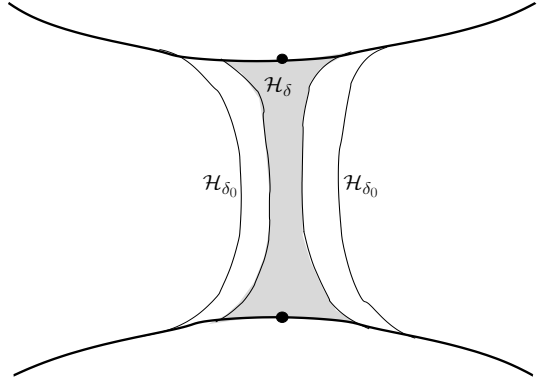}
 \\\caption{\label{fig2handles} Two handles}
\end{figure}Note that $\big(\{x{=}y{=}0\},Y|_{\{x{=}y{=}0\}}\big)$ is a Liouville subspace. This allows us to identify $\{x{=}y{=}0\}$ with the symplectization $\mathbb{R}\times\big(\Sigma^\delta{\cap}\{x{=}y{=}0\}\big)$. Moreover for $x{=}y{=}0$ and $z\leq \delta$ holds
  \[\psi_\delta(0,0,z) = \left(1+\frac{1+\veps}{\delta(1{+}2\veps)}\right) z -(1{+}\veps).\]
  As $\varphi^t$ satisfies $z(\varphi^t(q,p))=e^t{\cdot} z(q,p)$ for any $(q,p)\in\mathbb{R}^{2n}, \;t\in\mathbb{R}$, we can for $z\leq \delta$ express $\psi_\delta$ in symplectization coordinates $(r,p)$ on $\{x{=}y{=}0\}$ in the form
  \[\psi_\delta(r,p)=\mathfrak{a}_\delta{\cdot} e^r-(1{+}\veps).\]
  As $r=0$ corresponds to $\Sigma^\delta$ in $\{x{=}y{=}0\}$ and as the $z$-value on $\Sigma^\delta\cap\{x{=}y{=}0\}$ is smaller then $\delta$ and $\psi_\delta(\Sigma^\delta)=-1$, we have 
  \[-1=\mathfrak{a}_\delta{\cdot} e^0-(1{+}\veps)\qquad\Leftrightarrow\qquad \mathfrak{a}_\delta=\veps.\]
  This shows that $\psi_\delta$ has for any $\delta$ on $\Sigma^\delta\cap\{x{=}y{=}0\}$ the slope $\veps$ in radial direction.
 \end{itemize}
\end{remark}

\begin{remark}\label{dishandlesdifferentcontact}
 Let $\alpha_0=\lambda|_{\Sigma^-}$ be the contact form on $\Sigma^-$ and let $f:\Sigma^-\rightarrow\mathbb{R}$ be a smooth function. Then $\alpha_f=e^f{\cdot}\alpha_0$ is a contact form on $\Sigma^-$ defining the same contact structure. As discussed in \ref{secsetup}, $\alpha_f$ is the contact form on the hypersurface $\Sigma^-_f=\big\{(f(p),p)\,|\,p\in\Sigma^-\big\}$ in the symplectization of $(\Sigma^-,\alpha_0)$. Note that $(\Sigma_f^-,\alpha_f)$ is easily identified with a contact hypersurface in $\mathbb{R}^{2n}$ (also denoted by $\Sigma_f$) via the Liouville flows of $\partial_r$ on $\mathbb{R}{\times}\Sigma^-$ and of $Y$ on $\mathbb{R}^{2n}$.\\
 As being isotropic depends only on the contact structure, $S$ is isotropic in $\Sigma^-$ and $\Sigma^-_f$. More precisely, we identify $S$ with $S_f=\big\{(f(p),p)\big|\,p\in S{\subset} \Sigma^-\big\}$. We can define (as above) a handle $\mathcal{H}_f$ attached to $\Sigma^-_f$ along $S_f$ instead of $\Sigma^-$ along $S$. If $\mathcal{H}_f$ is sufficiently thin and $f|_S\leq 0$, we can understand $\mathcal{H}_f$ as a set inside $\mathcal{H}$ as follows: Denote by $\Sigma_f^+$ the other boundary of $\mathcal{H}_f$. Using again the flows of $\partial_r$ and $Y$, we can identify $\Sigma_f^+$ with the hypersurface $\big\{(f(p),p)\,\big|\,p\in\Sigma^+\big\}$ in the symplectization of $\Sigma^+$ and consequently with a hypersurface in $\mathbb{R}^{2n}$ which agrees with $\Sigma_f^-$ outside a compact set. The compact region bounded by $\Sigma_f^-$ and $\Sigma_f^+$ in $\mathbb{R}^{2n}$ is then identified with $\mathcal{H}_f$ (see Fig.\ \ref{fig2handles2contact}).
 \begin{figure}[!htb]
\centering
 \resizebox{7cm}{5cm}{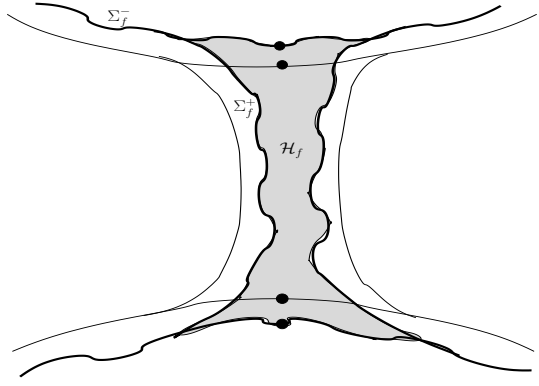}
 \\\caption{\label{fig2handles2contact} A handle for a different contact form}
\end{figure}
\end{remark}

\subsection{Linear extensions over the handle}\label{ExPsi}
In the proof of the Invariance Theorem, we need Hamiltonians $H\in Ad(W{\subset}W\cup\mathcal{H})$, i.e.\ $H|_W<0$, positive outside $W$ and linear on $[R,\infty)\times\partial(W{\cup}\mathcal{H})$ in $\widehat{W{\cup}\mathcal{H}}$ for $R$ large. As we saw in the previous section, it is easy to extend a Hamiltonian $H$ on $W$ to $W{\cup}\mathcal{H}$, if it is linear near $\partial W$. However, the extension of $H$ to $\widehat{W{\cup}\mathcal{H}}$ is more subtle, in particular if $X_H$ should have outside $W$ only few 1-periodic orbits with negative action, possibly only one. Unfortunately, Cieliebak's fundamental article on subcritical handle attachment, \cite{Cie}, has a gap exactly at this point. In particular, it does not address the following difficulties:
\begin{remark} \label{dis1}
Any Hamiltonian $H$ on $\widehat{W{\cup} \mathcal{H}}$ has at least one critical point outside $W$, without loss of generality at the center $c$ of $\mathcal{H}$. Note that $H>0$ outside of $W$, so that $c$ is a constant $X_H$-orbit with negative action.\pagebreak[3]\\
Fix a hypersurface $\Sigma^+$ with contact form $\alpha^+$. This fixes cylindrical coordinates $(z,p)\in(0,\infty){\times}\Sigma^+$ on the cylindrical end of $\widehat{W{\cup}\mathcal{H}}$. Here, we switch again to $z=e^r$ as radial coordinate, such that $\Sigma^+{=}\{z{=}1\}$. In order to keep control of the 1-periodic orbits of $X_H$, we assume that outside a compact neighborhood of $\mathcal{H}$ holds that $H$ is of the form
\[H(z,p)=\mathfrak{a}\cdot z -\mathfrak{b},\]
where $\mathfrak{b}>0$ is a constant such that $H$ restricted to $\partial W=\{z{=}z^-\}$ is negative. Moreover, we may assume that on the cocore of $\mathcal{H}$ (the set $\{x{=}y{=}0\})\quad H$ is cylindrical, i.e.\ $H(z,p)=h(z)$ for a function $h:\mathbb{R}^+\rightarrow\mathbb{R}^+$. As $H$ is differentiable, we find that $h(0):=\lim_{z\rightarrow 0}h(z)=H(c)$ extends $h$ continuously. Note that $H(c)$ and $\mathfrak{b}$ may be arbitrarily large, if $H$ is part of an increasing cofinal sequence of Hamiltonians.\\
As $H=\mathfrak{a}{\cdot} z\, {-}\mathfrak{b}$ holds everywhere away from $\mathcal{H}$, we find a $z^+$, such that $h(z)=\mathfrak{a}{\cdot} z\, {-}\mathfrak{b}\;\forall z\geq z^+$. Without loss of generality, we may assume $z^+$ to be minimal. Using the Mean Value Theorem on $h$, we find some $z\in(0,z^+)$ such that:
\begin{align*}
  h(z^+)-h(0) &= h'(z)\cdot(z^+-0) &&\Leftrightarrow& h'(z)\cdot z^++H(c)&=h(z^+)=\mathfrak{a}\cdot z^+-\mathfrak{b}\qquad\\
 &&&\Rightarrow& (h'(z)-\mathfrak{a})\cdot z^+ &=-(\mathfrak{b}+H(c)).\tag{$\ast$}
\end{align*}
Now, we consider several cases
\begin{itemize}
 \item If $H$ has a constant slope on $\{x{=}y{=}0\}$, i.e.\ if $h'(z){=}\mathfrak{a}\Leftrightarrow h(z){=}\mathfrak{a}z\,{-}\mathfrak{b}$ for all $z$, then $(\ast)$ yields a contradiction, as $\mathfrak{b}{+}H(c)\neq 0$.
 \item If no closed $X_H$-orbits exist in $\{x{=}y{=}0\}$, then $|h'(z)-\mathfrak{a}|< \mu$, where $\mu$ is a constant such that $dist.(\mathfrak{a}, spec(\alpha^+))< \mu$. Here, $(\ast)$ yields $z^+>\frac{b+H(c)}{\mu}$.
 \item If $|h'(z)-\mathfrak{a}|\geq \mu$ for some $z$, then there exists $\tilde{z}$ with $h'(\tilde{z})\in spec(\alpha^+)$, i.e.\ there exists a 1-periodic $X_H$-orbit $\gamma$ on $\{x{=}y{=}0\}$. Let $z^0$ be the smallest of these $z$-values. Then $h'(z)>h'(z^0)-\mu$ for all $z\in(0,z^0)$. This implies that
 \begin{align*}
  && h(z^0)&\geq h(0)+z^0\cdot(h'(z^0)-\mu)\\
  &\Rightarrow& \mathcal{A}^H(\gamma)=h'(z^0){\cdot}z^0-h(z^0)&\leq z^0{\cdot}\mu-h(0)=z^0{\cdot}\mu-H(c).
 \end{align*}
 If we require that the action of $\gamma$ is positive, then $\mathcal{A}^H(\gamma)>0\Rightarrow z^+>z^0>\frac{H(c)}{\mu}$.
\end{itemize}
The first two cases include the case that apart from $c$, no new 1-periodic $X_H$-orbits are created on the cocore, while the third case covers the situation that all new orbits have positive action. In both cases, we find that $z^+$ is arbitrarily large, as $\mathfrak{b}$ and $H(c)$ are arbitrarily large. As $z^+$ was chosen to be minimal, this implies for a cofinal sequence of Hamiltonians $(H_n)$ that the set $\{x{=}y{=}0\}\cap\{z{\leq}z^+\}$ where $H_n$ is not linear increases in $n$ and covers in the limit the cocore $\{x{=}y{=}0\}$. Up to now, the author does not see how this can be accomplished without creating new 1-periodic orbits near $\mathcal{H}$.
\end{remark}

The solution to the problem described above is to allow more closed $X_H$-orbits with negative action on the cocore $\{x{=}y{=}0\}$ of $\mathcal{H}$. Using the Lyapunov function $f$, we can show that we create new 1-periodic $X_H$-orbits only on $\{x{=}y{=}0\}$. These can be explicitly described and one can show that they do not contribute to the Symplectic Homology. For the construction of such an $H$, we need the following two technical lemma:
\begin{lemma}\label{lemHamilt}
 Consider $\mathbb{R}^{2n}$ with the standard symplectic structure, the Liouville vector field $Y$  and the functions $x,y,z$ with Hamiltonian vector fields $X_x,X_y,X_z$ as given in \ref{handle}. Let $\Sigma\subset\mathbb{R}^{2n}$ be a smooth hypersurface transverse to $Y$ (i.e.\ $\Sigma$ contact) such that its Reeb vector field $R$ is of the form 
 \[R=c_x{\cdot} X_x-c_y{\cdot} X_y + c_z{\cdot} X_z,\qquad\qquad c_x,c_y,c_z \in C^\infty(\Sigma),\qquad c_x,c_y,c_z>0.\]
 Consider the function $\tilde{h}_\Sigma(y,r)=\mathfrak{a}{\cdot} e^r+\mathfrak{b}$ on $\mathbb{R}\times\Sigma$ and let $h_\Sigma:=\tilde{h}_\Sigma\circ \Phi^{-1}$ be its pushforward onto $\mathbb{R}^{2n}$ by the symplectic embedding $\Phi:\mathbb{R}\times\Sigma\rightarrow\mathbb{R}^{2n}, (r,p)\mapsto \varphi^r(p)$ via the flow $\varphi^t$ of $Y$. Then, the Hamiltonian vector field $X_h$ of $h_\Sigma$ is of the form
 \begin{align*}
  &X_h=C_x{\cdot} X_x - C_y{\cdot} X_y + C_z{\cdot} X_z,\qquad C_x,C_y,C_z\in C^{\infty}(\mathbb{R}^{2n}),\qquad C_x,C_y,C_z>0.
 \end{align*}
\end{lemma}\pagebreak[1]
\begin{remark}
The assumptions on $\Sigma$ are satisfied, if $\Sigma=\psi^{-1}(c)$ for a function $\psi$ on $x,y,z$ with $\partial_x \psi\big|_\Sigma,\partial_z \psi\big|_\Sigma>0$ and $\partial_y \psi\big|_\Sigma<0$ and $0\not\in\Sigma$.
\end{remark}
\begin{proof}
 As $X_{\tilde{h}}=\mathfrak{a}{\cdot} R$ on $\mathbb{R}{\times}\Sigma$, it follows that on $\mathbb{R}^{2n}$ holds $X_h|_{\varphi^t(\Sigma)}=\mathfrak{a} e^t{\cdot} R_t$, where $R_t$ is the Reeb vector field on $\varphi^t(\Sigma)$. By assumption, the Reeb vector field $R$ on $\Sigma$ satisfies
 \begin{align*}
  R &= c_x X_x-c_y X_y+c_z X_z = \textstyle\sum_{j=1}^k\left( c_x\,2q_j\partial_{p_j}+c_y\,p_j\partial_{q_j}\right) + c_z\sum_{j=k+1}^n \left(q_j\partial_{p_j}-p_j\partial_{q_j}\right).
 \end{align*}
 Since $Y=\sum_{j=1}^k\big(2q_j\partial_{q_j}-p_j\partial_{p_j}\big)+\frac{1}{2}\sum_{j=k+1}^n\big(q_j\partial_{q_j}+p_j\partial_{p_j}\big)$, its flow $\varphi^t$ is
 \[\varphi^t(q,p)=\Big(\underbrace{...\,,e^{2t}\cdot q_j,e^{-t}\cdot p_j,...}_{j=1,...\,,k}\text{\huge,}\underbrace{...\,,e^{t/2}\cdot q_j,e^{t/2}\cdot p_j,...}_{j=k+1,...\,,n}\Big).\]
 As $\mathcal{L}_Y\lambda=\lambda$ and $\mathcal{L}_Y\omega=\omega$, we find for $R$ and any $\xi\in T_{\varphi^t(p)}\varphi^t(\Sigma)$ that
 \begin{align*}
  \lambda_{\varphi^t(p)}\big(D\varphi^t_p R\big) &=\big({\varphi^t}^\ast \lambda\big)_p(R)&&=e^t\cdot \lambda_p(R)&&=e^t,\\
  \omega_{\varphi^t(p)}\big(D\varphi^t_p R\,,\,\xi\big) &=\big({\varphi^t}^\ast\omega\big)_p\big(R,(D\varphi^t_p)^{-1}(\xi)\big)&&=e^t\cdot\omega_p\big(R,(D\varphi^t_p)^{-1}(\xi)\big)&&=0,
 \end{align*}
as $R$ is the Reeb vector field and $(D\varphi^t_p)^{-1}(\xi)\in T\Sigma$. This shows that $e^{-t}\cdot D\varphi^t R$ is the Reeb vector field $R_t$ of $\varphi^t(\Sigma)$. Hence, $X_h$ is of the announced form as
\[X_h|_{\varphi^t(\Sigma)}=\mathfrak{a} e^t{\cdot} R_t = \mathfrak{a}\cdot D\varphi^t(R) = \mathfrak{a} e^{-t}c_xX_x - \mathfrak{a} e^{2t}c_yX_y+\mathfrak{a} e^{t/2}c_zX_z.\]
\end{proof}
\begin{lemma} \label{interpolation}
Let $(\Sigma,\alpha)$ be a compact contact manifold with contact form $\alpha$ and symplectization $\big(\mathbb{R}{\times}\Sigma,\omega{=}d(e^r\alpha)\big)$ and let $||\cdot||$ denote a norm with respect to a metric given by an $\omega$-compatible almost complex structure. Let $\veps,\delta,c>0$ be constants.\\
Then there exists a smooth monotone increasing function $g:\mathbb{R}\rightarrow[0,1]$ such that 
 \[g(e^r)=0\quad\text{ for }\quad r\leq -\veps\qquad\text{ and }\qquad g(e^r)=1\quad\text{ for }\quad r\geq 0\tag{$\ast$}\]
 and for all $\phi,\psi\in C^1(\mathbb{R}\times\Sigma)$ with $\phi|_{\{0\}\times\Sigma}=\psi|_{\{0\}\times\Sigma}$ and $|\partial_r\phi(r,p)-\partial_r\psi(r,p)|<c$ for all $(r,p)\in[-\veps,0]{\times}\Sigma$, holds that the Hamiltonian vector fields $X_\phi, X_\psi$ satisfy
 \[\sup_{r\in\Sigma \atop p\in\mathbb{R}}\left|\left|X_{\phi+(\psi-\phi)g}(r,p)-\Big(X_\phi(r,p)+\big(X_\psi(r,p)-X_\phi(r,p)\big)\cdot g(e^r)\Big)\right|\right|\leq \delta.\tag{$\ast\ast$}\]
 In other words, we can interpolate between $\phi$ and $\psi$ along $[-\veps,0]\times\Sigma$, such that the Hamiltonian vector field $X_{\phi+(\psi-\phi)g}$ of the interpolation is arbitrary close to the interpolation of the Hamiltonian vector fields $X_\phi$ and $X_\psi$.
\end{lemma}
\begin{proof}
 As the Hamiltonian vector field of $e^r$ is the Reeb vector field $R$, we calculate
 \[X_{\phi+(\psi-\phi)g}(r,p)=X_\phi(r,p)+\big(X_\psi-X_\phi\big)(r,p)\cdot g(e^r)+\big(\psi-\phi\big)(r,p)\cdot g'(e^r)\cdot R(p).\]
Therefore, $(\ast\ast)$ translates to
\[ \left|\left|\big(\psi-\phi\big)(r,p)\cdot g'(e^r)\cdot R(p)\right| \right| \leq \delta \qquad\forall (r,p)\in[-\veps,0]{\times}\Sigma.\]
Using $\phi|_{\Sigma\times\{0\}}=\psi|_{\{0\}\times\Sigma}$, we can estimate the left hand side by
\begin{align*}
 \left|\left|\big(\psi{-}\phi\big)(r,p) g'(e^r) R(p)\right| \right|&=\left|\left|-\mspace{-7mu}\int_r^0\mspace{-7mu}\partial_s\big(\psi{-}\phi\big)(s,p)\,ds\cdot g'(e^r) R(p)\right| \right|\leq -rc g'(e^r) ||R||_\infty.
\end{align*}
 Writing $z=e^r$, we find that $(\ast\ast)$ is satisfied, if $0\leq g'(z) \leq \frac{-\delta}{c||R||\log z}$ for all $z\in[e^{-\veps},1]$. As $\int_{e^{-\veps}}^1 \frac{-\delta}{c||R||\log z}\,dz=\infty$, we can choose a smooth function $\tilde{g}$ satisfying
 \[0\leq \tilde{g}(z)\leq \frac{-\delta}{c||R||_\infty \log z},\quad \tilde{g}\equiv 0 \;\text{ for } \;z\leq e^{-\veps}\;\text{ or }\;z\geq 1 \quad \text{ and }\; \int_{e^{-\veps}}^1 \tilde{g}(z)\,dz = 1.\]
Setting $\displaystyle g(e^r)=g(z):=\int_{e^{-\veps}}^z\tilde{g}(s)\,ds$ then gives the desired function. 
\end{proof}\pagebreak[3]
Now, we construct $H$ in two steps, first extending an admissible Hamiltonian $H|_W$ on $W$ via $\psi_\delta$ to $W{\cup} \mathcal{H}_\delta$, then constructing a linear extension to $\widehat{W{\cup} \mathcal{H}_\delta}$. For simplicity, we assume that $H|_W=1{\cdot} e^r{-}2$ near $\Sigma^-=\{r\}{\times}\partial W$. Then $\Sigma^-=\{H|_W{=}\,{-}1\}$ and $H|_W$ has slope 1 on $\Sigma^-$. For the general case $H|_W=\mathfrak{a}{\cdot} e^r{+}\mathfrak{b}$, we take $H$ as constructed below and extend $H|_W$ with $\mathfrak{a}{\cdot} H {+}\mathfrak{b}{+}2\mathfrak{a}$. Explicitly, the two steps are as follows:
\begin{itemize}
 \item \underline{Step 1}: Recall that the isotropic sphere $S\subset\Sigma^-=\{\phi{=}\,{-}1\}$ is given by
\[ S:=\{x{=}z{=}0,\; y{=}1\}.\]
For a small neighborhood $U$ of $S$ identify the isotropic setup $\big(U,\omega,Y,\Sigma^-\cap U,S\big)$ with an isotropic setup $\big(U,\omega,Y,(\{r\}{\times}\partial W)\cap U,S\big)$ in $W$ (by abuse of notation, we use the same letters for identified objections on $W$ resp. $\mathbb{R}^{2n}$). Consider the following linear Hamiltonian \[\tilde{h}^-_\Sigma:\mathbb{R}\times\Sigma^-\rightarrow\mathbb{R},\quad \tilde{h}_\Sigma^-(,p)=1{\cdot} e^r-2\]
and its pushforward $h_\Sigma^-$ onto $\mathbb{R}^{2n}$ defined by $h^-_\Sigma=\tilde{h}_\Sigma^-\circ \Phi^{-1}$, where $\Phi(r,p)=\varphi^r(p)$ is the symplectic embedding provided by the flow $\varphi^t$ of $Y$. Note that $H|_W$ coincides with $h_\Sigma^-$ under the identification of isotropic setups.\\
As the Reeb vector field $R_{\Sigma^-}$ of $(\Sigma^-,\lambda|_{T\Sigma^-})$ coincides with the Hamiltonian vector field $X_\phi$ on $S$, we find $X_{h_\Sigma^-}=R_{\Sigma^-}=X_\phi$ and hence $dh_\Sigma^-=d\phi$ on $S$. As also $h_\Sigma^-(\Sigma^-)=\phi(\Sigma^-)=-1$, we find that $h_\Sigma^-$ and $\phi$ coincide up to first order on $S$. Therefore, given any neighborhood $U^\delta$ of $S$, there exists a function $\hat{\phi}\in C^\infty(\mathbb{R}^{2n})$ and a neighborhood $\hat{U}^\delta\subset U^\delta$, such that $\hat{\phi}\equiv h_\Sigma^-$ on $\mathbb{R}^{2n}\setminus U^\delta,\; \hat{\phi}\equiv\phi$ on $\hat{U}^\delta$ and $\hat{\phi}$ is arbitrarily $C^1$-close to $h_{\Sigma}^-$. As $X_\phi=X_x-X_y+X_z$ and $X_{h_\Sigma^-}=C_x^- X_x-C_y^- X_y+C_z^- X_z$ with $C_x^-,C_y^-,C_z^->0$ by Lem.\ \ref{lemHamilt}, we can additionally arrange that
\[X_{\hat{\phi}} = \hat{C}_x{\cdot} X_x - \hat{C}_y{\cdot} X_y + \hat{C}_z{\cdot} X_z\quad\text{ with }\quad \hat{C}_x,\hat{C}_y,\hat{C}_z>0.\]
As the $X_x$- and $X_z$-part of $X_\phi$ and $X_{h^-_\Sigma}$ are both 0 on $\{x{=}z{=}0\}$, we can make $U^\delta$ and $\hat{U}^\delta$ arbitrarily thin in the $x$- and $z$-direction, while keeping a fixed size in the $y$-direction. This allows us to choose in (\ref{eqXX}) for the definition of $\psi_\delta$ one $\veps$ for all handles. Fix $\veps$ sufficiently small and choose $\delta$ depending on $U^\delta$ so small such that the lower boundary $\mathcal{H}_\delta{\cap}\Sigma^-=\Sigma^-\setminus\Sigma^\delta$ lies in $\hat{U}^\delta$. Then set
\begin{align*}
 \hat{H}& : W\cup \mathcal{H}_\delta\rightarrow \mathbb{R},\qquad\qquad \hat{H}=
 \begin{cases}\psi_\delta & \text{ on }\big(\hat{U}^\delta\cap\{\phi{\leq}{-}1\}\big)\cup \mathcal{H}_\delta\\ 
 \hat{\phi} & \text{ on }\big(U^\delta\cap\{\phi{\leq}{-}1\}\big)\setminus \hat{U}^\delta\\
 H|_W & \text{ on }W\setminus U^\delta
\end{cases}.
\end{align*}
Since $\psi_\delta=\hat{\phi}$ outside a small neighborhood of $\mathcal{H}_\delta$, $\phi=\hat{\phi}$ on $\hat{U}^\delta$ and $\hat{\phi}=h_\Sigma^-=K$ outside $U^\delta$, we find that $\hat{H}$ is smooth on its domain. Moreover, as $\psi_\delta,\hat{\phi}$ and $h_\Sigma^-$ satisfy $(\psi1)$ on $U\cup \mathcal{H}_\delta$, so does $\hat{H}$, i.e.\ there exist smooth $\hat{C}_x,\hat{C}_y,\hat{C}_z>0$ such that
\[X_{\hat{H}}=\hat{C}_x{\cdot} X_x-\hat{C}_y{\cdot} X_y+\hat{C}_z{\cdot} X_z.\]
See Figure \ref{fig6} for the areas where $\hat{H}$ is defined.
\begin{figure}[ht]
\centering
\begin{minipage}[htb]{6cm}
 \resizebox{6cm}{!}{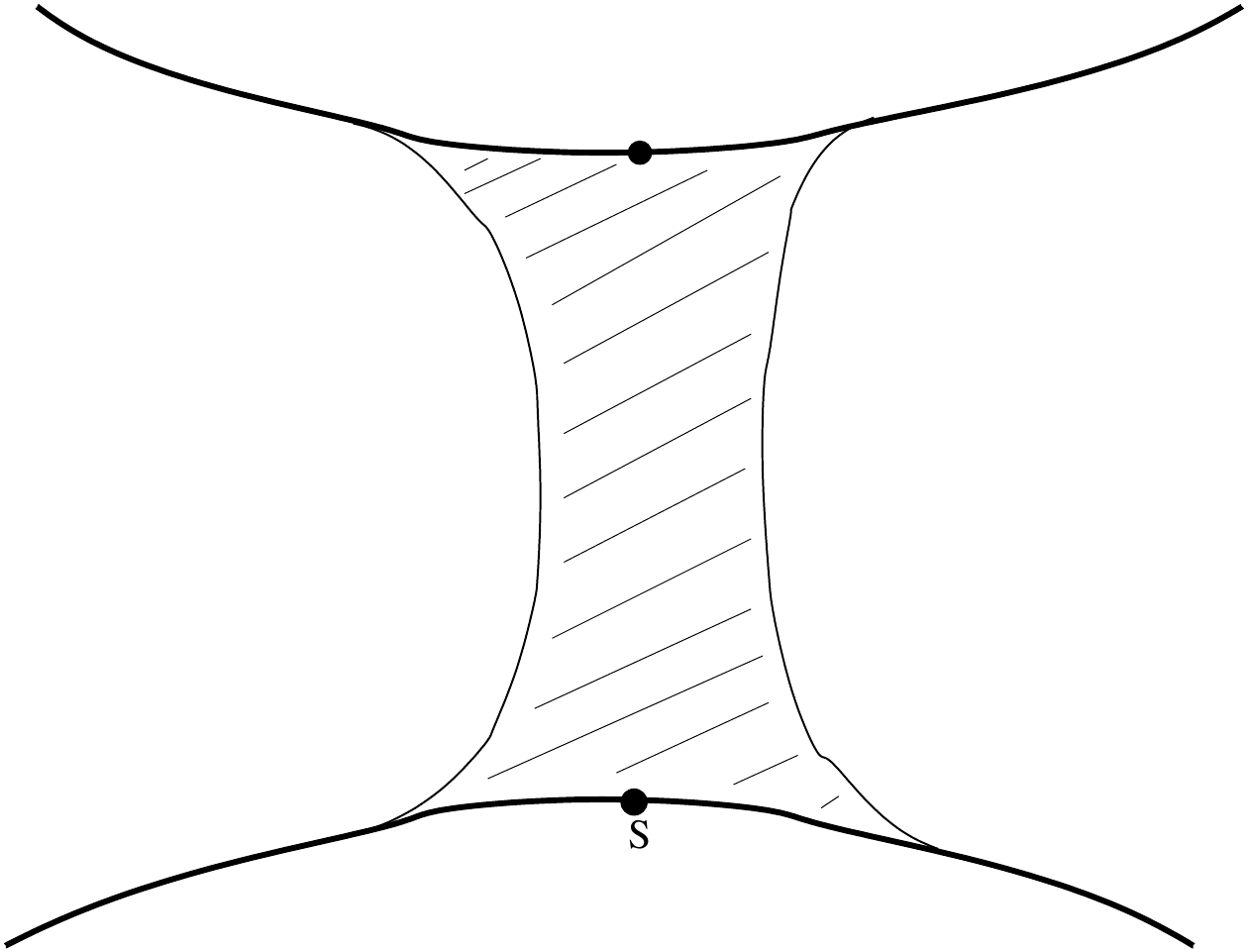}
 \centering{The handle $\mathcal{H}_\delta$}
 \end{minipage}
 \begin{minipage}[htb]{6cm}
  \resizebox{6cm}{!}{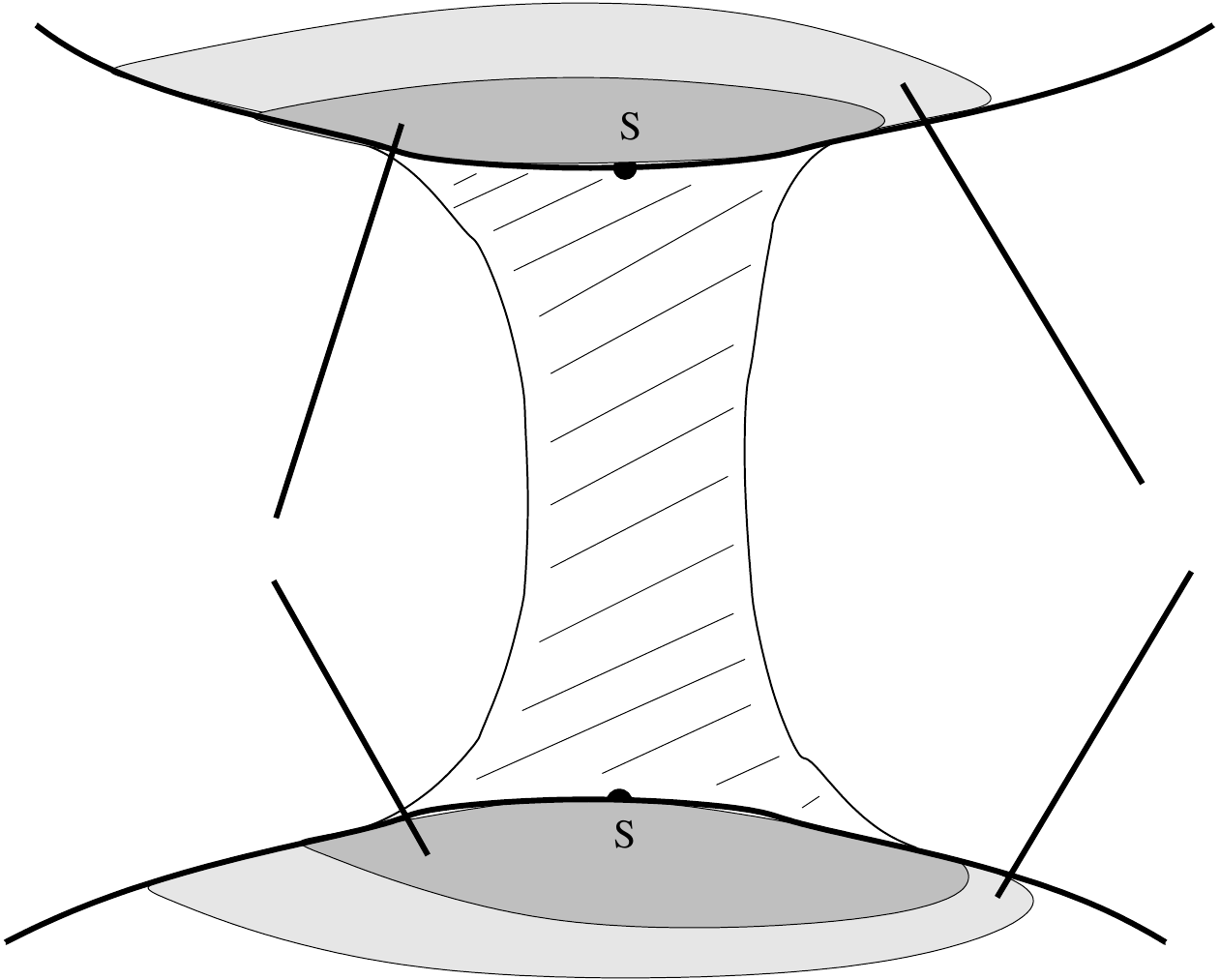}
 \centering{$U^\delta\cap\{\phi{\leq}\,{-}1\}$ and $\hat{U}^\delta\cap\{\phi{\leq}\,{-}1\}$}
 \end{minipage}
\begin{minipage}[htb]{6cm}
 \resizebox{6cm}{!}{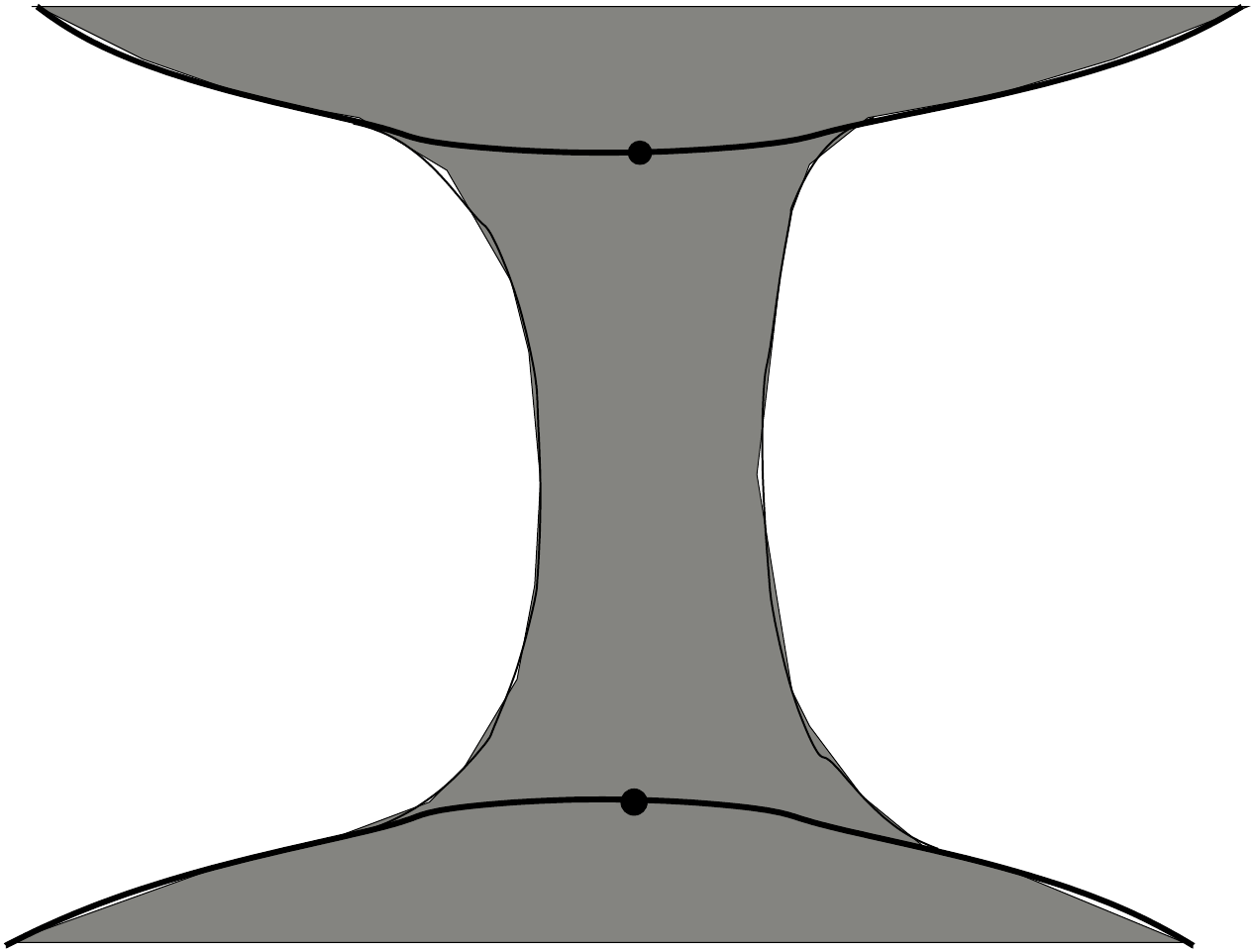}
 \centering{The area, where $\hat{H}$ is defined}
\end{minipage}
\begin{minipage}[htb]{6cm}
 \resizebox{6cm}{!}{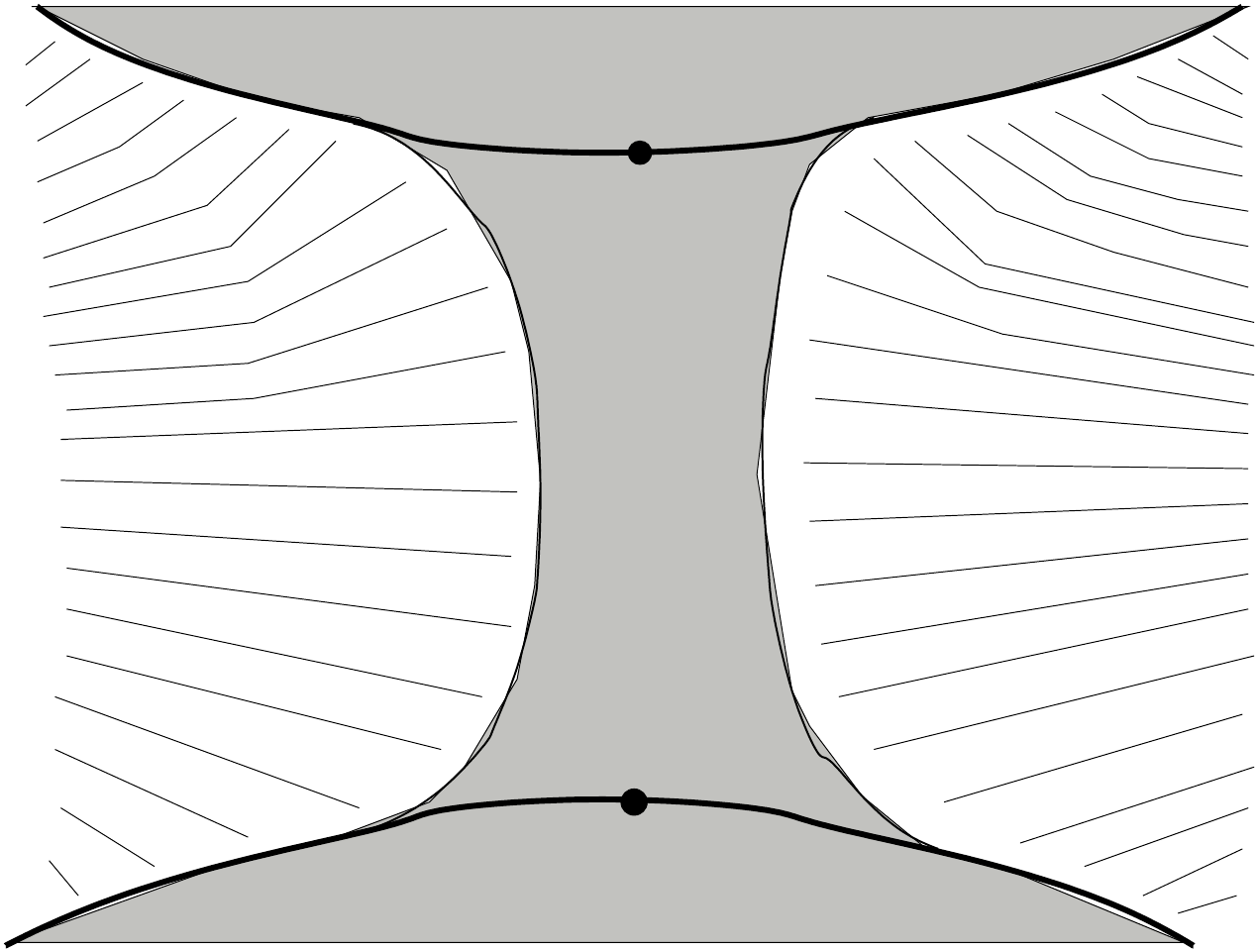}
 \centering{The completion of $[0,\infty)\times\Sigma^\delta$}
\end{minipage}
 \caption{\label{fig6}Areas, where $\hat{\psi}$ is defined}
\end{figure}
\item \underline{Step 2} Consider on $\mathbb{R}\times\Sigma^\delta$ the linear function $h^\delta$ given by $\;\tilde{h}^\delta(r,p)=1{\cdot} e^r-2$ \hfill and its pushforward $h^\delta=\tilde{h}^\delta$ to $\mathbb{R}^{2n}$ by the flow of $Y$.\\
On $\widehat{W{\cup} \mathcal{H}_\delta}\setminus\big(W{\cup} \mathcal{H}_\delta\big)$ we define $H$ by $H=h^\delta$.\\
On $W{\cup} \mathcal{H}_\delta$ we define $H$ as an interpolation between $\widehat{H}$ and $h^\delta$, i.e.\
\begin{align*}H:=\hat{H}+\big(h^\delta-\hat{H}\big)\cdot g(h^\delta),
\end{align*}
where $g$ is a function given by Lem.\ \ref{interpolation} such that for $\tau>0$ small holds
\[g\big(e^r{-}2\big)= 0\quad \text{ for } \quad r\leq -\tau\qquad\text{ and }\qquad g\big(e^r{-}2\big) = 1\quad \text{ for } \quad r\geq 0\]
and $X_H$ is arbitrarily close to the interpolation of $X_{\hat{H}}$ and $X_{h^\delta}$.\\
On $\mathbb{R}{\times}\big(\partial W{\setminus} U^\delta\big)$ holds $H|_W=h_\Sigma^-=e^r{-}2=h^\delta$. Hence we interpolate in this area between the same Hamiltonians. However on $\mathbb{R}{\times}\big(\partial(W{\cup} \mathcal{H}_\delta)\setminus(\partial W{\setminus} U^\delta)\big)$, we have by Lem.\ \ref{lemHamilt} and the construction of $\widehat{H}$ that
\begin{align*}
 X_{\hat{H}} &= \hat{C}_x{\cdot} X_x-\hat{C}_y{\cdot} X_y + \hat{C}_z{\cdot} X_z, & \hat{C}_x,\hat{C}_y,\hat{C}_z&>0\\
 X_{h^\delta} &= C_x^\delta{\cdot} X_x-C_y^\delta{\cdot} X_y + C_z^\delta{\cdot} X_z, & C_x^\delta,C_y^\delta,C_z^\delta&>0\\
 X_{g(h^\delta)}&=g'(h^\delta)\big(C_x^\delta{\cdot} X_x-C_y^\delta{\cdot} X_y + C_z^\delta{\cdot} X_z\big).
\end{align*}
Hence $X_H$ is also of this form and satisfies $(\psi1)$, as it is close to the interpolation of $X_{\hat{H}}$ and $X_{h^\delta}$. It follows from $(\psi1)$ with the help of the Lyapunov function $f$ that 1-periodic orbits of $X_H$ are either 1-periodic orbits of $X_H|_W$ inside $W$ or the constant orbit at the center of $\mathcal{H}_\delta$, non-constant orbits on the handle near $\Sigma^\delta\cap\mathcal{H}_\delta$ or Hamiltonian orbits that pass over the handle, go into $W$ and come back. The last type of orbits is not 1-periodic if the handle is chosen sufficiently thin, so that the only new orbits after attaching the handle are the constant orbit at the center of $\mathcal{H}_\delta$ and those on the outer boundary of $\mathcal{H}_\delta$.
\end{itemize}

\subsection{New closed orbits and their Conley-Zehnder indices}\label{sec5.4}
We saw at the end of the last paragraph that all 1-periodic orbits of $X_H$ that are not orbits of $X_H|_W$ lie on the handle. Moreover, we saw that there exist functions $C_x, C_y, C_z>0$ such that near the handle $X_H=C_x{\cdot} X_x-C_y{\cdot} X_y+C_z{\cdot} X_z$. Hence all new 1-periodic orbits are contained in $\{x{=}y{=}0\}$. There we have
\[H(0,0,z)=\psi_\delta(0,0,z)=\Big(1+\frac{1+\veps}{\delta(1{+}2\veps)}\Big)\cdot z-(1{+}\veps)=\veps{\cdot} e^r-(1{+}\veps)\]
in symplectization coordinates adapted to $\Sigma^\delta$ (see Rem.\ \ref{remarksonpsi}). As $h_\delta=e^r{-}2$, we find that on $\{x{=}y{=}0\}$ the function $C_z$ is an increasing function in $z=\frac{1}{2}\sum_{j=k+1}^n (q_j^2+p_j^2)$ interpolating between the constants $\big(1+\frac{1+\veps}{\delta(1+2\veps)}\big)$ and $\frac{1}{\veps}\big(1+\frac{1+\veps}{\delta(1+2\veps)}\big)$.\pagebreak[1]\\
For the proof of the Invariance Theorem, we need to determine the 1-periodic Hamiltonian orbits of $\mathfrak{a}{\cdot} H {+}\mathfrak{b}{+}2\mathfrak{a}$ for $\mathfrak{a},\mathfrak{b}\in\mathbb{R}$. The resulting Hamiltonian vector field is simply $\mathfrak{a}{\cdot} X_H$. Let $\varphi_H^t$ denote its flow. Then we calculate
\[\frac{d}{dt}z(\varphi^t_H)= dz(\mathfrak{a}{\cdot} X_H)=\mathfrak{a} C_z\cdot dz(X_z)=0.\]
It follows that $z$ is constant along flow lines of $\varphi_H^t$ and as $C_z$ is a function of $z$ on $\{x{=}y{=}0\}$ it follows that $C_z$ is also constant along flow lines of $\varphi_H^t$ on $\{x{=}y{=}0\}$. Introducing the complex coordinates $z_j=q_j+i\cdot p_j$, we find that 
\[ X_z=\Big(\underbrace{0,\dots,0}_{j=1,...,k}\, \textbf{,}\underbrace{\dots,iz_j,\dots}_{j=k+1,...,n}\Big).\]
As $X_x=X_y=0$ on $\{x{=}y{=}0\}$, the flow $\varphi^t_H$ on this set is given by
\begin{equation}\label{eqA}
  \varphi^t_H(0,\dots,0,z_{k+1},\dots,z_n)=\Big(0,\dots,0,e^{i\mathfrak{a} C_z t}\cdot z_{k+1},\dots, e^{i\mathfrak{a} C_z t}\cdot z_n\Big).
 \end{equation}
A non-constant orbit of $\mathfrak{a}{\cdot} X_H$ is hence 1-periodic if and only if $\mathfrak{a} C_z(z)\in 2\pi\mathbb{Z}$. As $z=\frac{1}{2}\sum_{j=k+1}^n(q_j^2{+}p_j^2)$, the 1-periodic orbits form families that are diffeomorphic to the standard sphere $S^{2(n-k)-1}$, except for the constant orbit at $z=0$.\bigskip\pagebreak[3]\\
Next, we calculate $\mu_{CZ}$ for these orbits. Let $\gamma$ denote a 1-periodic orbit of $\mathfrak{a}{\cdot} X_H$. By identifying $T_{\gamma(t)}\mathbb{R}^{2n}$ with $\mathbb{R}^{2n}$ in the obvious way, we obtain a path of matrices $\Phi_\gamma(t):=D\varphi^t_H(0)\in Sp(2n)$. Differentiating yields:
\begin{align*}
 \frac{d}{dt}\Phi_\gamma(t)&=\frac{d}{dt}D\varphi^t_H(0)=D\left(\frac{d}{dt}\varphi^t_H(0)\right)=D\; \mathfrak{a}{\cdot} X_H\big(\varphi^t_H(0)\big)\\
 &=\mathfrak{a}\cdot D\big(C_x X_x- C_y X_y+ C_z X_z\big)\big(\varphi^t_H(0)\big)\tag{cf.\ (\ref{eqXX})}\\
&=\mathfrak{a}\cdot diag\bigg(\underbrace{\dots,\Big(\begin{smallmatrix} 0 & C_y\\2C_x& 0\end{smallmatrix}\Big),\dots}_{j=1,...\,,k}\text{\Huge,}\underbrace{\dots, i C_z,\dots}_{j=k+1,...,n}\bigg)\circ\Phi_\gamma(t).
\end{align*}
Note that no derivatives of $C_x$ or $C_y$ are involved, as $X_x{=}X_y{=}0$ on $\{x{=}y{=}0\}$. Thus, $\Phi_\gamma$ is of block form $\Phi_\gamma = diag\big(\Phi^1_\gamma,...\,,\Phi_\gamma^n\big)$, where $\Phi^j_\gamma$ are paths of $2{\times}2$-matrices which are solutions of initial value problems with $\Phi^j_\gamma(0)=\mathbbm{1}$ and
\begin{align*}
 {\textstyle\frac{d}{dt}}\Phi^j_\gamma(t)&=\mathfrak{a}\begin{pmatrix}0 & C_y\\ 2C_x & 0\end{pmatrix}\Phi^j_\gamma(t) = \begin{pmatrix} 0 & -1\\ 1& 0\end{pmatrix}\mathfrak{a}\begin{pmatrix}  2C_x & 0\\0& -C_y\end{pmatrix}\Phi^j_\gamma(t)& j&=1,...\,,k\\
 {\textstyle\frac{d}{dt}}\Phi^j_\gamma(t)&=i\mathfrak{a} C_z\cdot\Phi^j_\gamma(t) & j&=k{+}1,...\,, n.
\end{align*}
As $\mathfrak{a}\left(\begin{smallmatrix}2C_x & 0\\ 0 & - C_y\end{smallmatrix}\right)$ has for all $t$ signature zero, it follows with (CZ0) that $\mu_{CZ}(\Phi^j_\gamma)=0$ for $1{\leq} j{\leq} k$. For $k{+}1{\leq} j{\leq} n$, we find $\Phi^j_\gamma(t)=e^{i\mathfrak{a} C_z\cdot t}$ and by (CZ1) that
\[\mu_{CZ}\big(\Phi^j_\gamma\big)= \left\lfloor\frac{\mathfrak{a} C_z}{2\pi}\right\rfloor+\left\lceil\frac{\mathfrak{a} C_z}{2\pi}\right\rceil,\qquad j=k{+}1,...\,n.\]
Using the direct sum property, the Conley-Zehnder index of $\gamma$ is given by
\[\mu_{CZ}(\gamma)=\mu_{CZ}(\Phi_\gamma) =\sum_{j=1}^n \mu_{CZ}\big(\Phi^j_\gamma\big)
=(n{-}k)\cdot \left(\left\lfloor\frac{\mathfrak{a} C_z}{2\pi}\right\rfloor+\left\lceil\frac{\mathfrak{a} C_z}{2\pi}\right\rceil\right).\]
As $C_z{\geq} 1$, we have $\mu_{CZ}(\gamma){\rightarrow} \infty$ for $\mathfrak{a}{\rightarrow}\infty$.\pagebreak[1]\\ Choosing a Morse function with 2 critical values on $S^{2(n-k)-1}$, we find that the Morse-Bott index for non-constant $\gamma$ is given by
\begin{align*}
 \mu(\gamma)&=(n{-}k)\left(\left\lfloor\frac{\mathfrak{a} C_z}{2\pi}\right\rfloor+\left\lceil\frac{\mathfrak{a} C_z}{2\pi}\right\rceil\right)-\frac{1}{2}\big(2(n{-}k)\,{-}1\big)+\frac{1}{2}+\left\lbrace\begin{smallmatrix}0\\2(n-k)-1\end{smallmatrix}\right.\\
 &=(n{-}k)\left(\left\lfloor\frac{\mathfrak{a} C_z}{2\pi}\right\rfloor+\left\lceil\frac{\mathfrak{a} C_z}{2\pi}\right\rceil\right)-(n{-}k)+1+\left\lbrace\begin{smallmatrix}0\\2(n-k)-1\end{smallmatrix}\right..
\end{align*}
Since for 1-periodic $\gamma$ holds $\mathfrak{a} C_z(z) \in 2\pi\mathbb{Z}$, their Morse-Bott indices are of the form
\begin{equation}\label{indexonhandle}
 \mu(\gamma)=(n{-}k)(2l{-}1)+\left\lbrace\begin{smallmatrix}1\phantom{(n-k)}\\2(n-k)\end{smallmatrix}\right.\qquad\text{ for } l\in\mathbb{N}.
\end{equation}

\subsection{Handle attachment and Symplectic Homology}\label{secHandleAt}
In this section, we use our construction of Hamiltonians on subcritical $k$-handles $\mathcal{H}$ to prove the Invariance Theorem, Cor.\ \ref{corinvsur} for $SH^+$ and Prop. \ref{propcofinalHamilt}, which states that the a.f.g.\ property is also invariant under handle attachment.
\begin{proof}[\textnormal{\textbf{Proof of Theorem \ref{theoinvsur} (Invariance Theorem)}}]~\\
We have to show that $SH_\ast(W)\cong SH_\ast(V)$ and $SH^\ast(W)\cong SH^\ast(V)$, if $V$ is obtained from $W$ by attaching a subcritical handle $\mathcal{H}$. In fact, we show that the transfer maps
\[\pi_\ast(W,V): SH_\ast(V)\rightarrow SH_\ast(W)\quad\text{ and }\quad \pi^\ast(W,V): SH^\ast(W)\rightarrow SH^\ast(V),\]
for the embedding $W\subset W{\cup}\mathcal{H}= V$ are isomorphisms. By Sec.\ \ref{sectransfer}, these maps are the truncation maps $SH_\ast(V)\rightarrow SH_\ast^{\geq0}(W{\subset}V)$ and $SH^\ast_{\geq0}(W{\subset}V)\rightarrow SH^\ast(V)$ composed with the isomorphisms $SH_\ast^{\geq0}(W{\subset}V)\cong SH_\ast(W)$ and $SH^\ast_{\geq0}(W{\subset}V)\cong SH^\ast(W)$. Therefore, we only have to prove that the truncation maps are isomorphisms. The idea is to construct a cofinal sequence of Hamiltonians $(H_l)\subset Ad(V)\cap Ad(W{\subset}V)$ for which we can directly show
\begin{equation}\label{proofinvtheo}
 \begin{aligned}
  SH^{\geq0}_\ast(W{\subset}V) = \varinjlim_{l\rightarrow\infty} FH^{>0}_\ast(H_l)&\overset{(1)}{\simeq}\varinjlim_{l\rightarrow\infty} FH_\ast(H_l)=SH_\ast(V)\\
 SH_{\geq0}^\ast(W{\subset}V) = \varprojlim_{l\rightarrow\infty} FH_{>0}^\ast(H_l) &\overset{(2)}{\simeq}\varprojlim_{l\rightarrow\infty} FH^\ast(H_l) = SH^\ast(V).
 \end{aligned}
\end{equation}
We assume that all closed Reeb orbits of the contact form $\alpha=\lambda|_{\partial W}$ on $\partial W$ are transversely non-degenerate and that the attaching area does not intersect any closed Reeb orbit. This is generically satisfied (see \cite{FraCie}, appendix B).\\
Now fix sequences $(\mathfrak{a}_l), (\veps_l)\subset\mathbb{R}$, such that $\mathfrak{a}_l\not\in spec(\partial W,\alpha)$, $\mathfrak{a}_{l+1}\geq\mathfrak{a}_l$, $\mathfrak{a}_l\rightarrow\infty$ and $\veps_l\rightarrow 0$. Then choose an increasing sequence of non-degenerate Hamiltonians $H_l$ on $W$ that is on $\partial W\times (-\veps_l,0]$ of the form
\[H_l|_{\partial W\times(-\veps_l,0]} = \mathfrak{a}_l{\cdot} e^r-\big(1{+}\veps_l\big)\]
and extend $H_l$ over the handle by a function $\psi$ with $\mathfrak{a}=\mathfrak{a}_l$ and $\mathfrak{b}=-(1{+}\veps_l)$ as described in \ref{handle} and \ref{ExPsi}. For each $l$ choose the handle $\mathcal{H}_l$ so thin such that each trajectory of $X_{H_l}$ which leaves and reenters the handle has length greater than 1. This is possible as the attaching sphere is isotropic with $\dim S=k<n=\frac{1}{2}\dim V$ and can hence be chosen to have no Reeb chords (see \ \ref{appendixB}).\pagebreak[1]\\
Note that $\widehat{W{\cup}\mathcal{H}_l}$ is the same symplectic manifold for all choices of $\mathcal{H}_l$, as choosing different handles is to be understood as choosing different parametrizations on the cylindrical part of the completion. Moreover, if $\mathcal{H}_l$ is thinner than $\mathcal{H}_{l-1}$, then $\partial(W{\cup}\mathcal{H}_l)$ lies inside $W{\cup}\mathcal{H}_{l-1}$ (see Rem.\ \ref{remarksonpsi}). This guarantees that $H_l\geq H_{l-1}$ everywhere, since inside $W$ this holds by assumption, on the handle this holds by construction and on the cylindrical part of the symplectization the slope of $H_l$ in coordinates adapted to $\partial(W{\cup}\mathcal{H}_{l-1})$ is at least $\mathfrak{a}_l>\mathfrak{a}_{l-1}$, as $\partial(W{\cup}\mathcal{H}_l)$ lies inside $W{\cup}\mathcal{H}_{l-1}$. Thus we obtain a cofinal admissible sequence $(H_l)$, whose 1-periodic orbits having positive action are all contained in $W$. Recall from \ref{sectrunc} the long exact sequences
\begin{align*}
 \dots\rightarrow FH^{\geq 0}_{j+1}(H_l)\rightarrow FH^{<0}_j(H_l)\rightarrow FH_j(H_l)\rightarrow FH_j^{\geq 0}(H_l)\rightarrow\dots\\
 \dots\rightarrow FH_{< 0}^{j-1}(H_l)\rightarrow FH_{\geq0}^j(H_l)\rightarrow FH^j(H_l)\rightarrow FH^j_{< 0}(H_l)\rightarrow\dots
\end{align*}
and note that $FH^{\geq0}_j(H_l)$ is generated by all 1-periodic orbits of $H_l$ inside $W$, while $FH^{< 0}_j(H_l)$ is generated by all other orbits. The orbits of negative action all lie on the handle and are explicitly given in (\ref{eqA}). They are all transversely non-degenerate and their Morse-Bott indices are given by $(n{-}k)\big(\frac{\mathfrak{a}_l C_z}{\pi}{-}1\big)+\left\lbrace\begin{smallmatrix}1\\ 2(n-k)\end{smallmatrix}\right.$. It follows that the possible values of $\mu(\gamma)$ increase to $\infty$ as the slope $\mathfrak{a}=\mathfrak{a}_l$ tends to $\infty$ (see \ref{sec5.4}). Therefore, $FH^{<0}_j(H_l)$ becomes eventually zero as $l$ increases and so does $FH^{<0}_{j+1}(H_l)$. This implies that $FH_j(H_l)\rightarrow FH^{\geq 0}_j(H_l)$ is an isomorphism for $l$ large enough. As the direct limit is an exact functor, these maps converge to an isomorphism in the limit, proving (\ref{proofinvtheo})(1).\\
For cohomology, the same arguments apply. Although the inverse limit is only left exact, it still takes the isomorphism $FH^j_{\geq 0}(H_n)\rightarrow FH^j(H_n)$ to an isomorphism in the limit (see \cite{eilenberg}, Thm.\ 5.4 or \cite{bourbaki2}, $\S6$, no.3, prop.\ 4). This proves (\ref{proofinvtheo})(2).
\end{proof}
\begin{proof}[\textnormal{\textbf{Proof of Corollary \ref{corinvsur}}}]~\\
 We have to describe $SH^+(V)$ in terms of $SH^+(W)$, if $V$ is obtained from $W$ by attaching a subcritical handle $\mathcal{H}$. Recall from Sec.\ \ref{sectransfer} that the transfer maps can be constructed to respect the action filtration. Thus, we find that they fit into the following commutative diagram of long exact sequences:
 \[\begin{aligned}\begin{xy}\xymatrix{\rightarrow H_{\ast+n}(V,\partial V)\ar[d]^{j^!_{\ast+n}} \ar[r] & SH_\ast(V)\ar[d]^{\pi_\ast(W,V)}\ar[r] & SH_\ast^+(V)\ar[d]^{\pi^+_\ast(W,V)} \ar[r] & H_{\ast-1+n}(V,\partial V) \ar[d]^{j^!_{\ast-1+n}} \rightarrow\\
 \rightarrow H_{\ast+n}(W,\partial W) \ar[r] & SH_\ast(W)\ar[r] & SH_\ast^+(W) \ar[r] & H_{\ast-1+n}(W,\partial W) \rightarrow }\end{xy}\end{aligned}\tag{$\ast$}\]
 The horizontal sequences are obtained by action filtration with $SH^0_\ast(V)\cong H_{\ast+n}(V,\partial V)$. The vertical maps are the transfer map $\tilde{\pi}_\ast(W,V)=\pi_\ast(W,V)$ or induced by the transfer map on the filtered complexes (see proof of Prop.\ \ref{prop2transmaps} and \cite{Vit}, Thm.\ 3.1 or \cite{Ritter}, Thm.\ 9.5 for more details).\pagebreak[1]\\
 In particular, $\tilde{\pi}_\ast(W,V): SH^0_\ast(V)\rightarrow SH^0_\ast(W)$ is the shriek map $j_{\ast+n}^!$ defined as follows. The inclusion $j:W\hookrightarrow V$ induces in singular cohomology the map $j^\ast: H^\ast(V)\rightarrow H^\ast(W)$ which induces $j^!_\ast: H_{2n-\ast}(V,\partial V)\rightarrow H_{2n-\ast}(W,\partial W)$ via the following commutative diagram using the Poincaré duality $PD$:
 \[\begin{xy}\xymatrix{H^\ast(V) \ar[d]^{j^\ast} \ar[r]^{\mspace{-40mu}\cong}_{\mspace{-40mu}PD} & H_{2n-\ast}(V,\partial V)\ar[d]^{j^!_\ast}\\ H^\ast(W) \ar[r]^{\mspace{-40mu}\cong}_{\mspace{-40mu}PD} & H_{2n-\ast}(W,\partial W).}\end{xy}\]
 As $V=W\cup\mathcal{H}$, it follows from excision that $H^\ast(V,W)\cong H^\ast(\mathcal{H},\partial^-\mathcal{H})$, where $\partial^-\mathcal{H}$ is the part of $\partial\mathcal{H}$ that is glued to $W$. As $\mathcal{H}$ is a 1-handle, the quotient $\mathcal{H}\big/\partial^-\mathcal{H}$ retracts to $S^1$. Using these facts with the long exact sequence for the pair $(V,W)$ and the assumption that $W$ has one connected component more than $V$, it is easy to show that $j^\ast$ and hence $j^!_{2n-\ast}$ are isomorphisms for $\ast\neq 0$ and injective for $\ast=0$.\\
 By Thm.\ \ref{theoinvsur}, we know that $\pi_\ast(W,V)$ is an isomorphism for all $\ast$. Hence applying the five-lemma (cf.\ \cite{hatcher}, Chap. 2.1) to $(\ast)$ shows that $\pi^+_\ast(W,V)$ is an isomorphism for $\ast\neq n,n{+}1$. For the remaining cases, we look at the following section of $(\ast)$:
 \[\begin{xy}\xymatrix@C-16pt@M=1pt{0\ar[r] & SH_{n+1}(V)\ar[d]^{\cong}\ar[r] & SH_{n{+}1}^+(V)\ar[d]^{inj} \ar[r] & \mathbb{Z}_2^{k-1} \ar[d]^{inj} \ar[r] & SH_n(V) \ar[d]^\cong \ar[r] & SH_n^+(V)\ar[d]^{surj} \ar[r] & H_{2n{-}1}(V,\partial V)\ar[d]^{\cong} \ar[r]^{\mspace{5mu}i_V}& SH_{n{-}1}(V)\ar[d]^\cong\\
 0 \ar[r] & SH_{n+1}(W)\ar[r] & SH_{n+1}^+(W) \ar[r] & \mathbb{Z}_2^k \ar[r] & SH_n(W) \ar[r] & SH_n^+(W) \ar[r] & H_{2n-1}(W,\partial W) \ar[r]^{\mspace{20mu}i_W} & SH_{n-1}(W)}\end{xy}\]
 Here, $k$ resp.\ $k{-}1$ are the number of connected components of $W$ resp.\ $V$, so that $H_{2n}(V,\partial V)\cong\mathbb{Z}_2^{k-1}$ and $H_{2n}(W,\partial W)=\mathbb{Z}^k_2$. Note that $H_{2n+1}(V,\partial V)=H_{2n+1}(W,\partial W)=0$ due to dimensions, while $\pi_n^+(W,V)$ surjective,  $\pi_{n+1}^+(W,V)$ injective follows from the five-lemma. The usual count of ranks gives 
 \begin{align*}
  0 &= \rk SH_{n+1}(V)-\rk SH_{n+1}^+(V)+1-\rk SH_n(V)+\rk SH_n^+(V)-\rk \ker i_V\\
  0 &= \rk SH_{n+1}(W)-\rk SH_{n+1}^+(W)+2-\rk SH_n(W)+\rk SH_n^+(W)-\rk \ker i_W.
 \end{align*}
As $\ker i_V\cong \ker i_W$ by commutativity and $SH_\ast(V)\cong SH_\ast(W)$, we obtain
\[\rk SH_{n+1}^+(W)-\rk SH_{n+1}^+(V)=\rk SH_n^+(W)-\rk SH_n^+(V)+1.\tag{$\ast\ast$}\]
$\pi_n^+(W,V)$ being surjective and $\pi_{n+1}^+(W,V)$ being injective yields the estimates
\[\rk SH_{n+1}^+(W)-\rk SH_{n+1}^+(V)\geq 0\qquad\text{ and }\qquad \rk SH_n^+(V)-\rk SH_n^+(W)\geq 0,\]
which yields together with $(\ast\ast)$ that
\begin{itemize}
 \item either
$\qquad\rk SH_{n+1}^+(W)=\rk SH_{n+1}^+(V)$ and $\rk SH_n^+(V)=\rk SH_n^+(W)\,{+}1$
\item or $\qquad\rk SH_{n+1}^+(W)\,{-}1=\rk SH_{n+1}^+(V)$ and $\rk SH_n^+(V)=\rk SH_n^+(W)$.
\end{itemize}\pagebreak[1]
\end{proof}
\begin{proof} [\textnormal{\textbf{Proof of Proposition \ref{propcofinalHamilt}}}]~\\
 We have to show that if $(\Sigma,\xi)$ with contact form $\alpha$ is a.f.g.\ in degree $j$ with bound $b_j(\xi)$, shown by sequences $f_l:\Sigma\rightarrow \mathbb{R}, (\mathfrak{a}_l)\subset\mathbb{R}$, then $(\widehat{\Sigma},\hat{\xi})$ obtained by contact surgery along an isotropic $(k{-}1)$-sphere $S$ is also a.f.g., provided that $S$ has no $R_l$-Reeb chords shorter than $\mathfrak{a}_l$.\\
 By Rem.\ \ref{remarksafg}, we may assume that $f_l(p)\leq f_{l-1}(p)\;\forall\, p\in\Sigma$ and $\mathfrak{a}_l\geq\mathfrak{a}_{l-1}$. Let ${W=\big((-\infty,0]{\times}\Sigma,e^r{\cdot}\alpha\big)}$ denote the negative symplectization of $(\Sigma,\alpha)$. As described in \ref{secconsur}, contact surgery along $S$ is realized by attaching a $k$-handle $\mathcal{H}$ along $S$ to $\partial W$. Let $\widehat{W}$ denote the completion of $W$, that is exactly the whole symplectization of $(\Sigma,\alpha)$. As in \ref{secsetup}, consider the subset $W_l$ of $\widehat{W}$ defined by $W_l{:=}\big\{(r,p)\,\big| r{\leq} f(p), p{\in}\Sigma\big\}$. We note that $f_l(p){\leq} f_{l-1}(p){\leq} 0\;\forall\,p$ implies $W_l{\subset} W_{l-1}{\subset} W$ and that $\partial W_l{=}\Sigma_l{=}\big\{(f(p),p)\,\big|\,p\in\Sigma\big\}$.\\
 For every $l$ choose the handle $\mathcal{H}_l$ attached to $W_l$ so thin, such that each Reeb trajectory that leaves the attaching area and reenters it later is longer then $\mathfrak{a}_l$. This is possible as the attaching sphere $S$ has no Reeb chords shorter than $\mathfrak{a}_l$. Moreover, choose $\mathcal{H}_l$ so thin that $W_l{\cup}\mathcal{H}_l$ lies inside $W_{l-1}{\cup}\mathcal{H}_{l-1}$ (see Rem.\ \ref{dishandlesdifferentcontact} for the interpretation of $W_l{\cup}\mathcal{H}_l$ as a subset of $\widehat{W_{l-1}{\cup}\mathcal{H}_{l-1}}$). Note that $\widehat{W_l{\cup}\mathcal{H}_l}$ and $\widehat{W{\cup}\mathcal{H}}$ are easily identified since $\partial\big(W_l{\cup}\mathcal{H}_l\big)$ can be understood as a contact hypersurface inside the negative symplectization of $\partial\big(W{\cup}\mathcal{H}\big)$ (see Rem.\ \ref{remarksonpsi} and Rem.\ \ref{dishandlesdifferentcontact}).\pagebreak[3]\\
 If we construct $\mathcal{H}_l$ with the help of a function $\psi$ as in Sec.\ \ref{handle} and \ref{ExPsi}, then there are no closed Reeb orbits on $\mathcal{H}_l$ except on the cocore and closed Reeb orbits that go over $\mathcal{H}_l$ and into $\partial W_l$ are longer than $\mathfrak{a}_l$ if $\mathcal{H}_l$ is sufficiently thin.\\
 Let $\beta_l$ denote the contact form on $\partial\big(W_l{\cup}\mathcal{H}_l\big)$. We find that $\beta_l$ agrees with $\alpha_l$ away from $\mathcal{H}_l$. Moreover, the closed Reeb orbits of $\beta_l$ shorter than $\mathfrak{a}_l$ are whose of $\alpha_l$ plus the ones on the cocore of the handle. Of the latter there exists by (\ref{indexonhandle}) at most one, if the degree $j$ equals $(n{-}k)N+\big\lbrace\begin{smallmatrix}1\phantom{(n-k)}\\2(n-k)\end{smallmatrix}$ for $N>0$ odd. This gives the bound $b_j(\hat{\xi})=b_j(\xi)$ or $b_j(\hat{\xi})=b_j(\xi){+}1$.\\
 As $W_l{\cup}\mathcal{H}_l\subset W_{l-1}{\cup}\mathcal{H}_{l-1}\subset W{\cup}\mathcal{H}$, we finally find that $\beta_l$ written in coordinates of the symplectization of $\partial\big(W{\cup}\mathcal{H}\big)$ in the form $\beta_l=e^{g_l}{\cdot}\beta$ with $g_l:\partial\big(W{\cup}\mathcal{H}\big)\rightarrow\mathbb{R}$ gives
 \[g_l(p)\leq g_{l-1}(p)\leq 0\qquad\forall\; p\in \big(W{\cup}\mathcal{H}\big).\]
 Hence $(g_l)$ and $(\mathfrak{a}_l)$ show that the contact structure on $\big(W{\cup}\mathcal{H}\big)$ is a.f.g..
\end{proof}

\section{Brieskorn manifolds}\label{secBries}
\subsection{General results}\label{secGenBrie}
In this subsection, we recall the construction of Brieskorn manifolds, their contact structures and fillings and the Morse-Bott index on these manifolds. Then we present Ustilovsky's infinitely many different contact structures on $S^{4m+1}$ and introduce finally the special exotic contact structures on $S^{2n-1}$ announced in Prop.\ \ref{propTwoConStr}.\pagebreak[1]\bigskip\\
Let $a = (a_0,a_1, ..., a_n)\in\mathbb{N}^{n+1}$ with $a_i\geq 2$ and define the complex polynomial 
\[f(z) = z_0^{a_0}+ z_1^{a_1} + ... + z_n^{a_n}.\]
Its level sets $V_a(t):=f^{-1}(t)$ are smooth complex hypersurfaces except for $V_a(0)$, which has a single singularity at zero. The link of this singularity $\Sigma_a:=V_a(0)\cap S^{2n+1}$ is the \textit{Brieskorn manifold} $\Sigma_a$. It is proven in \cite{LuMe} that the following 1-form $\lambda_a$ on $\mathbb{C}^{n+1}$ restricts to a contact form $\alpha_a:=\lambda_a|_\Sigma$ on $\Sigma_a$ defining the contact structure $\xi_a:=\ker \alpha_a$, the Reeb vector field $R_a$ and the Reeb flow $\varphi^t_a$:
\begin{align*}
 \lambda_a &=\frac{i}{8} \sum^n_{k=0} a_k (z_kd\bar{z}_k - \bar{z}_kdz_k),& R_a &= 4i\Big( \frac{z_0}{a_0},...,\frac{z_n}{a_n}{\Big)},\\ \varphi^t_a(z)&=\big(e^{4it/a_0} z_0,...,e^{4it/a_n} z_n{\big)}\mspace{-3mu}.
\end{align*}
 A Liouville filling $(V,\lambda)$ of $(\Sigma_a,\alpha_a)$ can be obtained as follows. Choose a monotone decreasing $\beta\in C^\infty (\mathbb{R})$ with $\beta(x)= 1$ for $x\leq \frac{1}{4}$ and  $\beta(x)=0$ for $x\geq \frac{3}{4}$ and define
\[ V_a:=V_\veps\cap B_1(0), \quad\text{ where }\quad V_\veps:=\Big\{z\in \mathbb{C}^{n+1}\,\Big|\, z_0^{a_0}+z_1^{a_1}+...+z_n^{a_n} = \veps\cdot\beta\big(||z||^2\big)\Big\}.\]
For $\veps$ small enough, $(V_a,\lambda_a)$ is a Liouville domain with boundary $(\Sigma_a,\alpha_a)$ and vanishing first Chern class $c_1(TV)$ (see \cite{Fauck1}, Prop.\ 2.1.3 or \cite{FauckThesis}, Prop.\ 99). By \cite{Fauck1}, Lem.\ 2.5, the Morse-Bott index of a closed $R_a$-trajectory $v$ of period $t=L\pi/2$ is given as
\begin{equation}\label{eqConZeh}
 \mu(v)= 2\cdot \sum_{k=0}^n\Big\lceil\frac{L}{a_k}\Big\rceil-2L+\mu_{Morse}(v)-(n{-}1).
\end{equation}
Here, we assume $sign(h''){=}\,{+}1$ for the slope of the Hamiltonian. Note that the grading convention in \cite{Fauck1} differs from the one used here by $\frac{1}{2}\dim\Sigma_a$.

\begin{proof}[\textnormal{\textbf{Proof of Proposition \ref{propTwoConStr}}}]~\\
We provide for every $k\in\mathbb{Z}$ an a.f.g.\ contact manifold $(\Sigma_k,\xi_k)$ with exact filling $(V_k,\lambda_k)$ such that $\rk SH_k^+(V_k)\neq 0$ and $\Sigma_k$ is homeomorphic to $S^{2n-1}$. These manifolds can be chosen as Brieskorn manifolds $(\Sigma_a,\xi_a)$ with their filling $(V_a,\lambda_a)$ as described above. Explicitly, we take 
\begin{align*}
 &&\Sigma_k^+:=\Sigma_k&=\Sigma_a \quad\text{ with } a=(2,2,2,a_3,....,a_n) &&\text{ for } k\geq n-1\\
 &\text{and }& \Sigma_k^-:=\Sigma_k&=\Sigma_a \quad\text{ with } a=(2,3,a_2,a_3,....,a_n) &&\text{ for } k\leq n-2,
\end{align*}
where $a_2<a_3<...<a_n$ are odd integers which are pairwise coprime and, depending on $k$, sufficiently large. The sign on $\Sigma^\pm_k$ just distinguishes the two cases within this proof. By \cite{Brie}, Satz 1, these Brieskorn manifolds are homeomorphic to $S^{2n-1}$. As the Reeb flow on Brieskorn manifolds is periodic, all $(\Sigma_k,\xi_k)$ are a.f.g.\ in every degree by Prop.\ \ref{PropafgExists}. It remains to calculate $SH_k^+(V_k)$. For that, recall the Reeb flow $\varphi^t_a$ on Brieskorn manifolds
\begin{align*}
 \text{on $\Sigma^+_k$ : } \varphi^t_a(z)&=\big(e^{2it}\cdot z_0,e^{2it}\cdot z_1,e^{2it}\cdot z_2,e^{4it/a_3}\cdot z_3,...,e^{4it/a_n}\cdot z_n\big)\\
 \text{on $\Sigma^-_k$ : } \varphi^t_a(z)&=\big(e^{2it}\cdot z_0,e^{4/3\cdot it}\cdot z_1,e^{4it/a_2}\cdot z_2,...,e^{4it/a_n}\cdot z_n\big).
\end{align*}
If the period $t=L\pi/2$ is smaller than $a_3\pi/2$, we find in the first case, $k\geq n{-}1$, that there are periodic orbits only for $L\in2\mathbb{Z}$ forming the manifolds
\[\mathcal{N}^L=\Big\{z\in\Sigma_k^+\Big|z_j{=}0, j{\geq}3\Big\}\cong\Big\{\mspace{-2mu}z=(z_0,z_1,z_2)\Big|z_0^2+z_1^2+z_2^2=0, \big|z\big|^2=1\Big\}\cong S^\ast S^2.\]
On $S^\ast S^2$, there exists a Morse function with exactly 4 critical points, having Morse indices 0,1,2,3. Let us denote these critical points on $\mathcal{N}^L$ by $L\gamma_0, L\gamma_1, L\gamma_2, L\gamma_3$. By (\ref{eqConZeh}), we find that their Morse-Bott index is for $L<a_3$ exactly
\[\mu(L\gamma_c)=\underbrace{3L}_{a_0,a_1,a_2=2}+\underbrace{2(n{-}2)}_{a_3,...,a_n}-2L+c-(n{-}1)=L+n-3+c.\]
The distribution of the critical points among the degrees $j$ can be visualized as
\begin{center}
\begin{tabular}{c||c|c|c|c|c|c|c|c}
 $j$ & $\mspace{-3mu}{<}n{-}1\mspace{-3mu}$ & $\mspace{-3mu}n{-}1\mspace{-3mu}$ & $\mspace{-3mu}n\mspace{-3mu}$ & $\mspace{-3mu}n{+}1\mspace{-3mu}$ & $\mspace{-3mu}n{+}2\mspace{-3mu}$ &  $\mspace{-3mu}$...$\mspace{-3mu}$ & $\mspace{-3mu}n{-}1+2l\mspace{-3mu}$ & $\mspace{-3mu}n{-}1+2l+1$ \\\hline
 $L\gamma_c$ & -- & $\mspace{-3mu}2\gamma_0\mspace{-3mu}$ & $\mspace{-3mu}2\gamma_1\mspace{-3mu}$ & $\mspace{-3mu}2\gamma_2, 4\gamma_0\mspace{-3mu}$ & $\mspace{-3mu}2\gamma_3, 4\gamma_1\mspace{-3mu}$ &  $\mspace{-3mu}$...$\mspace{-3mu}$ & $\mspace{-3mu}2l\gamma_2, (2l{+}2)\gamma_0\mspace{-3mu}$ & $\mspace{-3mu}2l\gamma_3, (2l{+}2)\gamma_1$
\end{tabular}
\end{center}
The first irregularity occurs for $a_3=L$. If $L\geq a_3$, we can estimate the index by
\[\mu(v)= 2 \sum_{k=0}^n\Big\lceil\mspace{-2mu}\frac{L}{a_k}\mspace{-2mu}\Big\rceil-2L+\mu_{Morse}(v)-(n{-}1)\geq3L+2(n{-}2)-2L+0-(n{-}1)\geq a_3+n{-}3.\]
So for $a_3,...,a_n$ sufficiently large, we find that $FC_k(V_k)$ is generated by 2 elements (as given in the table) and by $2\gamma_0$ resp.\ $2\gamma_1$ for $k=n{-}1$ resp.\ $k=n$. It is shown in the author's thesis, \cite{FauckThesis} Thm.\ 109, that the boundary operator between these generators is zero with respect to $\mathbb{Z}_2$-coefficients. Actually, the result was show for Rabinowitz-Floer homology, but the same argument can also be used to deduce this result for $SH^+$. Explicitly, the symplectic symmetry $(z_0,z_1,z_2,z_3,...,z_n)\mapsto (-z_0,-z_1,z_2,z_3,...,z_n)$ can be used to show that all Floer trajectories between these generators come in pairs and are therefore counted as 0 in $\mathbb{Z}_2$-coefficients (see also \cite{Ueb}, 3.3). Hence we can conclude
\[SH_k^+(V_k)\cong\begin{cases}(\mathbb{Z}_2)^2 &\text{ for } k\geq n+1\\\;\;\mathbb{Z}_2 &\text{ for }k=n{-}1, n\end{cases}.\]
In the second case, $k<n{-}1$, we find by the same reasoning that if $t=L\pi/2<a_2\pi/2$, then there are periodic orbits only for $L\in 6\mathbb{Z}$, forming the manifolds 
\[\mathcal{N}^L =\Big\{z\in\Sigma_k^-\Big|z_j=0, j\geq 2\Big\}\cong\Big\{(z_0,z_1)\Big|z_0^2+z_1^3=0, \big|(z_0,z_1)\big|^2=1\Big\}\cong S^1.\]
On $S^1$, there exists a Morse function with exactly 2 critical points with Morse index 0 or 1. Let them be denoted by $L\gamma_0$ and $L\gamma_1$. Again by (\ref{eqConZeh}), their Morse-Bott index is for $L<a_2$ given by
\[\mu(L\gamma_c)=\underbrace{2\cdot \frac{L}{2}}_{a_0}+\underbrace{2\cdot \frac{L}{3}}_{a_1}+\underbrace{2(n{-}1)}_{a_2,...,a_n}-2L+c-(n{-}1)=-\frac{L}{3}+(n{-}1)+c.\]
The distribution of these critical points among the degrees $j$ is as follows
\begin{center}
 \begin{tabular}{c|c|c|c|c|c|c|c|c|c}
   $j$ & ... & $n{-}1\,{-}2l$ & $n{-}1\,{-}2l\,{+}1$ & ... & $n{-}5$ & $n{-}4$ & $n{-}3$ & $n{-}2$ & $>n{-}2$\\\hline
   $L\gamma_c$ & ... & $6l\gamma_0$ & $6l\gamma_1$ & ... & $12\gamma_0$ & $12\gamma_1$ & $6\gamma_0$ & $6\gamma_1$ & --
  \end{tabular}\pagebreak[2]
\end{center}
The first irregularity occurs for $L=a_2$. If $L\geq a_2$, we can estimate the index by
\begin{align*}
 \mu(v)&= 2 \sum_{k=0}^n\Big\lceil\mspace{-2mu}\frac{L}{a_k}\mspace{-2mu}\Big\rceil\,{-}2L\,{+}\mu_{Morse}(v)\,{-}(n{-}1)&
 &\leq 2\sum_{k=0}^n\Big(\frac{L}{a_k}\,{+}1\Big)\,{-}2L\,{+}2n{-}1\,{-}(n{-}1)\\
 &= 2L\left(\sum_{k=0}^n\frac{1}{a_k}-1\right)+3n+2&
 &\mspace{-7mu}\overset{\substack{a_0=2\\a_1=3}}{=}2L\left(\sum_{k=2}^n\frac{1}{a_k}-\frac{1}{6}\right)+3n+2.
\end{align*}
Now if $\sum_{k=2}^n\frac{1}{a_k}<\frac{1}{6}$ and if $a_2,...,a_n$ are sufficiently large, then we find for all degrees $j>2a_2\big(\sum_{k=2}^n\frac{1}{a_k}{-}\frac{1}{6}\big){+}3n{+}2$ that $FC_j^+(V_k)$ is generated by the unique critical point listed in the table. Again, the boundary operator is zero in all these cases. This can be read off the table as whenever two critical points have index difference one, then either both lie on the same critical manifold $\mathcal{N}^L\cong S^1$, or the one with the bigger index has smaller action. In the first case the boundary operator is zero as the Morse boundary operator on $S^1$ is zero and in the second case there are no connecting Floer trajectories $u$, as the action has to decrease along $u$. Hence we conclude
\[SH_k^+(V_k)\cong \mathbb{Z}_2\qquad \text{ for } k\leq n{-}2.\]
\end{proof}
Let us leave the general picture and turn to the Brieskorn manifolds, where
\[a=(2,...,2,\,p)\in\mathbb{N}^{n+1},\; p \text{ odd and } n=2m{+}1 \text{ odd.}\]
For fixed $p$, we write $\Sigma_p, \xi_p, \alpha_p,\lambda_p, V_p$ instead of $\Sigma_a, \xi_a, \alpha_a, \lambda_a, V_a$. In \cite{Brie}, Satz 1, Brieskorn showed that $\Sigma_p$ is a $4m{+}1$-dimensional homotopy sphere and $\Sigma_p$ is diffeomorphic to $S^{4m+1}$ if and only if $p\equiv \pm 1\bmod 8$ and diffeomorphic to the Kervaire sphere otherwise. For $p\equiv \pm 1\bmod 8$, the $(\Sigma_p,\xi_p)$ are exactly the exotic contact structures on $S^{4m+1}$ given by Ustilovsky in \cite{Usti}. On $\Sigma_p$, the Reeb flow $\varphi_p^t$ is 
\[\varphi^t_p(z)=\left(e^{2it}{\cdot} z_0,\,...\,,e^{2it}{\cdot} z_{n-1}\,,\,e^{4it/p}{\cdot} z_n\right).\]
Writing $t=L\pi/2$, we see that $\varphi^t_p$ has periodic orbits exactly for $L\in 2\mathbb{Z}$. Note that we do not have closed orbits of period $t{=}k\frac{p}{2}\pi\Leftrightarrow L{=}kp$, since any $z\in\Sigma_p$ has at least two non-zero components. The closed orbits with period $t=L\pi/2$ form manifolds $\mathcal{N}^L$ of two types:
\begin{itemize}
 \item if $p\nmid L$, then $\mathcal{N}^L=\big\{z{\in}\Sigma_p\,\big|\,z_n{=}0\big\}=\big\{\,z{\in}\mathbb{C}^{n+1}\,\big|\,z_n{=}0,\,\sum z_k^2{=}0,\,||z||^2{=}1\big\}$, which is diffeomorphic to the unit cotangent bundle $S^\ast S^{n-1}$.
 \item if $p\;{\mid}\, L$, then $\mathcal{N}^L=\Sigma_p$, which is for $p\equiv \pm 1\bmod 8$ diffeomorphic to $S^{2n-1}$.
\end{itemize}
On $S^\ast S^{n-1}$ there exists a Morse function with exactly 4 critical points of Morse indices $0, n{-}2,n{-}1, 2n{-}3$. On $S^{2n-1}$ there exists a Morse function having exactly 2 critical points of Morse index 0 and $2n{-}1$. We denote these critical points by
\begin{align*}
 &L\gamma_0,\; L\gamma_{n-2},\; L\gamma_{n-1},\; L\gamma_{2n-3}\; &&\text{ on $\mathcal{N}^L\cong S^\ast S^{n-1} $ if $ p\nmid L$}\\
 &L\gamma_0,\; L\gamma_{2n-1}&&\text{ on $\mathcal{N}^L\cong S^{2n-1}\phantom{S} $ if $ p\mid L$.} 
\end{align*}
With (\ref{eqConZeh}) and $L\in2\mathbb{Z}$, we find that the Morse-Bott index $\mu(L\gamma_c)$ is given by
\begin{align*}
 &&\mu(L\gamma_c)&= nL+2\big\lceil L/p\big\rceil-2L+c-(n{-}1)=(n{-}2)(L{-}1)+2\big\lceil L/p\big\rceil+c-1\\
 &\Rightarrow&\mu(L\gamma_c)&=
\begin{cases}
 (L{-}1)(n{-}2)+2\big\lceil \textstyle\frac{L}{p}\big\rceil\;\;\;-1 & \text{for }c=0\\
 \phantom{(-1} L\phantom{)}(n{-}2)+2\big\lceil \textstyle\frac{L+1}{p}\big\rceil-1 & \text{for } c=\;\,n{-}2 \text{ and }p\nmid L\\
 \phantom{(-1} L\phantom{)}(n{-}2)+2\big\lceil \textstyle\frac{L+1}{p}\big\rceil & \text{for } c=\;\,n{-}1 \text{ and }p\nmid L\\
 (L{+}1)(n{-}2)+2\big\lceil \textstyle\frac{L+2}{p}\big\rceil & \text{for } c=2n{-}3 \text{ and }p\nmid L, L{+}1\\
 (L{+}1)(n{-}2)+2\big\lceil \textstyle\frac{L+2}{p}\big\rceil-2 & \text{for } c=2n{-}3 \text{ and }p\mid L{+}1\\
 (L{+}1)(n{-}2)+2\big\lceil \textstyle\frac{L+2}{p}\big\rceil & \text{for } c=2n{-}1 \text{ and }p\mid L
\end{cases}
\end{align*}
Define $f_p:\mathbb{Z}\rightarrow\mathbb{Z}$ to be the strictly increasing function $f_p(l)=(l{-}1)(n{-}2)+2\big\lceil l/p\big\rceil.$
\begin{theorem}[cf.\ \cite{Fauck1}, Thm. 2.7]\label{theoSHofSigmap}~\\
 For $n\geq 5$ and $k\geq n$ holds that $SH_k(V_p,\Sigma_p) = 0$, except
 \begin{align*}
 &\rk SH_k(V_p,\Sigma_p) = 1&&\text{ if }\quad k=f_p(l)+\left\lbrace\begin{smallmatrix}-1\\\phantom{-}0\end{smallmatrix}\right.&&\text{ for } l\in\phantom{2}\mathbb{N} \text{ and }p\nmid l{-}1,\\
 &\rk SH_k(V_p,\Sigma_p) = \;?&&\text{ if }\quad k=f_p(l)+\left\lbrace\begin{smallmatrix}-2\\-1\end{smallmatrix}\right.&&\text{ for } l\in2\mathbb{N} \text{ and } p\mid l{-}1.
\end{align*}
The ranks at the place of ? are not known, but are either 0 or 1.
\end{theorem}
\begin{proof}
 It follows directly from the above calculations, that the distribution of the critical points $L\gamma_c$ among the degrees $j$ looks as follows:
 \begin{center}
  \begin{tabular}{c||c|c|c|c|c|c|c|c|c|c}
   $j$ & $\mspace{-3mu}$...$\mspace{-3mu}$ & $\mspace{-3mu}f_p(l)\text{-}1\mspace{-3mu}$ & $\mspace{-3mu}f_p(l)\mspace{-3mu}$ & $\mspace{-3mu}$...$\mspace{-3mu}$ & $\mspace{-3mu}f_p(l{+}1)\text{-}1\mspace{-3mu}$ & $\mspace{-3mu}f_p(l{+}1)\mspace{-3mu}$ & $\mspace{-3mu}$...$\mspace{-3mu}$ & $\mspace{-3mu}f_p(l{+}2)\text{-}1\mspace{-3mu}$ & $\mspace{-3mu}f_p(l{+}2)\mspace{-3mu}$ & ...\\\hline
   $L\gamma_c$ & $\mspace{-3mu}0\mspace{-3mu}$ & $\mspace{-3mu}l\gamma_0\mspace{-3mu}$ & $\mspace{-3mu}(l\text{-}2)\gamma_{2n-3}\mspace{-3mu}$ & $\mspace{-3mu}0\mspace{-3mu}$ & $\mspace{-3mu}l\gamma_{n-2}\mspace{-3mu}$ & $\mspace{-3mu}l\gamma_{n-1}\mspace{-3mu}$ & $\mspace{-3mu}0\mspace{-3mu}$ & $\mspace{-3mu}(l{+}2)\gamma_0\mspace{-3mu}$ & $\mspace{-3mu}l\gamma_{2n-3}\mspace{-3mu}$ & 0
  \end{tabular}
 \end{center}
Here, $l\in 2\mathbb{Z}$ is even and $p\nmid l, l{+}1$. If $p\mid l$, then the middle block in the table is missing and the last entry is $l\gamma_{2-1}$ instead of $l\gamma_{2n-3}$ (but with the same index!). If $p\mid l{+}1$, then the last entry $l\gamma_{2n-3}$ has index $f_p(l{+}2){-}2$ instead of $f_p(l{+}2)$. As in the proof of Prop.\ \ref{propTwoConStr}, we can directly see that the boundary operator between these generators of $FC_\ast(V_p)$ is almost always zero. In fact for any pair of $L\gamma_c$ with index difference one, either both lie on the same critical manifold $\mathcal{N}^L$ or the one with bigger index has lower action. The only exceptions are the cases where $p\mid l{+}1$. As the assumption $k\geq n$ guarantees that $FC_k(V_p)$ is only generated by Reeb trajectories, we conclude that $\rk SH_k(V_p)$ is exactly as claimed.
\end{proof}

\subsection{Connected sums of Brieskorn spheres}\label{secConSum}
We saw that the Reeb flow $\varphi_p^t$ of $\alpha_p$ on $\Sigma_p$ is periodic. Unfortunately, we need for the contact surgery construction\footnote{The contact surgery construction can of course always be done, without any conditions.} in this article at least one point that does not lie on a closed Reeb orbit. For this purpose, we will perturb $\alpha_p$ to a new contact form $\alpha_p'$ defining the same contact structure $\xi_p$. We will use a perturbation due to Uebele, \cite{Ueb}, which is inspired by the one used by Ustilovsky in \cite{Usti}. This perturbation is different from the one described in the appendix. The advantage of the perturbation here is that the resulting bounds $b_k(\xi_p')$ on the number of generators for Symplectic Homology are smaller. First, we make the following change of coordinates
\begin{equation}\label{UstiPer}
\begin{pmatrix}w_0\\w_1\end{pmatrix}=\frac{1}{\sqrt{2}}\begin{pmatrix}1 & i\\1& -i\end{pmatrix}\begin{pmatrix}z_0\\z_1\end{pmatrix},\quad w_2=z_2,\quad\dots\quad,\quad w_n=z_n.
\end{equation}
In these coordinates $\displaystyle\Sigma_p=\Big\{w\in\mathbb{C}^{n+1}\,\Big|\, 2w_0w_1+w_2^2+...+w_{n-1}^2+w_n^p=0,\, ||w||^2=1\Big\}.$\bigskip\\
Next, we introduce a new contact form $\alpha_p':=\frac{1}{K} \alpha_p$, with $K(w):=||w||^2+\veps\big(|w_0|^2{-}|w_1|^2\big)$, where $\veps>0$ is a sufficiently small irrational number. As $\alpha_p'$ is obtained from $\alpha_p$ by multiplication with a positive function, they define the same contact structure. In fact, $\alpha_p'$ should be though of as the restriction of $\lambda_p$  to the hypersurface $\Sigma_p'\subset \mathbb{R}{\times}\Sigma_p$ inside $\widehat{V_p}$ defined by $\Sigma_p':=\Big\{\big({-}\log K(y),y\big)\;\Big|\;y\in\Sigma_p\Big\}$.\\
Ustilovsky shows in \cite{Usti}, Lem.\ 4.1, that the Reeb vector field of $\alpha_p'$ is
\[R_p'=\Big(\;2i(1{+}\veps)w_0, 2i(1{-}\veps)w_1,2iw_2,\,...,\,2iw_{n-1},\frac{4i}{p}w_n\Big).\]
Hence, the closed Reeb orbits form 4 different types of manifolds:
\begin{itemize}
 \item $\widetilde{\mathcal{N}}^L=\Big\{ w\in \Sigma_p\,\Big|\,w_0{=}w_1{=}w_n{=}0\Big\}\cong S^\ast S^{n-3}$ of $\frac{L\pi}{2}$-periodic orbits for $p\nmid L$
 \item $\widetilde{\mathcal{N}}^L=\Big\{ w\in \Sigma_p\,\Big|\,w_0{=}w_1{=}0\Big\}\cong S^{2n-5}$ of $\frac{L\pi}{2}$-periodic orbits for $p\mid L$
 \item $\widetilde{\mathcal{N}}^L_+=\Big\{ (w_0,0,...,0)\in\Sigma_p\Big\}\cong S^1$ of $\frac{L\pi}{2(1+\veps)}$-periodic orbits
 \item $\widetilde{\mathcal{N}}^L_-=\Big\{(0,w_1,0,...,0)\in\Sigma_p\Big\}\cong S^1$ of $\frac{L\pi}{2(1-\veps)}$-periodic orbits.
\end{itemize}
Note that the flow of $R_p'$ is no longer periodic. In particular points $w\in \Sigma_p$ with $w_0, w_1, w_2$ all non-zero do not lie on any closed Reeb orbit. The Conley-Zehnder index for any Reeb trajectory $v$ of $\alpha'_p$ of length $L\pi/2$ can be computed analogous to \cite{Usti}, Lem.\ 4.2, and is given by:
\begin{equation}\label{muCZperturbed}
 \mu_{CZ}(v)=\left\lfloor{\mspace{-3mu}\textstyle\frac{L(1+\veps)}{2}}\mspace{-3mu}\right\rfloor\mspace{-4mu}+\mspace{-4mu}\left\lceil{\mspace{-3mu}\textstyle\frac{L(1+\veps)}{2}}\mspace{-3mu}\right\rceil+\left\lfloor{\mspace{-3mu}\textstyle\frac{L(1-\veps)}{2}}\mspace{-3mu}\right\rfloor\mspace{-4mu}+\mspace{-4mu}\left\lceil{\mspace{-3mu}\textstyle\frac{L(1-\veps)}{2}}\mspace{-3mu}\right\rceil+\mspace{-3mu}\sum_{k=2}^{n-1}\bigg(\mspace{-6mu}\Big\lfloor{\mspace{-3mu}\textstyle\frac{L}{2}}\mspace{-3mu}\Big\rfloor\mspace{-4mu}+\mspace{-4mu}\Big\lceil{\textstyle\mspace{-3mu}\frac{L}{2}}\mspace{-3mu}\Big\rceil\mspace{-6mu}\bigg)+\left\lfloor{\mspace{-3mu}\textstyle\frac{L}{p}}\mspace{-3mu}\right\rfloor\mspace{-4mu}+\mspace{-4mu}\left\lceil{\textstyle\mspace{-3mu}\frac{L}{p}}\mspace{-3mu}\right\rceil-\Big(\mspace{-5mu}\big\lfloor \mspace{-3mu}L\mspace{-3mu}\big\rfloor\mspace{-4mu}+\mspace{-4mu}\big\lceil \mspace{-3mu}L\mspace{-3mu}\big\rceil\mspace{-5mu}\Big).
\end{equation}
We can choose Morse-functions on $\widetilde{\mathcal{N}}^L, \widetilde{\mathcal{N}}^L_\pm$ having critical points $L\gamma_c, L\gamma_c^\pm$ with indices $0, n{-}4, n{-}3, 2n{-}7$ for $p\nmid L$ or $0, 2n{-}5$ for $p\mid L$ on $\widetilde{N}^L$ or $0,1$ on $\widetilde{\mathcal{N}}^L_\pm$. For $\veps$ fixed and $L\veps< 1$, we can estimate the Gau{\ss} brackets in (\ref{muCZperturbed}) and calculate the Morse-Bott index as:
\begin{align*}
\shortintertext{$\bullet$ on $\widetilde{\mathcal{N}}^L_+$ for period $\frac{L}{1+\veps}{\cdot} \frac{\pi}{2}, L\in2\mathbb{N}$ and $c\in\{0,1\}$}
  \mu(L\gamma_c^+)&=\textstyle\left\lfloor\mspace{-3mu}\frac{L(1+\veps)}{2(1+\veps)}\mspace{-3mu}\right\rfloor\mspace{-5mu}+\mspace{-5mu}\left\lceil\mspace{-3mu}\frac{L(1+\veps)}{2(1+\veps)}\mspace{-3mu}\right\rceil+\left\lfloor\mspace{-3mu}\frac{L(1-\veps)}{2(1+\veps)}\mspace{-3mu}\right\rfloor\mspace{-5mu}+\mspace{-5mu}\left\lceil\mspace{-3mu}\frac{L(1-\veps)}{2(1+\veps)}\mspace{-3mu}\right\rceil+ \sum_{k=2}^{n-1}\left(\left\lfloor\mspace{-3mu}\frac{L}{2(1+\veps)}\mspace{-3mu}\right\rfloor\mspace{-5mu}+\mspace{-5mu}\left\lceil\mspace{-3mu}\frac{L}{2(1+\veps)}\mspace{-3mu}\right\rceil\right)\\
 &\phantom{\,+\,}\textstyle +\left\lfloor\mspace{-3mu}\frac{L}{p(1+\veps)}\mspace{-3mu}\right\rfloor\mspace{-5mu}+\mspace{-5mu}\left\lceil\mspace{-3mu}\frac{L}{p(1+\veps)}\mspace{-3mu}\right\rceil-\left(\left\lfloor \mspace{-3mu}\frac{L}{(1+\veps)}\mspace{-3mu}\right\rfloor\mspace{-5mu}+\mspace{-5mu}\left\lceil \mspace{-3mu}\frac{L}{(1+\veps)}\mspace{-3mu}\right\rceil\right)-\frac{1}{2}\dim\widetilde{\mathcal{N}}^L_++\frac{1}{2}+c\\
 &= \textstyle 2 \frac{L}{2}+\left(\frac{L}{2}{-}1{+}\frac{L}{2}\right)+(n{-}2)\left(\frac{L}{2}{-}1{+}\frac{L}{2}\right)+2\left\lceil\frac{L}{p}\right\rceil -1 -(L{-}1{+}L)+c\\
 &=(n{-}2)(L{-}1)+\textstyle 2\left\lceil\frac{L}{p}\right\rceil + \left\lbrace\begin{smallmatrix}-1\\\phantom{-}0\end{smallmatrix}\right.\allowdisplaybreaks
   \shortintertext{$\bullet$ on $\widetilde{\mathcal{N}}^L$ for period $L{\cdot}\frac{\pi}{2},\; L\in2\mathbb{N},\;p\mid L$ and $c\in\{0, 2n{-}5\}$}
 \mu(L\gamma_c) &= \textstyle \left(\frac{L}{2}{+}\frac{L}{2}{+}1\right)+\left(\frac{L}{2}{-}1{+}\frac{L}{2}\right)+(n{-}2)2\frac{L}{2}+2\left\lceil\frac{L}{p}\right\rceil -2L-(n{-}3)+c\\
 &=\begin{cases}
    \bigg.(n{-}2)(L{-}1)+\textstyle 2\left\lceil\frac{L}{p}\right\rceil + 1 & c=0\\
    (n{-}2)(L{+}1)+\textstyle 2\left\lceil\frac{L+2}{p}\right\rceil -2 & c=2n{-}5 \text{ as }p\mid L
   \end{cases}
   \allowdisplaybreaks
   \shortintertext{$\bullet$ on $\widetilde{\mathcal{N}}^L$ for period $L{\cdot}\frac{\pi}{2},\; L\in2\mathbb{N},\;p\nmid L$ and $c\in\{0,n{-}4,n{-}3, 2n{-}7\}$}
 \mu(L\gamma_c) &=\textstyle\left\lfloor\mspace{-3mu}\frac{L(1+\veps)}{2}\mspace{-3mu}\right\rfloor\mspace{-5mu}+\mspace{-5mu}\left\lceil\mspace{-3mu}\frac{L(1+\veps)}{2}\mspace{-3mu}\right\rceil+\left\lfloor\mspace{-3mu}\frac{L(1-\veps)}{2}\mspace{-3mu}\right\rfloor\mspace{-5mu}+\mspace{-5mu}\left\lceil\mspace{-3mu}\frac{L(1-\veps)}{2}\mspace{-3mu}\right\rceil+ \sum_{k=2}^{n-1}\Big(\mspace{-5mu}\left\lfloor\frac{L}{2}\right\rfloor\mspace{-5mu}+\mspace{-5mu}\left\lceil\frac{L}{2}\right\rceil\mspace{-5mu}\Big)\\
 &\phantom{\,+\,}\textstyle +\left\lfloor\frac{L}{p}\right\rfloor\mspace{-5mu}+\mspace{-5mu}\left\lceil\frac{L}{p}\right\rceil-\left(\left\lfloor L\right\rfloor\mspace{-5mu}+\mspace{-5mu}\left\lceil L\right\rceil\right)-\frac{1}{2}\dim\widetilde{\mathcal{N}}^L+\frac{1}{2}+c\\
 &= \textstyle \left(\frac{L}{2}{+}\frac{L}{2}{+}1\right)+\left(\frac{L}{2}{-}1{+}\frac{L}{2}\right)+(n{-}2)2\frac{L}{2}+2\left\lceil\frac{L}{p}\right\rceil\mspace{-4mu} {-}1 -2L-(n{-}4)+c\\
 &=\begin{cases}
    (n{-}2)(L{-}1)+\textstyle 2\left\lceil\frac{L}{p}\right\rceil + 1 & c=0\\
    \bigg.(n{-}2)\phantom{(}L\phantom{-1)}+\textstyle 2\left\lceil\frac{L+1}{p}\right\rceil + \left\lbrace\begin{smallmatrix}-1\\\phantom{-}0\end{smallmatrix}\right. & c=n{-}4, n{-}3, \text{ as }p\nmid L\\
    (n{-}2)(L{+}1)+\textstyle 2\left\lceil\frac{L+2}{p}\right\rceil -2 & c=2n{-}7 \text{ if }p\nmid L{+}1\\
    \bigg.(n{-}2)(L{+}1)+\textstyle 2\left\lceil\frac{L+2}{p}\right\rceil -4 & c=2n{-}7 \text{ if }p\mid L{+}1
   \end{cases}
   \allowdisplaybreaks
   \shortintertext{$\bullet$ on $\widetilde{\mathcal{N}}^{L}_-$ for period $\frac{L}{(1-\veps)}{\cdot}\frac{\pi}{2}, L\in2\mathbb{N}$ and $c\in\{0,1\}$}
 \mu\Big(L\gamma_c^-\Big)&=\textstyle\left\lfloor\mspace{-3mu}\frac{L(1+\veps)}{2(1-\veps)}\mspace{-3mu}\right\rfloor\mspace{-5mu}+\mspace{-5mu}\left\lceil\mspace{-3mu}\frac{L(1+\veps)}{2(1-\veps)}\mspace{-3mu}\right\rceil+\left\lfloor\mspace{-3mu}\frac{L(1-\veps)}{2(1-\veps)}\mspace{-3mu}\right\rfloor\mspace{-5mu}+\mspace{-5mu}\left\lceil\mspace{-3mu}\frac{L(1-\veps)}{2(1-\veps)}\mspace{-3mu}\right\rceil+\sum_{k=2}^{n-1}\textstyle \left(\left\lfloor\mspace{-3mu}\frac{L}{2(1-\veps)}\mspace{-3mu}\right\rfloor\mspace{-5mu}+\mspace{-5mu}\left\lceil\mspace{-3mu}\frac{L}{2(1-\veps)}\mspace{-3mu}\right\rceil\right)\\
 &\phantom{\,+\,}\textstyle+\left\lfloor\mspace{-3mu}\frac{L}{p(1-\veps)}\mspace{-3mu}\right\rfloor\mspace{-5mu}+\mspace{-5mu}\left\lceil\mspace{-3mu}\frac{L}{p(1-\veps)}\mspace{-3mu}\right\rceil -\left(\left\lfloor \mspace{-3mu}\frac{L}{(1-\veps)}\mspace{-3mu}\right\rfloor\mspace{-5mu}+\mspace{-5mu}\left\lceil \mspace{-3mu}\frac{L}{(1-\veps)}\mspace{-3mu}\right\rceil\right) -\frac{1}{2}\dim\widetilde{\mathcal{N}}^L_-+\frac{1}{2}+c\\
 &= \textstyle \left(\frac{L}{2}{+}\frac{L}{2}{+}1\right)+2\frac{L}{2}+(n{-}2)\left(\frac{L}{2}{+}\frac{L}{2}{+}1\right)+\textstyle 2\left\lceil\mspace{-3mu}\frac{L}{p(1-\veps)}\mspace{-3mu}\right\rceil \mspace{-4mu}{-}1-(L{+}L{+}1)+c\\
 &=\begin{cases}(n{-}2)(L{+}1)+2\left\lceil\frac{L+2}{p}\right\rceil + \left\lbrace\begin{smallmatrix}-1\\\phantom{\,-}0\end{smallmatrix}\right. &\text{if } p\nmid L{+}1\\
 \bigg.(n{-}2)(L{+}1)+2\left\lceil\frac{L+2}{p}\right\rceil + \left\lbrace\begin{smallmatrix}-3\\-2\end{smallmatrix}\right. &\text{if } p\mid L{+}1\end{cases}.
\end{align*}
From the above calculations, we can directly read off the number $b_k(\xi_p')$ of closed Reeb orbits of $\alpha_p'$ having Morse-Bott index $k$ as:
\begin{proposition}\label{perturbedbounds} $b_k(\xi_p') = 0$, unless
 \begin{align*}
 &b_k(\xi_p') = \begin{cases}	1 & \text{if } k \in f_p(l)+\{-2, +1\}\\
				2 & \text{if } k \in f_p(l)+\{-1,0\}
                \end{cases} && \text{for } p\nmid l{-}1 \text{ and } 2\mid l\\
 &b_k(\xi_p') = 1 \qquad \text{if } k \in f_p(l)+ \{-1,0\} &&\text{for } p\nmid l{-}1 \text{ and } 2\nmid l\\
 &b_k(\xi_p') = 1 \qquad \text{if } k \in f_p(l)+ \{-4,-3,-2,-1,0,+1\} &&\text{for } p\mid l{-}1 \text{ and } 2\mid l.
\end{align*}
\end{proposition}
Note that $b_k(\xi'_p)=0$ around $k=f_p(l)$ if $2p\mid l{-}1$. In particular holds
\begin{equation}\label{indextable}
 b_k(\xi'_p)=0	\quad\text{ if }	\begin{cases}	\quad\, k\in f_p(l{-}1)+\{+2,+3,+4,+5\}\\
				\text{or } k\in f_p(l)+\{-3,-2,-1,0\}\\
				\text{or } k\in f_p(l{+1})+\{-6,-5,-4,-3\},
				\end{cases}
\end{equation}
due to the following estimate (using $n{-}2\geq 3$)
 \begin{align*}
  f_p(l{+}1)\,{-}2=l(n{-}2)\,{+}2\big\lceil (l{+}1)/p\big\rceil\,{-}2 &\geq (l{-}1)(n{-}2)\,{+}3\,{+}2\big\lceil l/p\big\rceil\,{-}2 &&\mspace{-20mu}= f_p(l)\,{+}1\\
  &\geq(l{-}2)(n{-}2)\,{+}6\,{+}2\big\lceil(l{-}1)/p\big\rceil &&\mspace{-20mu}= f_p(l{-1})\,{+}6.
 \end{align*}
\begin{proof}[\textnormal{\textbf{Proof of Theorem \ref{differentcontacttheo}}}]~\\
We have to show that if $(j,p)\neq (i,q)$ for $i,j,p,q\in\mathbb{N}$ and $p,q\equiv \pm 1\bmod 8$, then $j{\cdot} \Sigma_p$ and $i{\cdot}\Sigma_q$ are different contact structures on $S^{4m+1}$ for $m\geq 2$. From the Invariance Theorem (Thm.\ \ref{theoinvsur}), we obtain that
 \begin{align*}
   SH_k\big(j{\cdot} V_p,j{\cdot} \Sigma_p\big)= {\textstyle \bigoplus_j} SH_k(V_p,\Sigma_p).\tag{$\ast$}
 \end{align*}
 Moreover, $\xi_p$ is asymptotically finitely generated in every degree $k$ with bound $b_k(\xi_p')$ for a contact form $\alpha_p'$ not having a periodic Reeb flow. With Prop.\ \ref{propcofinalHamilt}, it follows that $j{\cdot}\xi_p$ is also a.f.g.\ in every degree $k$ with bound $b_k\big(j{\cdot}\xi_p\big)=j\cdot b_k(\xi_p')$ if $k$ is not of the form $(2l{-}1)(n{-}1)+\left\lbrace\begin{smallmatrix}1\phantom{(n-1)}\\2(n-1)\end{smallmatrix}\right.,\, l\in\mathbb{N}$. For simplicity, we also exclude $k=n-1$, which allows us to rewrite the assumptions on $k$ as
 \[k\neq (2l{-}1)(n{-}1)+\left\lbrace\begin{smallmatrix}0\\1\end{smallmatrix}\right.,\qquad l\in\mathbb{N}.\tag{$\ast\ast$}\]
 For such $k$, Prop.\ \ref{PropAsympfiniteGene} gives for any Liouville filling $V$ of $(S^{2n-1}, j{\cdot}\xi_p)$ the estimate
 \[\rk SH_k\big(V,j{\cdot}\Sigma_p\big)\leq j\cdot b_k(\xi_p').\tag{$\ast{\ast}\ast$}\]
 We assume without loss of generality $p<q$ and distinguish the following cases:\\
 \underline{If $q>2p{+}1$}, then (since $q$ is odd) $q\geq2p{+}3$. Here, we look at the index
 \begin{align*}
  k=(2p{+}1)(n{-}2)+2&=(2p{+}1)(n{-}2)+2\left\lceil\frac{2p{+}2}{q}\right\rceil\phantom{-4} = f_q(2p{+}2)\\
  &= (2p{+}1)(n{-}2)+2\left\lceil\frac{2p{+}2}{p}\right\rceil{-}4=f_p(2p{+}2)-4.
 \end{align*}
 On one hand, we know by the additivity of $SH$ (see ($\ast$)) and Thm.\ \ref{theoSHofSigmap} that
 \[\rk SH_k\big(i{\cdot} V_q,i{\cdot} \Sigma_q\big)=i\cdot \rk SH_k(V_q,\Sigma_q) = i.\]
 On the other hand, we claim that $k$ is not of the form $(\ast\ast)$. Otherwise, as $k$ is odd (as $n$ is odd), there exists $l\in\mathbb{N}$ such that
 \begin{align*}
  &&(2p{+}1)(n{-}2)+2=k&=(2l{-}1)(n{-}1)+1 \\
  &\Leftrightarrow& 2p(n{-}2)+n&=(2l{-}2)(n{-}1)+n &&\Leftrightarrow& \frac{p(n{-}2)}{n{-}1}&=l{-}1.
 \end{align*}
But this is impossible as $p$ and $n$ are odd. Therefore $(\ast{\ast}\ast)$ and (\ref{indextable}) imply for any filling $V$ of $j\cdot \Sigma_p'$ that
 \[\rk SH_k\big(V,j{\cdot}\Sigma_p\big)\leq j\cdot b_k(\xi_p') =j\cdot 0=0\neq i =\rk SH_k\big(i{\cdot} V_q,i{\cdot} \Sigma_q\big).\]
 This implies that $ j\cdot \xi_p$ and $ i\cdot \xi_q$ cannot be contactomorphic.\\
 \underline{If $q<2p{+}1$}, then $p<q\leq 2p{-}1$ and we look at the indices
 \begin{align*}
  && k=2p(n{-}2)+4&=2p(n{-}2)+2\left\lceil\frac{2p{+}1}{q}\right\rceil \phantom{-2}= f_q(2p{+}1)\\
  &&&=2p(n{-}2)+2\left\lceil\frac{2p{+}1}{p}\right\rceil {-}2=f_p(2p{+}1)-2\\
  &\text{and}& k=2p(n{-}2)+3&=f_q(2p{+}1)-1=f_p(2p{+}1)-3
 \end{align*}
In both cases, we know by the additivity of $SH$ and Thm.\ \ref{theoSHofSigmap} that
 \[\rk SH_k\big(i{\cdot} V_q,i{\cdot} \Sigma_q\Big)=i\cdot \rk SH_k(V_q,\Sigma_q) = i.\]
 Now, the first value for $k$ is even, while the second one is odd. So if $k$ were of the form $(\ast\ast)$, then there exists some $l\in\mathbb{N}$ such that
 \begin{align*}
  &&2p(n{-}2)+4&=(2l{-}1)(n{-}1) \\
  &\Leftrightarrow& 2p(n{-}2)+4 &= 2l(n{-}1)-(n{-}1) &&\Leftrightarrow& p(n{-}2)+\frac{n{+}3}{2}&=l(n{-}1)\\
  &\text{or}&2p(n{-}2)+3&=(2l{-}1)(n{-}1)+1 \\
  &\Leftrightarrow& 2p(n{-}2)+2&=2l(n{-}1)-(n{-}1) &&\Leftrightarrow& p(n{-}2)+\frac{n{+}1}{2} &=l(n{-}1).
 \end{align*}
 The first case is impossible for $n\equiv 1\bmod 4$, while the second is impossible for $n\equiv 3\bmod 4$. So for $k$ appropriately chosen, we find again $\rk SH_k\left(V,j\cdot\Sigma_p\right)=0$ for any filling $V$ of $j\cdot \Sigma_p'$ by (\ref{indextable}) and ($\ast{\ast}\ast$). Hence $j\cdot \xi_p$ and $i\cdot \xi_q$ cannot be contactomorphic.\\
 \underline{If $q=2p{+}1$}, we look at the index
 \begin{align*}
  k=(4p{+}1)(n{-}2)+4&=(2q{-}1)(n{-}2)+2\left\lceil\frac{2q}{q}\right\rceil\phantom{+-6}=f_q(2q)\\
  &=(4p{+}1)(n{-}2)+2\left\lceil\frac{4p{+}2}{p}\right\rceil{-}6=f_p(4p{+}2)-6.
 \end{align*}
 We find again that $\rk SH_k\big(i{\cdot} V_q, i{\cdot} \Sigma_q\big) = i.$ Moreover, $k$ is odd, so if $k$ were of the form $(\ast\ast)$ for some $l\in\mathbb{N}$, then we would have
 \begin{align*}
  && (4p{+}1)(n{-}2)+4 &= (2l{-}1)(n{-}1)+1 \\
  &\Leftrightarrow& 4p(n{-}2)+n+2 &=(2l{-}2)(n{-}1)+n &&\Leftrightarrow& \frac{2p(n{-}2)+1}{n{-}1} &= l{-}1,
 \end{align*}
 which is impossible, as the numerator is odd, while the denominator is even. Thus, for any filling $V$ of $j\cdot \Sigma_p'$ follows $\rk SH_k\left(V,j\cdot\Sigma_p\right)=0$ 
 and hence that $j\cdot \xi_p$ and $i\cdot \xi_q$ are not contactomorphic.
\end{proof}
\begin{proof}[\textnormal{\textbf{Proof of Corollary \ref{CorMeanEul}}}]~\\
On $S^{4m+1}$, we have to find $k$ different contact structures with the same mean Euler characteristic $\chi_m$. In \cite{Esp}, 8.3, Espina calculated $\chi_m$ for $(\Sigma_p,\xi_p)$ as 
 \[\chi_m(\Sigma_p)=\frac{1}{2}\cdot\frac{(n{-}1)p+1}{(n{-}2)p+2}.\tag{$\ast$}\]
 Alternatively, calculations similar to whose in \cite{KwKo}, 5.8, can be used to obtain the same result. Moreover, in \cite{Esp}, Cor.\ 5.7, or \cite{BourOan2}, one finds the following formula for the mean Euler characteristic of the boundary connected sum of two Liouville domains $V_1,V_2$ with $\dim V_i=2n$ :
 \[\chi_m(V_1\# V_2)=\chi_m(V_1)+\chi_m(V_2)+(-1)^n{\textstyle\frac{1}{2}}.\tag{$\ast\ast$}\]
 Applying $(\ast)$ and $(\ast\ast)$ to the $l$-fold connected sum $l{\cdot}\Sigma_p$ of $\Sigma_p$ with itself, we obtain:
 \[\chi_m(l{\cdot}\Sigma_p)=\frac{l}{2}\left(\frac{(n{-}1)p+1}{(n{-}2)p+2}+(-1)^n\right)-\frac{(-1)^n}{2}.\]
 Note that $\frac{(n{-}1)p+1}{(n{-}2)p+2}>1$, as $p>1$, so that the expression in brackets is always positive. Now for any $k$ positive integers $p_1, ..., p_k$ with $p_i\equiv \pm 1\pmod 8$ and $p_i\neq p_j$ for $i\neq j$, it is not difficult to find integers $l_1,..., l_k$ such that $\chi_m(l_1\cdot\Sigma_{p_1})=...=\chi_m(l_k\cdot\Sigma_{p_k})$. We just have to require for all $2\leq i\leq k$ that
 \[l_i=l_1\cdot\frac{\frac{(n{-}1)p_1+1}{(n{-}2)p_1+2}+(-1)^n}{\frac{(n{-}1)p_i+1}{(n{-}2)p_i+2}+(-1)^n}.\]
 For that, it suffices to choose $l_1$ such that the product on the right side is an integer for all $2\leq i\leq k$. Then all $l_i\cdot\Sigma_{p_i}$ have the same mean Euler characteristic, but nevertheless they are not contactomorphic by Thm.\ \ref{differentcontacttheo}.
\end{proof}

\appendix
\section{Perturbing periodic Reeb flows}\label{appendixA}
The purpose of this appendix is to show that most contact manifolds which admit a contact form with periodic Reeb flow are asymptotically finitely generated in every degree, in particular via contact forms with a non-periodic Reeb flow. We need such contact forms for the estimate of the Symplectic Homology of connected sums with these manifolds (see Prop.\ \ref{propcofinalHamilt}).
\subsection{Setup}
Let $(\Sigma,\alpha)$ be a compact $(2n{-}1)$-dimensional contact manifold with contact structure $\xi=\ker \alpha$ and Reeb vector field $R$ such that the flow $\varphi^t$ of $R$ is periodic, i.e. there exists a time $T_0>0$ with $\varphi^{T_0}=Id$. Hence, any length $T\leq T_0$ of a closed Reeb orbit must satisfy $T_0/T \in\mathbb{N}$. Moreover, the length of a closed Reeb orbit is bounded from below, as $R(p)\neq 0$ for all $p\in\Sigma$. Thus there are only finitely many $0<T_1<...<T_N<T_0$ such that $T_j$ is the length of a closed Reeb orbit.\pagebreak[1]\\
The set of diffeomorphisms $G:=\{\varphi^{T_j}\,|\,0\leq j\leq N\}$ clearly carries the structure of a finite abelian group which acts on $\Sigma$. Let $Fix(\varphi^{T_j})=\{p{\in}\Sigma\,|\,\varphi^{T_j}(p){=}p\}$ denote the fixed point set of $\varphi^{T_j}$. Using the exponential map of a $G$-invariant Riemannian metric $g$, we can find an $\varphi^{T_j}$-invariant neighborhood $U$ around any $p\in Fix(\varphi^{T_j})$ such that the action of $\varphi^{T_j}$ on $U$ is conjugated to the action of $D_p\varphi^{T_j}$ on $T_p\Sigma$ restricted to a sufficiently small neighborhood $B$ of the origin. We find that $U\cap Fix(\varphi^{T_j})$ is conjugated to $B\cap Fix(D_p\varphi^{T_j})=B\cap \ker (D_p\varphi^{T_j}{-}\mathbbm{1})$. This shows that $Fix(\varphi^{T_j})$ is a submanifold of $\Sigma$ and that $T_pFix(\varphi^{T_j})=\ker (D_p\varphi^{T_j}{-}\mathbbm{1})$.\pagebreak[1]\\
As $Fix(\varphi^{T_j})=\mathcal{N}^{T_j}=\mathcal{N}^{T_j+k\cdot T_0}, k\in\mathbb{N}$, are exactly the sets formed by the periodic orbits of $R$, we find that every contact manifold $(\Sigma,\alpha)$ with periodic Reeb flow satisfies the Morse-Bott assumption (\ref{CondMB}). Moreover, $\xi$ is an a.f.g.\ contact structure in almost all degrees. Indeed, if the mean index $\Delta(\Sigma)$ of the principal orbit $\Sigma=Fix(\varphi^{T_0})=\mathcal{N}^{T_0}$ is non-zero, then we find by the iterations formula (\ref{eqIteration}) that the Conley-Zehnder index of the orbits $\mathcal{N}^{T_j+k\cdot T_0}$ growths approximately linear in $k$ with slope $\Delta(\Sigma)$. Hence there are only finitely many closed Reeb orbits having any given Morse-Bott index $k$. If however $\Delta(\Sigma)=0$, then the iterations formula implies that the Morse-Bott index of any closed Reeb orbit stays in the interval $(-2n{+}1,2n)$, so no closed Reeb orbit has Morse-Bott index $k$ outside $(-2n{+}1,2n)$. Hence, $\xi=\ker \alpha$ is a.f.g.\ in any degree $k$ if $\Delta(\Sigma)\neq0$ and in degrees $k\in\mathbb{Z}\setminus[-2n{+}2,2n{-}1]$ if $\Delta(\Sigma)=0$. As sequences $(f_l)$ and $(\mathfrak{a}_l)$ to show that $\xi$ is a.f.g., one can choose $f_l\equiv 0$ and any increasing sequence $(\mathfrak{a}_l)$  with $\mathfrak{a}_l\rightarrow \infty$.

\subsection{An explicit perturbation}
The following construction was first described in \cite{Bour}, Sec.\ 2 (see also \cite{KwKo}, Lem.\ 5.17). We present it here to construct a sequence $\alpha_l:=e^{f_l}\cdot\alpha,\; f_l:\Sigma\rightarrow\mathbb{R}$ of perturbations of $\alpha$ where the flow of the Reeb vector fields $R_l$ is not periodic and to add some details missing in the original argument.\\
In order to perturb $\alpha$, we first construct positive Morse functions $\bar{f}_1,..., \bar{f}_N,\bar{f}_0$ on the orbit spaces $Q_j:=Fix(\varphi^{T_j})\big/S^1$. Note that $Q_j$ is in general a symplectic orbifold as the Reeb flow does not act freely on $Fix(\varphi^{T_j})$ if there exists an $i<j$ such that $T_i|T_j$. We construct the $\bar{f}_j$ by the following inductive procedure:
\begin{enumerate}
 \item $Q_1$ is a smooth manifold. Thus pick any positive Morse function $\bar{f}_1$ on $Q_1$.
 \item For $Q_j$, $j>1$, the singular set of $Q_j$ is exactly $S_j=\bigcup_{T_i|T_j} Q_i$. At first extend for every $i<j$ with $T_i|T_j$ the functions $\bar{f}_i$ to a function $\tilde{f}_j$ on a small tubular neighborhood of $S_j$ via the quadratic function $v\mapsto ||v||^2_g$ on the normal bundles to $Q_i$ in $Q_j$ (here, $g$ is a fixed invariant metric on $\Sigma$). Then extend $\tilde{f}_j$ to a positive Morse function $\bar{f}_j$ on $Q_j$. If $Q_j$ has no singular set, pick any positive Morse function.
 \item Repeat this procedure for all $j$, in particular up to $Q_0=Fix(\varphi^{T_0})\big/S^1=\Sigma\big/S^1$.
\end{enumerate}
Now lift $\bar{f}_0$ to $\Sigma$ to obtain a function $f$ which is invariant under the Reeb flow $\varphi^t$. As $R\in\ker df$ for any $\varphi^t$-invariant function $f$ and as $\xi=\ker \alpha$ and $\alpha\wedge(d\alpha)^{n-1}$ is non-degenerate, we can associate to any $\varphi^t$-invariant function $f$ a unique Hamiltonian vector field $X_f$ by the two conditions
\[X_f(p)\in\xi(p)\quad \forall p\in\Sigma\qquad\text{ and }\qquad d\alpha(\cdot, X_f)=df.\]
Set $\alpha_\varepsilon:=(1{+}\varepsilon f){\cdot} \alpha$ as the perturbed contact form. Its Reeb vector field $R_\varepsilon$ is
\begin{align*}
R_\varepsilon&:=\frac{1}{1{+}\varepsilon f}\cdot R-\frac{\varepsilon}{(1{+}\varepsilon f)^2}\cdot X_{f}, \qquad\text{ since}\\
 \alpha_\varepsilon(R_\varepsilon) &= (1{+}\varepsilon f)\alpha\Big(\frac{1}{1{+}\varepsilon f}R-\frac{\varepsilon}{(1{+}\varepsilon f)^2}X_{f}\Big)=\frac{1{+}\varepsilon f}{1{+}\varepsilon f}\alpha(R)=1\allowdisplaybreaks\\
 d\alpha_\varepsilon (R_\varepsilon,\cdot) &=\big((1{+}\varepsilon f)d\alpha+\varepsilon df{\wedge}\alpha\big)\Big({\textstyle\frac{1}{1+\varepsilon f}}R-{\textstyle \frac{\varepsilon}{(1+\varepsilon f)^2}} X_{f} , \cdot\Big)\\
 &= d\alpha(R,\cdot)-{\textstyle \frac{\varepsilon}{1+\varepsilon f}} d\alpha(X_{f},\cdot)+{\textstyle \frac{\varepsilon}{1+\varepsilon f}}(df{\wedge}\alpha)(R,\cdot)-{\textstyle\frac{\varepsilon^2}{(1+\varepsilon f)^2}}(df{\wedge}\alpha)(X_{f},\cdot)\\
 &= 0 + {\textstyle \frac{\varepsilon}{1+\varepsilon f}}df-{\textstyle\frac{\varepsilon}{1+\varepsilon f}}df-0 =0,
\end{align*}
as $R\in\ker df$ and $X_{f}\in\xi=\ker \alpha$ and $df(X_{f})=d\alpha(X_{f},X_{f})=0$.\\ Finally, we calculate for the Lie bracket of $\frac{1}{1+\varepsilon f}R$ and $\frac{-\varepsilon}{(1+\varepsilon f)^2}X_{f}$ that
\begin{align*}
 \alpha\Big(\big[{\textstyle\frac{1}{1+\varepsilon f}}R, {\textstyle\frac{-\varepsilon}{(1+\varepsilon f)^2}} X_{f}\big]\Big) &=\mathcal{L}_{\frac{R}{1+\varepsilon f}} \underbrace{\alpha\big({\textstyle\frac{-\varepsilon}{(1+\varepsilon f)^2}} X_{f}\big)}_{=0}-\iota_{\frac{-\varepsilon X_{f}}{(1+\varepsilon f)^2}}\mathcal{L}_{\frac{R}{1+\varepsilon f}} \alpha\\
 &=-\iota_{\frac{-\varepsilon X_{f}}{(1+\varepsilon f)^2}}\Big(\underbrace{\iota_{\frac{R}{1+\varepsilon f}} d\alpha}_{=0}+d\Big(\underbrace{\alpha\big({\textstyle\frac{1}{1+\varepsilon f}}R\big)}_{=\frac{1}{1+\varepsilon f}}\Big)\Big)\\
 &={\textstyle \frac{-\varepsilon}{(1+\varepsilon f)^2}} df\Big({\textstyle \frac{-\varepsilon}{(1+\varepsilon f)^2}} X_{f}\Big)=0\allowdisplaybreaks\\
 d\alpha\Big(\big[{\textstyle\frac{1}{1+\varepsilon f}}R, {\textstyle\frac{-\varepsilon}{(1+\varepsilon f)^2}} X_{f}\big],\cdot\Big) &=\mathcal{L}_{\frac{R}{1+\varepsilon f}} d\alpha\big({\textstyle\frac{-\varepsilon}{(1+\varepsilon f)^2}} X_{f},\cdot\big)-\iota_{\frac{-\varepsilon X_{f}}{(1+\varepsilon f)^2}}\mathcal{L}_{\frac{R}{1+\varepsilon f}} d\alpha\\
 &=\mathcal{L}_{\frac{R}{1+\varepsilon f}} {\textstyle\frac{\varepsilon}{(1+\varepsilon f)^2}} df-\iota_{\frac{-\varepsilon X_{f}}{(1+\varepsilon f)^2}}\Big(\iota_{\frac{R}{1+\varepsilon f}}\underbrace{dd\alpha}_{=0}+d\big(\underbrace{\mspace{-3mu}d\alpha({\textstyle\frac{1}{1+\varepsilon}}R,\cdot)\mspace{-3mu}}_{=0}\big)\mspace{-5mu}\Big)\\
 &=\iota_{\frac{R}{1+\varepsilon f}}d\big({\textstyle\frac{\varepsilon}{(1+\varepsilon f)^2}}df\big)+d\Big({\textstyle\frac{\varepsilon}{(1+\varepsilon f)^2}} \underbrace{df\big({\textstyle\frac{1}{1+\varepsilon f}}R\big)}_{=0}\Big)\\
 &=\iota_{\frac{R}{1+\varepsilon f}}\Big({\textstyle\frac{\varepsilon}{(1+\varepsilon f)^2}}\underbrace{dd f}_{=0}-{\textstyle\frac{-2\varepsilon^2}{(1+\varepsilon f)^3}}\underbrace{df\wedge df}_{=0}\Big)=0
 \end{align*}
 \[\Rightarrow\;\;\big(\alpha\wedge(d\alpha)^{n-1}\big)\Big(\big[{\textstyle\frac{1}{1+\varepsilon f}}R, {\textstyle\frac{-\varepsilon}{(1+\varepsilon f)^2}} X_{f}\big],\cdot\Big) =0.\]
As $\alpha{\wedge}(d\alpha)^{n-1}$ is a volume form, we conclude that $\big[{\textstyle\frac{1}{1+\varepsilon f}}R, {\textstyle\frac{-\varepsilon}{(1+\varepsilon f)^2}} X_{f}\big]=0$ and hence that the flows of $\frac{1}{1+\varepsilon f}R$ and $\frac{-\varepsilon}{(1+\varepsilon f)^2}X_f$ commute. This implies that the flow of $R_\varepsilon$ is the composition of these two flows.\\
Compositions of closed orbits of $\frac{1}{1+\varepsilon f}R$ and $\frac{-\varepsilon}{(1+\varepsilon f)^2}X_f$ are hence closed orbits of $R_\varepsilon$, if both orbits have the same starting point and the same period. In particular, we have closed orbits of $R_\varepsilon$ through any critical point $p$ of $f$. These correspond exactly to critical points of $\bar{f}_0$. As $\bar{f}_0$ is Morse, we find that the Hessian of $f$ in the normal direction to the closed $R$-orbit through $p$ is non-degenerate and hence that these closed orbits of $R_\varepsilon$ are non-degenerate.\\
More general, closed orbits of $R_\varepsilon$ are in one-to-one correspondence with $\frac{-\varepsilon}{(1+\varepsilon f)^2}X_f$-chords of $\frac{1}{1+\varepsilon f}R$-orbits, where the length of the chord equals the flowtime from one endpoint to the other along the $\frac{1}{1+\varepsilon f}R$-trajectory.
By Prop.\ \ref{smallHamilt} below, we can choose for any given $\mathfrak{a}>0$ the constant $\varepsilon$ so small that the only $\frac{-\varepsilon}{(1+\varepsilon f)^2}X_f$-chords of length $\mathfrak{a}$ are the constant ones. This implies that for any given $\mathfrak{a}$ and $\varepsilon$ sufficiently small the only closed orbits of $R_\varepsilon$ with period less than $\mathfrak{a}$ are those over critical points of $\bar{f}_0$. As $\bar{f}_0$ is Morse on a compact set, we find that their number is finite and independent of $\varepsilon$. In \cite{Bour}, Lem.\ 2.4, the Conley-Zehnder index of these orbits through $p\in Fix(\varphi^{k\cdot T_j})$ was calculated as 
\[\quad\mu_{CZ}\Big(Fix(\varphi^{k\cdot T_j})\Big)-{\textstyle \frac{1}{2}}\dim\Big(Fix(\varphi^{k\cdot T_j})\big/ S^1\Big)+\mu_{Morse} \big(\bar{f}_{k\cdot T_j}([p])\big).\]
Here, $\bar{f}_{k\cdot T_j}=\bar{f}_i$ for $k T_j\equiv T_i \bmod T_0$. If the mean index $\Delta(\Sigma)$ of the principal orbit $Fix(\varphi^{T_0})=\Sigma$ is zero, then this implies that all Conley-Zehnder indices of closed orbits of $R_\varepsilon$ shorter then $\mathfrak{a}$ stay in the interval $(-2n{+}1,2n)$. If $\Delta(\Sigma)$ is non-zero, then only finitely many of these orbits have degree $k$ for any given integer $k$.\\
Hence, we can choose to any increasing sequences $\mathfrak{a}_l\rightarrow\infty$ a decreasing sequence $\varepsilon_l\rightarrow 0$ such that $(\mathfrak{a}_l)$ and $\big(f_l{:=}\log (1{+}\varepsilon_lf)\big)$ show that $\xi=\ker \alpha$ is asymptotically finitely generated in degree $k$, for any $k$ if $\Delta(\Sigma)\neq 0$ or for $k\in\mathbb{Z}\setminus[-2n{+}2,2n{-}1]$ if $\Delta(\Sigma)=0$. Moreover, the Reeb flow for any $\alpha_l:=e^{f_l}\cdot\alpha$ is not periodic.\bigskip\\
To finish the argument, it only remains to show the following contact version of the well-known fact that for any $C^2$-small Hamiltonian the only 1-periodic orbits of $X_H$ are the constant orbits at critical points of $H$.
\begin{proposition}\label{smallHamilt}
 Let $(\Sigma,\alpha)$ be a compact contact manifold with Reeb vector field $R$ and let $H:\Sigma\rightarrow\mathbb{R}$ be a smooth function which is invariant under the flow of $R$. Then there exists a constant $T>0$ such that each $X_H$-chord of a Reeb orbit of length less than $T$ is a constant orbit at a critical point of $H$. Here, $X_H$ is the Hamiltonian vector field of $H$ as defined above.
\end{proposition}
\begin{proof}
 We argue by contradiction: Assume that there exists a sequence of $X_H$-chords $(\gamma_l)$  of Reeb orbits with lengths $T_l>0$ and $\lim_{l\rightarrow\infty} T_l =0$. As $\Sigma$ is compact, we may assume by the Arzela-Ascoli Theorem that $(\gamma_l$) converges uniformly to a $X_H$-chord $\gamma$ of length $T=0$. Hence, $\gamma$ is a constant orbit at a critical point $p$ of $H$, as $d\alpha(\cdot, X_H)=dH$.\\
 Around $p$ we can find by the contact version of Darboux's Theorem (see \cite{Geiges}, Thm. 2.5.1) a neighborhood $U$ and coordinates $x_1,...,x_{n-1},y_1,...,y_{n-1},z$ such that $p=(0,...,0)$ and
 \[\alpha|_U = dz+\sum_{j=1}^{n-1} x_jdy_j.\]
 Note that $\partial_z$ is the Reeb vector field in these coordinates. Without loss of generality we may assume that $U=(-\veps,\veps)^{2n-1}$ for some small $\veps>0$. The quotient $Q$ of $U$ under the flow of the Reeb vector field $\partial_z$ is then easily identified with $(-\veps,\veps)^{2n-2}$ with the coordinates $x_j,y_j, j=1,...,n{-}1$. On $Q$ we have the symplectic form $\bar{\omega}=\sum_j dx_j{\wedge} dy_j$. If $\pi:U\rightarrow Q$ denotes the quotient map, it is easy to see that $\pi^\ast \bar{\omega}=d\alpha$. Moreover, as $H$ is invariant under the flow of $\partial_z$, we find that it descends to a function $\bar{H}$ on $Q$ with $\bar{H}\circ \pi= H$. Consequently, we find that the Hamiltonian vector fields satisfy $d\pi(X_H)=X_{\bar{H}}$. Thus, any $X_H$-chord $\gamma$ of a Reeb orbit in $U$ yields a closed (!) $X_{\bar{H}}$-orbit $\bar{\gamma}=\pi\circ \gamma$ in $Q$.\pagebreak[1]\\
 Now recall that on $\mathbb{R}^{2n-2}$ with the symplectic form $\sum_j dx_j{\wedge} dy_j$ it is a well-known fact that there exists for a $C^2$-bounded autonomous Hamiltonian $\bar{H}$ a constant $T$ such that any closed $X_{\bar{H}}$-orbit of period less then $T$ is constant (see \cite{Laud}, Lem.\ 2.2). Thus almost all of the $\bar{\gamma}_l$ are constant. This implies $X_H(\gamma_l)=\dot{\gamma}_l\in\mathbb{R}\partial_z=\mathbb{R} R$ for almost all $l$ and hence that almost all $\gamma_l$ are constant, as $X_H(\gamma_l)\in \xi$ --- a contradiction to our assumption on $\gamma_l$.
\end{proof}
We finish this appendix by providing explicit subcritical isotropic spheres $S\subset\Sigma$ satisfying the Reeb chord assumption from Prop.\ \ref{propcofinalHamilt}. To start, choose ${p\in\Sigma\setminus crit(f)}$ and fix, as in the proof of Prop.\ \ref{smallHamilt}, a neighborhood $U$ of $p$ and coordinates $x_1,...,x_{n-1},\linebreak[1] y_1,...,y_{n-1}, z$, such that in these coordinates $p=(0,...,0)$ and 
\[\alpha|_U=dz+\sum_{j=1}^{n-1}x_jdy_j,\qquad X_f(p)=\partial_{x_1} \qquad\text{ and }\qquad R=\partial_z.\]
Then consider the subcritical isotropic sphere
\[S:=\Big\{(x,y,z)\,\Big|\,\sum_{j=1}^k y_j^2=\veps^2, y_j=0, j>k, z=x_j=0\Big\},\]
with $\veps>0$ so small, that $S\subset U$. As $X_f(p)=\partial_{x_1}$, we can assume for $\veps$ small enough, that $X_f|_S$ is transverse to the plane $E:=\big\{(x,y,z)\,\big|\,z{=}x{=}0\big\}\supset S$. For each real number $\mathfrak{a}_l$ we can choose $\varepsilon_l$ so small that $S$ has no $R_\varepsilon$-Reeb chords of length shorter than $\mathfrak{a}_l$. This follows from the facts that the flow of $R_\varepsilon$ is the composition of the flows of $\frac{1}{1+\varepsilon f}R$ and $\frac{-\varepsilon}{(1+\varepsilon f)^2}X_f$, that the flow of $\frac{1}{1+\varepsilon f}R$ is a reparametrization of the (periodic) flow of $R$ and that the flow $\frac{-\varepsilon}{(1+\varepsilon f)^2} X_f$ near $S$ is transverse to $S$ and arbitrarily slow for $\varepsilon$ sufficiently small.

\section{Isotropic submanifolds avoiding Reeb chords}\label{appendixB}
The following appendix was thankfully communicated to the author by Dingyu Yang. We include it here as we could not find detailed proofs of the Theorems \ref{THMNOCHORDS} and \ref{FINALTHM} below anywhere in the literature.\bigskip\\
Let $(M, \xi)$ be a $2n{+}1$-dimensional cooriented contact manifold and let $S\subset M$ be a closed isotropic submanifold (i.e.\ $TS\subset \xi|_S$) of dimension $k<n$ (subcritical). Let $\iota_0: S\to M$ denote the inclusion. Let $\lambda$ be a contact 1-form ($\xi=\ker\lambda$) with the associated Reeb vector field $R_\lambda$.\\
Recall the neighborhood theorem of an isotropic submanifold in a contact manifold. The normal bundle $N_S M$ is isomorphic to $\R R_\lambda\oplus (TS)^{d\lambda}/TS\oplus \xi|_S/(TS)^{d\lambda_0}$, where $(TS)^{d\lambda}:=\big\{v\in \xi|_{S}\,\big|\, d\lambda (v, \tilde v)=0 \text{ for all } \tilde v\in TS\big\}$ and $d\lambda$ induces a symplectic vector bundle structure on $(TS)^{d\lambda}/TS$. Pick an almost complex structure $J$ on $\xi$ and denote the Riemannian metric $g:=\lambda\otimes \lambda+d\lambda\circ (Id_\xi{\times} J)$. Then $(TS)^{d\lambda}/TS$ can be identified with $E:=(\R R_\lambda \oplus TS\oplus JTS)^\perp$ a subbundle of $TM|_{S}$. Denote the normal exponential map $\exp^\perp: E\to M$, by restricting the exponential map given by $g$. Let $U(E)$ be a (poly)disk bundle on which $\exp^\perp$ is a diffeomorphism onto the image. Then the neighborhood theorem says that an open neighborhood of $S$ in $M$ is contactomorphic to a neighborhood of the zero section in $\big((\R_z{\times} T^\ast S)\oplus U(E), dz\,{-}\lambda_0\,{+}(\exp^\perp)^\ast \lambda \big)$\footnote{\label{LAMBDAE}Here, instead of $(\exp^\perp)^\ast\lambda$, any other 1-form $\lambda_E$ on $E$ as a manifold such that $\lambda_E|_{(TE|_{S})}=0$ and $d\lambda_E|_{(\Lambda^2 TE|_{S})}=d\lambda|_{(\Lambda^2 TE|_{S})}\equiv d\lambda|_{(\Lambda^2 E|_{S})}\oplus 0$ ($S$ is identified with the zero section) will also work.}. Here $\R_z\times T^\ast S$ is the 1-jet bundle over $S$, $\lambda_0$ is the canonical 1-form on $T^\ast S$, and the Whitney sum bundle $(\R_z\times T^\ast S)\oplus U(E)$ over $S$ is regarded as a manifold in this local model. This version of the neighborhood theorem is not explicitly written, but can be directly deduced from Weinstein's neighborhood theorem between diffeomorphic isotropic submanifolds in two contact manifolds with isomorphic normal data (see e.g.\ \cite{Geiges}, Thm.\ 2.5.8, 6.2.2).\\
To make this appendix self-contained without the use of the h-principle as a black box, but sufficient for the main body of this article, we make the following assumption on $S$. 

\begin{defn} An isotropic submanifold $S$ in a contact manifold $(M, \xi)$ is said to have a \textbf{trivial neighborhood}, if $(TS)^{d\lambda}/TS$ (or equivalently $E$ as above) is a trivial bundle for some contact form $\lambda$ (and $J$ on $\xi$). This condition depends only on $\xi$ and $S$.
\end{defn} 

Due to the neighborhood theorem above (and footnote \ref{LAMBDAE}), if $S$ has a trivial neighborhood, then there exists an open neighborhood $W$ of $S$ in $(M, \lambda)$ and a contactomorphism
\begin{align}\label{MapPhi}
\Phi: (W, \lambda)\rightarrow \big((-\delta, \delta)_z\times T_\delta^\ast \big(S{\times} (-\delta, \delta)^{n-k}\big), dz-\lambda_{\text{st}}\big)
\end{align}
of $(W,\lambda)$ onto an open neighborhood of $S\times \{0\}^{n-k}$ in the 1-jet bundle of $S\times (-\delta, \delta)^{n-k}$ such that $\Phi^\ast (dz-\lambda_{\text{st}})=h \lambda$, where $\lambda_{\text{st}}$ is the canonical 1-form in the cotangent bundle and $h$ is a nowhere vanishing function\footnote{$h$ can be taken to be 1 by using the Poincar\'e Lemma, but the given form suffices.}. Here, the disc cotangent bundle $T^\ast_\delta(S{\times} (-\delta, \delta)^{n-k})$ is defined using a Riemannian metric $g_{S}\times g_0$ on $S\times (-\delta, \delta)^{n-k}$ with $g_0$ standard. 

\begin{remark} For the handle attachment along a subcritical isotropic submanifold $S$ in the main article, the assumption of having a trivial neighborhood will always be satisfied. Slightly more generally, if $(TS)^{d\lambda}/TS$ admits a Lagrangian subbundle $\mathcal{L}$, e.g. if $S$ can be embedded into a Legendrian submanifold of $(M, \xi)$, one can still argue as below with $\mathcal{L}$ replacing $S\times (-\delta, \delta)^{n-k}$ with minor modifications. One can also get Thm.\ \ref{FINALTHM} for finitely many contact forms for general $(TS)^{d\lambda}/TS$, by carefully using the argument below to achieve (CT) together with the h-principle. It is plausible that full Thm.\ \ref{FINALTHM} in general might be obtainable by more deeply combining h-principle or jet transversality with the argument below, but it is beyond the scope of this short appendix.
\end{remark}

We will assume that $S$ has a trivial neighborhood from now on, unless stated otherwise. Let $\varphi^t$ be the Reeb flow of $R_\lambda$. For some $T\in (0,\infty]$, we introduce a mild condition below for $S\subset M$, called $(\lambda, T)$-transversality condition. We aim to show for $S\subset M$ satisfying the $(\lambda, T)$-transversality condition, that for a generic perturbation $\iota$ of the inclusion $\iota_0: S\rightarrow M$ (in a suitable perturbation space) which remains an isotropic embedding from $S$ into $M$ and is isotropically isotopic to $\iota_0$, there is no Reeb chord starting and ending on $\iota(S)$ of length $\leq T$.

\begin{defn}\label{defnLT} The submanifold $S\subset M$ is said to satisfy the \textbf{$(\lambda, T)$-transversality condition} if for some open neighborhood $U(S)$ of $S$ in $M$, the map $\mathring{\textnormal{ev}}: U(S){\times} (0, T]\to M{\times} M, (x,t)\mapsto (x,\varphi^t(x))$ is transverse to the diagonal $\Delta_M$. We can assume without loss of generality that for $W$ in the domain of $\Phi$ in (\ref{MapPhi}) above, we have $W\subset U(S)$.
\end{defn}

The condition is only non-vacuous over $(x, x)\in M{\times} M$ with $\varphi^t(x)=x\in U(S)$ and $t\in (0,T]$. The upper bound $T$ here makes $(\lambda, T)$-transversality easier to achieve and it suffices for applications. Lem.\ \ref{THICKENING} and Cor.\ \ref{GENERICITY} below do not use $T<\infty$ and work equally for $T=\infty$. By openness of transversality, we only need to require for the $(\lambda, T)$-transversality  that $d{\mathring{\textnormal{ev}}}\big(TM|_{S}{\times} T((0,T])\big)$ and $T\Delta_M$ span $T(M{\times} M)|_{\Delta_S}$.

\begin{remark} The $(\lambda, T)$-transversality condition is further equivalent to the following: For each fixed $t\in (0,T]$, the graph $\text{gr}(d\varphi^t|_{\xi|_S}):=\{(v, d\varphi^t(v)\,|\, v\in \xi|_S\}$ and the diagonal $\Delta_\xi:=\{(v,v)\,|\, v\in \xi\}$ span $(\xi{\times} \xi)|_{\Delta_S}=\{(v,w)\,|\, v, w\in \xi_z, z\in S\}$. This is seen by freezing the $x$-variable, and note that $d\mathring{\textnormal{ev}}(\frac{\partial\;}{\partial t})$ is $(0,1)\in \mathbb{R} R_\lambda{\times} \mathbb{R} R_\lambda\subset T(M{\times} M)$ and transverse to the diagonal in $\mathbb{R} R_\lambda\times \mathbb{R} R_\lambda$.
\end{remark}

\begin{example}\label{FORFREE} The following examples satisfy the $(\lambda, T)$-transversality condition:
\begin{enumerate}
\item When $S$ avoids periodic Reeb orbits of periods $T'\leq T$. More generally, if $\lambda$ is such that all closed orbits of periods $T'\leq T$ form submanifolds of dimensions $< 2n{+}1{-}k$ and $\{\text{period } T'\}$ is a countable subset of $(0, T]$, then a perturbation by a generic normal vector field in the trivial neighborhood (see \ref{MapPhi}) will displace $S$ away from closed orbit submanifolds while the result is still an isotropic submanifold and isotropically isotopic to $S$\footnote{This argument is a baby version of the argument in Lem.\ \ref{THICKENING} and Cor.\ \ref{GENERICITY} below, and can be easily done in finite dimensions, as only one $S$-factor is involved.}.
\item If $\lambda$ is non-degenerate, then $(\lambda, \infty)$-transversality holds by the equivalent definition in the above remark, and there is no restriction on the submanifold $S\subset M$. Note that $(\lambda,T)$-transversality actually requires that $\lambda$ is (not just transversely) non-degenerate at the intersections of $S$ with $T'$-periodic Reeb orbits, $T'\leq T$.
\end{enumerate}
\end{example}

The relevant condition for the problem of avoiding Reeb chords is the following: 

\begin{defn}\label{CTTRANSV}
An embedding $\iota: S\to M$ is said to satisfy the \textbf{chord transversality condition $(CT)$} if $\mathring{\text{\textnormal{ev}}}_{\iota}: S\times S\times (0,T]\to M\times M, (x, y, t)\mapsto (\iota (x), \varphi^t (\iota (y)))$ is transverse to the diagonal $\Delta_M$ in $M\times M$. Note that $\varphi^0=Id$.
\end{defn}

\begin{corollary}\label{NOCHORD} Let $\iota: S\to M$ be an embedding of a closed $k$-manifold into a $(2n{+}1)$-dimensional contact manifold $(M,\xi)$ such that $\iota$ satisfies $(CT)$. If $k<n$, then $\iota(S)$ has no Reeb chord of length $\leq T$.
\end{corollary}

\begin{proof} A point $(x,y,t)\in(\mathring{\text{\textnormal{ev}}}_\iota)^{-1}(\Delta_M)$ in the definition $(CT)$ specifies the starting and ending points $(\iota(y),\iota(x))$ of a Reeb chord on $\iota(S)$ and its length $t$. By $(CT)$ and a dimension count, one finds that $(\mathring{\text{\textnormal{ev}}}_\iota)^{-1}(\Delta_M)$ is a manifold of dimension $2k{+}1+2n{+}1\linebreak[1]-2(2n{+}1)\big)=2(k-n)<0$, thus empty.\end{proof}

Let us recall the definition of a thickening:

\begin{defn} Let $g: U\to Y$ be a map. A \textbf{thickening} $\tilde g: \tilde U\to Y$ of $g$ consists of a fiber bundle projection $\tilde U\to U$ with a section $i_U$ and a map $\tilde g: \tilde U\to Y$ with $\tilde g\circ i_U=g$. $\tilde U$ will be a polydisk Banach trivial bundle with the fiber as a perturbation parameter space.
\end{defn}

In order to perturb $\iota_0$ through isotropic embeddings to an $\iota$ satisfying $(CT)$, we will proceed in two steps: given $\iota_0$ satisfying the $(\lambda, T)$-transversality condition, one finds a thickening $\mathcal{I}$ of $\iota_0$ satisfying a variant of the $(CT)$ condition and $\mathcal{I}(\cdot, c)$ is isotropically isotopic to $\iota_0$ (see Lem.\ \ref{THICKENING}). Then from $\mathcal{I}$ one gets a generic perturbation $\iota_c:=\mathcal{I}(\cdot, c)$ of $\iota_0$ such that $\iota_c$ satisfies the $(CT)$ condition (see Cor.\ \ref{GENERICITY}).\bigskip\\
Since $S$ has a trivial neighborhood, we have a local model given by $\Phi$ in (\ref{MapPhi}).

Embed $(-\delta, \delta)^{n-k}$ in $T^{n-k}:=(\R/\rho \mathbb{Z})^{n-k}$ for $\rho>2\delta$. Define the Floer-$\mathbf{\epsilon}$ space
$$\mathcal{C}:=\Big\{f\in \mathcal{C}^\infty_\epsilon(S\times T^{n-k}, \R)\;\Big|\; f|_{S\times (T^{n-k}\backslash (-\delta, \delta)^{n-k})}=0\Big\},$$ where $\epsilon=(\epsilon_n)_n$ with $\epsilon_0=\epsilon_1=1$ is a sequence of non-negative reals sufficiently fast decreasing to $0$ such that $\mathcal{C}$ is a Banach space with the norm $\|\cdot\|_{\mathcal{C}^\infty_\epsilon(S\times T^{n-k},\R)}$ (see \cite{FLOERSPACE}). Denote $\mathcal{C}_\delta:=\big\{f\in \mathcal{C}\,\big|\, \|f\|_{\mathcal{C}^\infty_\epsilon(S\times T^{n-k}, \R)}<\delta\big\}$, a Banach manifold. We will simply use $f$ to denote $f|_{S\times (-\delta, \delta)^{n-k}}$ and also write $f(a,b)$ for $a\in S$ and $b\in (-\delta, \delta)^{n-k}$. 

Similarly, define the Floer-$\epsilon$ space $\mathcal{C}^\infty_\epsilon(S,\R)$ with the norm $\|\cdot\|_{\mathcal{C}^\infty_\epsilon(S,\R)}$ which is a Banach space, and denote $\mathcal{C}^\infty_{\epsilon, \delta}(S,\R)$ with norm less than $\delta$ which is a Banach manifold. Let $\mathcal{N}:=\{s:=(s_1, \cdots, s_{n-k})\circ \text{diag}\;|\; (s_i)_i\in (\mathcal{C}^\infty_{\epsilon,\delta}(S,\R))^{n-k}\}\equiv (\mathcal{C}^\infty_{\epsilon,\delta}(S,\R))^{n-k}\circ \text{diag}$, where $\text{diag}: S\to S^{n-k}$ is the diagonal map. $\{\text{gr}_s: x\mapsto (x, s(x))\;|\; s\in\mathcal{N}\}$ is regarded as the space of smooth sections (with norm control) of the trivial bundle $S\times (-\delta, \delta)^{n-k}$.

For fixed $(f,s)\in \mathcal{C}_\delta\times \mathcal{N}$, we define an isotropic embedding $S\to M$ as follows: 
\begin{align*}
\iota_{f,s}: S\to M, x \mapsto \Phi^{-1}(f(x, s(x)), (x, s(x)), d_{(x,s(x))}f) &\equiv  \Phi^{-1}(f(\text{gr}_s(x)), \text{gr}_s(x), d_{\text{gr}_s(x)}f) \\
&\equiv  \Phi^{-1}(f(\text{gr}_s(x)), \text{gr}_{df}(\text{gr}_s(x))),
\end{align*}
where $\text{gr}_s(x):=(x, s(x))$ and we also write $df(\tilde x):=d_{\tilde x}f$ for $\tilde x\in S\times (-\delta,\delta)^{n-k}$.
It is clearly an embedding as $\iota_0$ is perturbed by a `normal' vector field and also isotropic via the following argument: For $\text{gr}_{df}\circ \text{gr}_s \equiv \big(x\mapsto (\text{gr}_s(x), d_{\text{gr}_s(x)}f)\big)$
 as a map into the cotangent bundle restricted to $\text{gr}_s(S)$, we have $(\text{gr}_{df}\circ \text{gr}_s)^\ast \lambda_{\text{st}}=\text{gr}_{df}\circ \text{gr}_s$,
 which is canceled by the differential of the $z$-coordinate, thus $\iota_{f,s}^\ast \lambda=\frac{1}{h}(\Phi\circ \iota_{f,s})^\ast(dz-\lambda_{\text{st}})=0$. By first isotoping $\iota_{f, (1-t) s}$, then $\iota_{(1-t) f,0}$ and smoothing out, clearly $\iota_{f,s}$ is isotopic to $\iota_0\equiv \iota_{0,0}$ via isotropic embeddings. So we have constructed a (injective) map $$\mathfrak{emb}: \mathcal{C}_\delta\times \mathcal{N}\to \text{Emb}^{\text{isot}}_{\iota_0}(S, M), (f,s)\mapsto \iota_{f,s}$$ into the space of isotropic embeddings $S\to M$ that are isotropically isotopic to $\iota_0$.\pagebreak[1]\\ Define $\text{EV}: S\times \text{Emb}^{\text{isot}}_{\iota_0}(S, M)\to M, (x, \iota )\mapsto \iota(x)$ and $$\mathcal{I}:=\text{EV}\circ (Id_S\times \mathfrak{emb}): S\times (\mathcal{C}_\delta\times\mathcal{N})\to M, (x,(f,s))\mapsto \iota_{f,s}(x).$$

\begin{lemma}\label{THICKENING} Let $S$ be a closed isotropic submanifold having a trivial neighborhood in a contact manifold $(M,\xi)$ such that $S\subset M$ satisfies the $(\lambda, T)$-transversality condition. Then the thickening $\mathcal{I}: S\times (\mathcal{C}_\delta\times \mathcal{N}) \to M$ of $\iota_0$ defined above has the following properties:
\begin{enumerate}
\item[(i)] $\mathcal{I}(\cdot, c): S\to M$ is an isotropic embedding and isotropically isotopic to $\iota_0$.
\item[(ii)] $\mathring{\text{\textnormal{ev}}}_{\mathcal{I}}: (S\times S\times (0,T])\times (\mathcal{C}_\delta\times \mathcal{N})\to M\times M$, $$((x, y, t), c)\mapsto \mathring{\text{\textnormal{ev}}}_{\mathcal{I}}((x,y,t),c):=\mathring{\text{\textnormal{ev}}}_{\mathcal{I}(\cdot, c)} (x,y, t)\equiv (\mathcal{I}(x,c), \varphi^t(\mathcal{I}(y,c))),$$ is transverse to $\Delta_M$ in $M\times M$, here $c=(f, s)\in \mathcal{C}_\delta\times \mathcal{N}$.
\end{enumerate}
\end{lemma}

Note that condition $(ii)$ says that $\mathring{\text{\textnormal{ev}}}_\mathcal{I}$ is transverse to $\Delta_M$ even with the \textit{same} perturbation parameter for both $S$ factors (essential for carrying out Cor.\ \ref{GENERICITY}).
\begin{proof}
We have already established property $(i)$.

To prove property $(ii)$, we discuss the cases $x\not=y$ and $x=y$ separately.\medskip\\
\underline{For $x\not=y$} at the point $((x, y, t), c)$ with $c=(f, s)\in \mathcal{C}_\delta\times\mathcal{N}$, we consider the map 
\begin{align*}& \Psi: (S\times S)\times (\mathcal{C}_\delta\times\mathcal{N})\to \Phi(W)\times \Phi(W),\text{ where $\Phi$ is defined in (\ref{MapPhi})},\\
& ((x,y), (f,s)) \mapsto \Big(\big( f(x,s(x)),(x,s(x)), df(x, s(x))\big), \big( f(y,s(y)), (y, s(y)), df(y, s(y))\big)\Big)
\end{align*} and first show that $d\Psi$ generates $(0,(a,w))$ for any vector $(a,w)$ in the second factor
\begin{align*} T_{(f(\text{gr}_s(y)), \text{gr}_{df}(\text{gr}_s(y)))}\big(\mspace{-2mu}(-\delta,\delta)\times T_\delta^\ast (S{\times} (-\delta, \delta)^{n-k})\mspace{-2mu}\big)=\R\times T_{\text{gr}_{df}(\text{gr}_s(y))}\big(\mspace{-2mu}T_\delta^\ast (S{\times} (-\delta, \delta)^{n-k})\mspace{-2mu}\big)\mspace{-2mu}.
\end{align*} Let $\pi$ denote the projection $T^\ast_{\delta}(S{\times} (-\delta, \delta)^{n-k})\to S\times (-\delta, \delta)^{n-k}$. Denote $$v:=d\pi (w)\in T_{(y, s(y))}(S{\times} (-\delta, \delta)^{n-k})=T_y S{\times} \mathbb{R}^{n-k}$$ and its splitting into factors by $(v_1, v_2)$. We can find $s_v\in T_{s}\mathcal{N}=(\mathcal{C}^\infty_\epsilon(S, \mathbb{R}))^{n-k}\circ\text{diag}$ which is $0$ near $x$ and is the vector $v_2-(d_y s)(v_1)\in \mathbb{R}^{n-k}$ at $y$; so for a path $\gamma_{v_1}$ in $S$ with $\gamma_{v_1}(0)=y$ and $\gamma_{v_1}'(0)=v_1$, the path $\text{gr}_{s+\tau s_v}(\gamma_{v_1}(\tau))$ generates $v=(v_1, v_2)$. For clarity, let $\mathfrak{s}_{df}$ denote the section of the bundle $T^\ast_\delta (S{\times}(-\delta,\delta)^{n-k})$ defined by the exact 1-form $df$, so $\mathfrak{s}_{df}(y, s(y))=\big((y, s(y)), df(y, s(y))\big)$. Then $$w-(d_{(y, s(y))}\mathfrak{s}_{df})(v)\in \ker d_{\mathfrak{s}_{df}(y, s(y))}\pi\cong T^\ast _{(y,s(y))}(S{\times}(-\delta,\delta)^{n-k}),$$ where the last identification uses $T_p F\cong F$ canonically for any vector space $F$.\pagebreak[1]\\ We can find $f_{a,w}\in T_f \mathcal{C}_\delta=\mathcal{C}$ such that $f_{a,w}$ is $0$ near $(x, s(x))$, and $f_{a,w}(y,s(y))=a$, $df_{a,w}(y,s(y))=w-(d_{(y, s(y))}\mathfrak{s}_{df})(v)$. Then the path $$\tau\mapsto \text{gr}_{df+\tau df_{a,w}}(\text{gr}_{s+\tau s_v}(\gamma_{v_1}(\tau)))$$ generates $w$. Moreover, for any $(a, w)$, the path $$\tau\mapsto \big((f+\tau f_{a,w})(\text{gr}_{s+\tau s_v}(\gamma_{v_1}(\tau))), \text{gr}_{df+\tau df_{a,w}}(\text{gr}_{s+\tau s_v}(\gamma_{v_1}(\tau)))\big)$$ 
\begin{align*}
\text{generates }\quad (a+(d_{(y, s(y))}f)(v), w)&\in T_{(f(\text{gr}_s(y)),\text{gr}_{df}(\text{gr}_s(y)))} \big((\text{-}\delta, \delta)\times T^\ast_\delta(S{\times} (\text{-}\delta,\delta)^{n-k})\big)\quad\\
&=\mathbb{R}\times T_{\text{gr}_{df}(\text{gr}_s(y))}\big(T^\ast_\delta (S{\times} (\text{-}\delta,\delta)^{n-k})\big).
\end{align*}
If $\gamma_0$ is a constant path at $x$, then the following path generates $0\in T_x M$: $$\tau\mapsto \big((f+\tau f_{a,w})(\text{gr}_{s+\tau s_v}(\gamma_0(\tau))), \text{gr}_{df+\tau df_{a,w}}(\text{gr}_{s+\tau s_v}(\gamma_0(\tau)))\big).$$
Therefore, we have shown that $$\big(d_{((x,y),(f,s))}\Psi\big) \big((0, v_1), (f_{a',w},s_v)\big)=(0,(a,w)),$$ where $a':=a-(d_{(y,s(y))}f)(v)$. Namely, 
\[0 \times T_{(f (\text{gr}_s(y)), \text{gr}_f(\text{gr}_s(y)))}\big(\mspace{-3mu}(-\delta, \delta){\times} T^\ast_\delta(S{\times} (-\delta, \delta)^{n-k})\mspace{-3mu}\big){\subset} \big(d_{((x,y,t),(f,s))}\Psi\big)(0{\times} TS{\times} T(\mathcal{C}_\delta{\times}\mathcal{N})).\]
$d\Phi^{-1}: T\big((-\delta, \delta)\times T^\ast_\delta(S{\times} (-\delta, \delta)^{n-k})\big) \to TW$ and $d\varphi^t: TM \to TM$ for fixed $t$ are isomorphisms. So we have shown that the image of 
\begin{align*}
d\mathring{\text{\textnormal{ev}}}_\mathcal{I}|_{TS{\times} TS{\times} 0{\times} T(\mathcal{C}{\times}\mathcal{N})}=\big(Id\times (d\varphi^t|_{TM{\times} 0})\big)\circ  (d\Phi^{-1}\times d\Phi^{-1})\circ (d\Psi|_{TS{\times} TS{\times} T(\mathcal{C}{\times}\mathcal{N})})
\end{align*} contains $(\mathring{\text{\textnormal{ev}}}_{\mathcal{I}})^\ast (0{\times} TM)$, which together with $T\Delta_M$ spans $T(M{\times} M)$, as soon as\linebreak[4] $\mathring{\text{\textnormal{ev}}}_{\mathcal{I}}(x,y,t,c)\in \Delta_M$. Note that in this case transversality is achieved even without using the $t$-direction.\medskip\\ 
\underline{For $x=y$}, $\mathring{\text{\textnormal{ev}}}_{\mathcal{I}}((x,x, t), c)=(\mathcal{I}(x,c), \varphi^t(\mathcal{I}(x,c)))$, and $d\mathring{\text{\textnormal{ev}}}_\mathcal{I}|_{\Delta_{TM}\times T((0,T])\times T(\mathcal{C}_\delta\times \mathcal{N})}$ and the diagonal $T\Delta_M$ span $T(M{\times} M)|_{\Delta_M}$, because of the $(\lambda, T)$-transversality together with the fact that $d_{(x,c)}\mathcal{I}$ is a submersion onto $\mathcal{I}^\ast (TM)$ (using a part of the argument from the previous case), as $\mathcal{I}(x,c)\in W{\subset} U(S)$. So we have established property $(ii)$.
\end{proof} 

To get a perturbation of $\iota_0$ from its thickening, we pick a generic value in the perturbation parameter space independent of $S$ and \textit{equal} for both $S$-factors in the evaluation:

\begin{corollary}\label{GENERICITY} Consider the same setting as in Lem.\ \ref{THICKENING}. Write $\mathcal{D}:=S{\times} S{\times} (0,T]$, $\mathcal{E}:=\mathcal{C}_\delta{\times} \mathcal{N}$. $Z:=(\mathring{\text{\textnormal{ev}}}_{\mathcal{I}})^{-1}(\Delta_M)$ is a submanifold of $\mathcal{D}{\times}\mathcal{E}$ and $pr_2|_{Z}: Z\subset \mathcal{D}{\times}\mathcal{E}\to \mathcal{E}$ Fredholm. Then for a regular value $c$ of $pr_2|_Z$ (existing generically by Sard-Smale) defining $\iota_c:=\mathcal{I}(\cdot, c): S\to M$, $\mathring{\text{\textnormal{ev}}}_{\iota_c}\equiv \mathring{\text{\textnormal{ev}}}_{\mathcal{I}}(\cdot, c)$ is transverse to $\Delta_M$ in $M\times M$. 
\end{corollary}
\begin{proof} This genericity trick is classical (and underlies almost all transversality arguments). To quickly see it,  $d_z\mathring{\text{\textnormal{ev}}}_{\mathcal{I}}$ at points $z\in Z$ induces an isomorphism 
\[{T_z(\mathcal{D}{\times} \mathcal{E})/T_z Z}\to T_{\hat z}(M{\times} M)/T_{\hat z}\Delta_M\]
where $\hat z:=\mathring{\text{\textnormal{ev}}}_{\mathcal{I}}(z)$. If $c=pr_2(z)$ is a regular value, then $d_z (pr_2) (T_zZ)=T_{c}\mathcal{E}$. Therefore, $d_z\mathring{\text{\textnormal{ev}}}_{\iota_c}=(d_z{\mathring{\text{\textnormal{ev}}}_{\mathcal{I}}})|_{T_z\mathcal{D}}$ induces an isomorphism $T_z\mathcal{D}/d_z (pr_1)(T_zZ)\to T_{\hat z}(M{\times} M)/T_{\hat z}\Delta M$, which implies the conclusion.\end{proof}

\begin{theorem}\label{THMNOCHORDS}
Let $S$ be a closed subcritical isotropic submanifold having a trivial neighborhood in a contact manifold $(M, \xi)$ such that $S\subset M$ satisfies the $(\lambda, T)$-transversality condition for some fixed $\lambda$, e.g.\ obtained from Ex.\ \ref{FORFREE}. Then one can find an isotropic submanifold $S_c=\iota_c(S)$ with $\iota_c$ isotropically isotopic to $\iota_0$, such that $S_c$ has no Reeb chord of $R_\lambda$ of length $\leq T$.
\end{theorem}

\begin{proof} Lem.\ \ref{THICKENING} and Cor.\ \ref{GENERICITY} imply the existence of an isotropic embedding $\iota_c: S\to M$ isotropically isotopic to $\iota_0$ such that $\iota_c$ satisfies $(CT)$. The isotropic submanifold $\iota_c(S)$ then has no Reeb chord of length $\leq T$ by Cor.\ \ref{NOCHORD}.
\end{proof}

\begin{theorem}\label{FINALTHM} Let $(T_l)\subset\mathbb{R}, l\in\mathbb{N}$ and let $\lambda_l$ be a sequence of contact forms such that $\ker\lambda_l=\xi$ and $\lambda_l=e^{f_l}\lambda_1$ and $S\subset M$ satisfies $(\lambda_l, T_l)$-transversality for all $l$ with the same $U(S)$ (e.g.\ $\lambda_l$ is non-degenerate for all $l$ as in Ex.\ \ref{FORFREE} (2) with $T=\infty$).\\
Let $\mathcal{C}_\delta{\times}\mathcal{N}$ and $\mathcal{I}$ be as in Lem.\ \ref{THICKENING}. Then there exists a comeager set $\mathcal{R}\subset \mathcal{C}_\delta{\times}\mathcal{N}$ such that for all $c\in\mathcal{R}$ holds that $\iota_c=\mathcal{I}(\cdot,c):S\rightarrow M$ is an isotropic embedding isotropic isotopic to $\iota_0$, such that $S_c:=\iota_c(S)$ avoids all Reeb chords of $R_{\lambda_l}$ of lengths $\leq T_l, \forall \,l$. 
\end{theorem}
\begin{proof} 
Consider from Lem.\ \ref{THICKENING} the construction of a thickening $\mathcal{I}:S{\times}(\mathcal{C}_\delta{\times}\mathcal{N})\rightarrow M$ of $\iota_0$ with respect to $\lambda_1$. Namely, start with the tubular neighborhood 
\[\Phi: (W, \lambda_1)\to \Big((-\delta,\delta)_z\times T^\ast_\delta\big(S{\times} (-\delta, \delta)^{n-k}\big), dz-\lambda_{\text{st}}\Big)\]
and $W\subset U(S)$ (from Defn.\ \ref{defnLT}). Note that $\mathcal{I}(\cdot, c)$ is still an isotropic embedding isotropically isotopic to $\iota_0$ in $(M, \lambda_l)$ for all $l$, as
\[\mathcal{I}(\cdot, c)^\ast\lambda_l=\mathcal{I}(\cdot, c)^\ast (e^{f_l}\lambda_1)=e^{f_l(\mathcal{I}(\cdot, c))}(\mathcal{I}(\cdot, c)^\ast \lambda_1)=0.\]
$\mathring{\text{\textnormal{ev}}}_{\mathcal{I},l}$ (where the subscript $l$ signifies the dependence on the Reeb flow $\varphi^t_{\lambda_l}$ of $R_{\lambda_l}$) is still transverse to the diagonal $\Delta_M$ in $M{\times} M$ for all $l$, as the above argument only uses $d\varphi_{\lambda_l}$ being an isomorphism. Thus the desired transversality holds verbatim as in Lem.\ \ref{THICKENING}.\\
We apply Cor.\ \ref{GENERICITY} to $pr_2|_{Z_l}: Z_l\to \mathcal{C}_\delta{\times}\mathcal{N}$ where $Z_l:=(\mathring{\text{\textnormal{ev}}}_{\mathcal{I},l})^{-1}(\Delta_M)$. Then the countable intersection $\mathcal{R}:=\bigcap_l\mathcal{R}_l$ of the sets $\mathcal{R}_l$ of regular values of $pr_2|_{Z_l}$ is comeager and non-empty in $\mathcal{C}_\delta{\times}\mathcal{N}$. Any $c\in\mathcal{R}$ gives an isotropic submanifold $S_c:=\iota_c(S)$ which avoids all Reeb chords of $R_{\lambda_l}$ of lengths $\leq T_l$ for all $l$ and is isotropically isotopic to $S$.
\end{proof}

\begin{remark} Let $S$ be a closed Legendrian submanifold ($k=n$) such that $S\subset M$ satisfies the $(\lambda, T)$-transversality with $T<\infty$, e.g.\ Ex.\ \ref{FORFREE}. Trivially it has a trivial neighborhood, as $E=0$ here. The same construction allows us to find a Legendrian embedding  ${\iota_c: S\rightarrow M}$ that satisfies (CT) and is Legendrian isotopic to $\iota_0$. Applying the neighborhood theorem to $\tilde S=\iota_c(S)$, we see that for small $0<t \leq \epsilon_0$, the Reeb flow $\varphi^t$ of $\frac{\partial\;}{\partial z}$ maps $\tilde S$ within that neighborhood and increases $z$ coordinate from $z=0$ in $(\mathbb{R}_z{\times} T^\ast\tilde S, dz\,{-}\lambda_{\text{st}})$, thus there is no Reeb chord on $\tilde S$ of length $\leq \epsilon_0$. Since $\mathring{\text{ev}}_{\iota_c}: S{\times} S{\times} (0,T]\to M{\times} M$ is transverse to $\Delta_M$, $(\mathring{\text{ev}}_{\iota_c})^{-1}(\Delta_M)$ is a $0$-dimensional manifold by the same argument as in Cor.\ \ref{NOCHORD}; points in $(\mathring{\text{ev}}_{\iota_c})^{-1}(\Delta_M)$ correspond to Reeb chords and are just shown to lie in the compact set $S{\times} S{\times} [\epsilon_0, T]$, hence form a finite set. So, one can Legendrian isotope a closed Legendrian submanifold satisfying $(\lambda, T)$-transversality with $T<\infty$ to one that admits only finitely many Reeb chords of lengths $\leq T$. This is useful for setting up Legendrian SFT.
\end{remark}

\section*{Acknowledgments}
The author wishes to thank Klaus Mohnke, Chris Wendl, Kai Cieliebak, Urs Frauenfelder, Mark, Kegel, Alexandru Oancea and Peter Uebele for fruitful discussions which contributed to this article. Special thanks go to Dingyu Yang for providing the second appendix and to the anonym reviewer for suggesting several improvements. Moreover, this paper is partly based on the authors PhD-thesis, during which he was supported by the Studienstiftung des deutschen Volkes, the SFB 647 ``Raum, Zeit, Materie'' and the Berlin Mathematical School.

\bibliographystyle{plain}
\bibliography{../References/References.bib}

\end{document}